\documentclass{LMCS}

\def\dOi{11(2:8)2015}
\lmcsheading%
{\dOi}
{1--38}
{}
{}
{Jun.~29, 2014}
{Jun.~15, 2015}
{}

\ACMCCS{[{\bf Theory of computation}]:  Computational complexity and
  cryptography---Complexity theory and logic\,/\,Proof complexity}
\subjclass{Complexity theory and logic, Proof complexity}

\usepackage{hyperref, amssymb, multirow}
\usepackage{amsthm}
\theoremstyle{plain}

\newtheorem*{theorems:exp}{Theorem \ref{theorems:exp}}

\newtheorem*{theoremss:pcptt}{Theorem \ref{pcp}}

\begin{document}




\title[Logical strength of complexity theory]{Logical strength of complexity theory and a formalization of the PCP theorem in bounded arithmetic}


\author[J.~Pich]{J\'{a}n Pich}
\address{Department of Algebra\\ Faculty of Mathematics and Physics\\ Charles University in Prague\\ Sokolovska 83, Prague, CZ-186 75, The Czech Republic}


\begin{abstract}

We present several known formalizations of theorems from computational complexity in bounded arithmetic and formalize the PCP theorem in the theory $PV_1$ (no formalization of this theorem was known). This includes a formalization of the existence and of some properties of the $(n,d,\lambda)$-graphs in $PV_1$.


\end{abstract}

\maketitle






\def\bbb{

\section{Introduction}
\label{}

We investigate the provability of polynomial circuit lower bounds in weak fragments of arithmetic like $S^1_2(bit)$ or $APC_1$. These theories are sufficiently strong to prove many important results in Complexity Theory. In fact, they can be considered as formalizations of feasible mathematics. Our motivation behind the investigation of these theories is the general question whether existential quantifiers in complexity-theoretic statements can be witnessed feasibly. 
\medskip

Intuitively, if the statement expressing $n^k$-size circuit lower bounds for SAT was such a feasibly witnessed statement, for any $n^k$-size circuit with $n$ inputs we could efficiently find a formula of size $n$ on which the circuit fails to decide SAT. We present a natural formalization of $n^k$-size circuit lower bounds for SAT denoted $LB(SAT,n^k)$ and observe that its provability in $S^1_2(bit)$ gives us such error witnessing. One could hope to use the witnessing algorithm to derive a contradiction with some established hardness assumption, however, Atserias and Kraj\'{i}\v{c}ek (private communication) noticed that certain cryptographic conjectures imply the same form of witnessing, see Proposition \ref{AK}.
\medskip

We do not know how to obtain the unprovability of SAT circuit lower bounds in $S^1_2(bit)$ but we can do it basically for any weaker theory with stronger witnessing properties.
\medskip

In weaker theories the situation is less natural because they cannot fully reason about p-time concepts. In particular, $LB(SAT,n^k)$ is equivalent to a formula $LB_2(SAT,n^k)$ (defined in Section 5) in $S^1_2(bit)$ but not necessarily in weaker theories. Therefore, we need to consider these two formalizations separately. We present it in the case of theory $T_{NC^1}$ which is the true universal first-order theory in the language containing names for all uniform $NC^1$ algorithms.
\smallskip

If $T_{NC^1}$ proves $LB_2(SAT,n^k)$, there are uniform $NC^1$ circuits which for each $n^k$-size circuit $C$ with large enough $n$ find a formula $y$ of size $n$ and computation of $C$ on $y$ witnessing that $C$ decides SAT incorrectly on $y$. It is easy to show that in such case nonuniform $NC^1$ circuits could simulate $SIZE(n^k)$, see Proposition \ref{lb2}. Thus, a conditional unprovability of $LB_2(SAT,n^k)$ in $T_{NC^1}$ follows easily.
\smallskip

To prove $LB(SAT,n^k)$ in $T_{NC^1}$, the resulting uniform $NC^1$ circuits would need to output for each $n^k$-size circuit $C$ with large enough $n$ an error $y$ but they would not need to witness the computation of $C$ on $y$. In this sense, for $T_{NC^1}$ it is easier to reason about formalization $LB(SAT,n^k)$. We show that even $LB(SAT,n^{2kc})$ for $k\geq 1, c\geq 4$ is unprovable in $T_{NC^1}$ unless each $f\in SIZE(n^k)$ can be approximated by formulas $F_n$ of size $2^{O(n^{2/c})}$ with subexponential advantage: $P_{x\{0,1\}^n}[F_n(x)=f(x)]\geq 1/2+1/2^{O(n^{2/c})}$. The proof will be quite generic so, in particular, using known lower bounds on PARITY function, we will obtain that, unconditionally, $V^0$ cannot prove quasi polynomial ($n^{\log n}$-size) circuit lower bounds on SAT. Here, $V^0$ is a second-order theory used frequently in Bounded Arithmetic, see Section 5.
\medskip

To prove our main claim we firstly observe that by the KPT theorem \cite{k4} the provability of $LB(SAT,n^k)$ in universal theories like $T_{NC^1}$ gives us an $O(1)$-round Student-Teacher (S-T) protocol finding errors of $n^{2kc}$-size circuits attempting to compute SAT. Then, in particular, it works for $n^{2kc}$-size circuits encoding Nisan-Wigderson (NW) generators based on any functions $f\in SIZE(n^k)$ and suitable design matrices \cite{nis}. The interpretation of NW-generators as p-size circuits comes from Razborov \cite{r2}. In this situation we apply Kraj\'{i}\v{c}ek's proof that certain NW-generators are hard for $T_{PV}$ \cite{k3} which is the main technique we use. We show that it works in our context as well and allows us to use the S-T protocol to compute $f$ by subexponential formulas with a subexponential advantage.   
\smallskip

Perhaps the most significant earlier result of this kind was obtained by Razborov \cite{r}. Using natural proofs he showed that theory $S^2_2(\alpha)$ cannot prove polynomial circuit lower bounds on SAT unless strong pseudorandom generators do not exist. The second-order theory $S^2_2(\alpha)$ is however quite weak with respect to the formalization Razborov used. As far as we know his technique does not imply the unprovability of circuit lower bounds (formalized as here, see Section 2) even for Robinson's Arithmetic Q. In this respect, our proof applies to much stronger theories, basically to any theory weaker than $S^1_2(bit)$.
\smallskip

The paper is organized as follows. In Section 2 we formalize circuit lower bounds in the language of bounded arithmetic. In Section 3 we define theory $S^1_2(bit)$, state its properties and in Section 4 discuss the provability of circuit lower bounds in $S^1_2(bit)$. Section 5 defines subtheories of $S^1_2(bit)$ for which we prove our main unprovability results in Section 6.

\section{Formalization}

The usual language of arithmetic contains well known symbols: $0, S, +, \cdot, =, \leq$. To encode reasoning about computation it is natural to consider also symbols $\lfloor \frac{x}{2} \rfloor$, $|x|$ for the length of binary representation of $x$ and $\#$ with the intended meaning $x\# y=2^{|x|\cdot |y|}$. Theories of bounded arithmetic are defined using language $L=\{0, S, +, \cdot, =, \leq, \lfloor x/2 \rfloor, |x|, \#\}$. We will consider also language $L_{bit}$ which contains in addition symbol $x_i$ for the $i$-th bit of the binary representation of $x$. The basic properties of symbols from $L_{bit}$ are captured by a set of basic axioms $BASIC(bit)$ which we will not spell out, cf. \cite{b,k}.
\medskip

$\Sigma^b_0$ denotes the set of all formulas in the language $L$ with all quantifiers sharply bounded: $\exists x, x\leq |t|$, $\forall x, x\leq |t|$ where $t$ is a term not containing $x$. All relations defined by $\Sigma^b_0$ formulas are p-time computable. $\Sigma^b_i$ resp. $\Pi^b_i$ for $i>0$ are sets of formulas constructed from sharply bounded formulas by means of $\wedge$, $\vee$, sharply bounded, and existential bounded quantifiers: $\exists y\ y\leq t$ resp. universal bounded quantifiers: $\forall y\ y\leq t$ for $x$ not occurring in $t$. All NP resp. coNP are representable by $\Sigma^b_1$ resp. $\Pi^b_1$ formulas, cf. \cite{l1,l2,l3}.
\medskip

Define $\Sigma^b_i(bit), \Pi^b_i(bit)$ for $i\geq 0$ as above but in the language $L_{bit}$. For $i\geq 1$, $\Sigma^b_i(bit)$-formulas are actually equivalent to $\Sigma^b_i$-formulas in theory $PV_1$, cf. \cite{c3,k}, see also Section 3. Analogously, for $\Pi^b_i$-formulas with $i\geq 1$.
\smallskip

We will now express circuit lower bounds in $L_{bit}$.
\medskip

Firstly, denote by $Comp(C,y,w)$ a $\Sigma^b_0(bit)$-formula saying  that $w$ is a computation of circuit $C$ on input $y$. Such a formula can be constructed in many ways and our results work for any $\Sigma^b_0(bit)$ formalization. For simplicity, we present here a less efficient one where $C$ represents a directed graph on $|w|$ vertices.       
\smallskip

Let $E_C(i,j)$ be $C_{[i,j]}$ for pairing function $[i,j]=(i+j)(i+j+1)/2+i$. $E_C(i,j)=1$, $i,j<|w|$ means that there is an edge in circuit $C$ going from the $i$-th vertex to the $j$-th vertex. For $k<|w|$, let $N_C(k)$ be the tuple of bits $(C_{[|w|,|w|]+2k},C_{[|w|,|w|]+2k+1})$ encoding the connective in the $k$-th node of circuit $C$, say $(0,1)$ be $\wedge$, $(1,0)$ be $\vee$, and $(1,1)$ and $(0,0)$ be $\neg$. Therefore, $|C|=[|w|,|w|]+2|w|$. Then let $Circ(C,y,w)$ be the formula stating that $C$ encodes a $|w|$-size circuit with $|y|$ inputs:\\ 

$\forall j<|w|, j\geq |y|$

\ \ \ $(N_C(j)=(1,0)\vee N_C(j)=(0,1)\rightarrow \exists i,k<j\ i\neq k\forall l<j, l\neq k, l\neq j$
 
\ \ \ \ \ \ $(E_C(i,j)=1\wedge E_C(k,j)=1\wedge E_C(l,j)=0) )\wedge$ 

\ \ \ $(N_C(j)=(1,1)\vee N_C(j)=(0,0)\rightarrow \exists i<j\forall l<j, k\neq i$

\ \ \ \ \ \ $(E_C(i,j)=1\wedge E_C(l,j)=0))$\\ 

which means that if the $j$-th node of $C$ is $\wedge$ or $\vee$, there are exactly two previous nodes $i,k$ of $C$ with edges going from $i$ and $k$ to $j$, if the $j$-th node of $C$ is $\neg$, there is exactly one previous node $i$ with an edge going from $i$ to $j$.

$Comp(C,y,w)$ says that for each $i<|y|$ the value of $w_i$ is the value of the $i$-th input bit of $y$ and each $w_j$ is an evaluation of the $j$-th node of circuit $C$ given $w_k$'s evaluating nodes connected to the $j$-th node:\\

$Circ(C,y,w)\wedge \forall i<|y|\ y_i=w_i\wedge \forall j,k,l<|w|[\\ (N_C(j)=(1,0)\wedge E_C(k,j)=1\wedge E_C(l,j)=1\rightarrow (w_j=1\leftrightarrow w_k=1\wedge w_l=1))\wedge\\ (N_C(j)=(0,1)\wedge E_C(k,j)=1\wedge E_C(l,j)=1\rightarrow (w_j=1\leftrightarrow w_k=1\vee w_l=1))\wedge\\ ((N_C(j)=(0,0)\vee N_C(j)=(1,1))\wedge E_C(k,j)=1\rightarrow (w_j=1\leftrightarrow w_k=0)) ]$\\

Formula $C(y;w)=1$ stating that $w$ is accepting computation of circuit $C$ on input $y$ will be $Comp(C,y,w)\wedge w_{|w|-1}=1$. Similarly for $C(y;w)=0$.
\medskip

Next, let $SAT(y,z)$ be a $\Sigma^b_0(bit)$-formula saying that $z$ is a satisfying assignment to the propositional 3-CNF formula $y$. 
\smallskip

To define it explicitly for each $i,j,k<2m$ we let $y_{[i,j,k]}=1$ if and only if the 3-CNF encoded in $y$ contains a clause of variables $v^p_i,v^p_j, v^p_k$ where $v^p_i$ is $v_i$ if $i<m$ and $\neg v_{i-m}$ if $i\geq m$. Here also $[i,j,k]=[i,[j,k]]$. Hence, the 3-CNF encoded in $y$ has $m$ variables $v_0,...,v_{m-1}$ and $|y|=[2m-1,2m-1,2m-1]+1$. We use $m$ implicitly given by $y$ in the formula $SAT(y,z)$:\\ 

$\forall i,j,k<2m \ [y_{i,j,k}=1\rightarrow$ 

\ \ \ \ \ \ $(i,j,k<m\rightarrow z_i=1\vee z_j=1\vee z_k=1)\wedge$

\ \ \ \ \ \ $(i,j<m\wedge k\geq m\rightarrow z_i=1\vee z_j=1\vee z_{k-m}=0)\wedge$

\ \ \ \ \ \ $ ... $

\ \ \ \ \ \ $ (i,j,k\geq m\rightarrow z_{i-m}=0\vee z_{j-m}=0\vee z_{k-m}=0)]$\\

Finally, for any $k$, a hardness of SAT for $n^k$-size circuits can be expressed as

$LB(SAT,n^k):$
\smallskip

$\forall 1^n>n_0\ \forall C\ \exists y,a\ |a|<|y|=n\ \forall w,z\ |w|\leq n^k, |z|<|y|\ [Comp(C,y,w)\rightarrow$

$(C(y;w)=1\wedge \neg SAT(y,z))\vee (C(y;w)=0\wedge SAT(y,a))]$\\

Here $n_0$ is a fixed constant which is not indicated in $LB(SAT,n^k)$. This should not cause any confusion. Whenever we say that $LB(SAT,n^k)$ is provable in a theory $T$ we mean that it is provable in $T$ for some $n_0$. Further, $\forall 1^n>n_0$ is a shortcut for $\forall m,n$ such that $|m|=n\wedge m>n_0$. Therefore, $y$ is feasible in $m$ and for each $n_0$ and $k$, $LB(SAT,n^k)$ is universal closure of a $\Sigma^b_2(bit)$ formula.    

\section{Feasible Mathematics}

If we obtain $n^k$-size circuit lower bounds for SAT but do not find any efficient method how to witness errors of potential $n^k$-size circuits for SAT, some of these circuits might work in practice like correct ones. We will now define theories of feasible mathematics where provability of $n^k$-size circuit lower bound for SAT implies the existence of such an error witnessing.
\medskip

Perhaps, the most prominent one is $S^1_2$ introduced by Buss \cite{b}. We will use its conservative extension $S^1_2(bit)$ which consists of $BASIC(bit)$ and polynomial induction for $\Sigma^b_1(bit)$-formulas $A$: $$A(0)\wedge \forall x (A(\lfloor x/2\rfloor)\rightarrow A(x))\rightarrow \forall x A(x)$$ An important property of $S^1_2(bit)$ is Buss's witnessing theorem:

\begin{thm}[Buss \cite{b}] If $S^1_2(bit)\vdash \exists y A(x,y)$ for $\Sigma^b_0(bit)$-formula $A$, then there is a p-time functions $f$ such that $A(x,f(x))$ holds for any $x$.
\end{thm}

$S^1_2(bit)$ admits also a useful kind of witnessing for $\Sigma^b_2(bit)$-formulas.

\begin{thm}[Kraj\'{i}\v{c}ek \cite{k5}]\label{skpt} If $S^1_2(bit)\vdash \exists y\forall z\leq t\ A(x,y,z)$ for $\Sigma^b_0(bit)$-formula $A$ and term $t$ depending only on $x,y$, then there is p-time algorithm S such that for any $x$ either $\forall z\leq t\ A(x,S(x),z)$ or for some $z_1$, $\neg A(x,S(x),z_1)$. In the latter case, either $\forall z\leq t\ A(x,S(x,z_1),z)$ or there is $z_2$ such that $\neg A(x, S(x,z_1), z_2)$. However after $k\leq poly(|x|)$ rounds of this kind, $\forall z\leq t\ A(x,S(x,z_1,...,z_k),z)$ holds for any $x$.
\end{thm}

Another theory with similar witnessing properties is $PV_1$ which is an extension of a theory $PV$ defined by Cook \cite{c3}, see also \cite{k}. The language of $PV_1$ consist of symbols for all functions given by a Cobham-like inductive definition of p-time functions (hence it contains $L_{bit}$). $PV_1$ defined in \cite{k4} is then a first-order theory axiomatized by equations defining all the function symbols and a derivation rule similar to polynomial induction for open formulas. It is a universal theory, i.e. it has an axiomatization by purely universal sentences, and because all function symbols of $PV_1$ have well-behaved $\Sigma^b_1$ and $\Pi^b_1$ definitions in $S^1_2(bit)$, $PV_1$ is contained in the extension of $S^1_2(bit)$ by these definitions which we denote also $S^1_2(bit)$. $PV_1$ proves induction and polynomial induction for $\Sigma^b_0(PV)$-formulas defined similarly as $\Sigma^b_0$-formulas but in the language of $PV_1$. There is also an interesting theory $APC_1$ introduced by Je\v{r}\'{a}bek \cite{e} which is an extension of $S^1_2(bit)$ capturing a subclass of BPP similarly as $S^1_2(bit)$ captures P.
\medskip

Theories $S^1_2(bit), PV_1$ and $APC_1$ are weak fragments of arithmetic but they are sufficiently strong to prove many important things. In \cite[chap. 15]{k} it is shown how to prove PARITY $\notin AC^0$ in $APC_1$. Razborov\cite{r3} argued that $S^1_2(bit)$ is the right theory capturing techniques from circuit complexity in 1995. We expect that $APC_1$ captures feasible reasoning so well that any provable statement about feasible concepts is provable in $APC_1$ assuming that feasible concepts intuitively correspond to BPP concepts . Of course, this does not contain, for instance, Shannon's argument if it is formalized so that it manipulates with exponentially big objects. 

\subsection{Equivalent formalizations of $LB(SAT,n^k)$}

There are more possible formalizations of circuit lower bounds that are essentially equivalent to $LB(SAT,n^k)$. For example, $SCE(SAT,n^k)$ meaning that for each $n^k$-size circuit there is a satisfiable formula of size $n$  such that the circuit will not find its satisfying assignment.\\

$SCE(SAT,n^k):$
\smallskip

$\forall 1^n>n_0\ \forall C\ \exists y,a |a|<|y|=n\ \forall w,z\ |w|\leq n^k, |z|<|y|$

\ \ \ $[SAT(y,a)\wedge (C(y;w)=z\rightarrow \neg SAT(y,z))]$\\

where $C(y;w)=z$ means that $w$ is a computation of circuit $C$ on input $y$ with output bits $z$. Formally, $Comp(C,y,w)\wedge\forall i<|z| (w_{|w|-i-1}=1\leftrightarrow z_i=1)$. SCE in $SCE(SAT,n^k)$ refers to ``search SAT counter example".
\medskip

Another formalization of circuit lower bounds is given by the following formula $DCE(SAT,n^k)$ where DCE refers to ``decision SAT counter example". Now circuits $C$ have again just one output but using self-reducibility they can be used to search for satisfying assignments of propositional formulas: If $C$ says that formula $y$ is satisfiable, we can set the first free variable in $y$ firstly to 1 and then to 0, and use $C$ to decide in which of these cases the resulting formula is satisfiable, then in the same manner continue searching for the full satisfying assignment. If no such $C$ can be used to find satisfying assignments of satisfiable propositional formulas, for each such $C$ there is a formula $y$ and a possibly partial assignment to its variables $a$ such that either $SAT(y,a)$ and $C$ says that $y$ is unsatisfiable or $\neg SAT(y,a)$ for full assignment $a$ of $y$ and $C$ says that $a$ satisfies $y$ or it happens that $C$ gets into a local inconsistency: for a partial assignment $a$ of $y$ it claims that $y$ assigned by $a$ is satisfiable but when we extend $a$ by setting the first of the remaining free variables by 1 and 0 in both cases $C$ claims that the resulting formula is unsatisfiable. Formally:\\

$DCE(SAT,n^k):$
\smallskip

$\forall 1^n>n_0\forall C\ \exists y,a |a|<|y|=n\ \forall w^0,...,w^4\ |w^0|,...,|w^4|\leq n^k[$ 

\ \ \ $(Comp(C,y,w^0)\rightarrow (C(y;w^0)=0\wedge SAT(y,a)))\vee$ 

\ \ \ $(Comp(C,y(a),w^1)\rightarrow (C(y(a);w^1)=1\wedge FA(a,y)\wedge\neg SAT(y,a)))\vee$ 

\ \ \ $(Comp(C,y(a),w^2)\rightarrow (C(y(a);w^2)=1\wedge PA(a,y)\wedge$ 

\ \ \ \ \ \ \qquad\qquad\qquad\qquad\qquad $(Comp(C,y(a1),w^3)\rightarrow C(y(a1);w^3)=0)\wedge$

\ \ \ \ \ \ \qquad\qquad\qquad\qquad\qquad $(Comp(C,y(a0),w^4)\rightarrow C(y(a0);w^4)=0)))] $\\

where $y(a)$ encodes formula $y$ assigned by $a$, $FA(a,y)$ resp. $PA(a,y)$ means that $a$ is full resp. partial assignment to variables in $y$ and $y(a1)$ resp. $y(a0)$ is $y$ assigned by extension of $a$ which set the first unassigned variable in $y$ by 1 resp. by 0. We leave details of these encodings to the reader.
\smallskip

$LB(SAT,n^k), SCE(SAT,n^k), DCE(SAT,n^k)$ are basically equivalent. We claim that this is provable already in $PV_1$ and hence also in $S^1_2(bit)$. 

\begin{prop} $PV_1$ proves the following implications

\centerline{$SCE(SAT,n^{2k})\rightarrow LB(SAT,n^k)$}
\centerline{$LB(SAT,n^{2k})\rightarrow SCE(SAT,n^k)$}
\centerline{$LB(SAT,n^{k})\rightarrow DCE(SAT,n^k)$}
\centerline{$DCE(SAT,n^{k})\rightarrow LB(SAT,n^k)$}
where $n_0$ is arbitrary but the same constant in the assumption and the conclusion of each implication.
\end{prop}

Proof: The first implication was observed in \cite{ck}: Assume $\neg LB(SAT,n^k)$, i.e. for a big enough $n$ there is an $n^k$-size circuit $C$ deciding SAT on instances of size $n$. Then there is a p-time function which given a circuit $C$ witnessing $\neg LB(SAT,n^k)$ produces an $n^{2k}$-size circuit $sC$ which outputs a satisfying assignment $sC(y)$ for every satisfiable formula $y$ of size $n$. For each $i$, the circuit $sC$ finds the $i$-th bit of the satisfying assignment by asking $C$ whether $y$ remains satisfiable if the $i$-th variable is set to 1, given the values it has previously found for the first $i-1$ variables. Then (assuming $\neg LB(SAT,n^k)$ and $SAT(y,a)$) $PV_1$ proves by $\Sigma^b_0(PV)$ induction on $i$ that $y$ instantiated by the first $i$ truth values is satisfiable according to $C$ and hence $\neg SCE(SAT,n^{2k})$.  
\smallskip

Concerning the second implication: If $\neg SCE(SAT,n^k)$, i.e. for a big enough $n$ there is an $n^k$-size circuit $C$ which outputs a satisfying assignment $C(y)$ for every satisfiable formula of size $n$, then there is a p-time function which given any such circuit $C$ produces an $n^{2k}$-size circuit $dC$ which decides SAT on instances of size $n$. Given a formula $y$, $dC$ outputs 1 if and only if $C(y)$ satisfies $y$. Assuming $\neg SCE(SAT,n^k)$ it follows in $PV_1$ that $(SAT(y,a)\rightarrow dC(y;w)=1)\wedge (dC(y;w)=1\rightarrow SAT(y,C(y)))$ for any $y,a$ of size $|a|<|y|=n$, hence $\neg LB(SAT,n^{2k})$. 
\smallskip 

Next, in $PV_1$, if circuit $C$ witnesses $\neg DCE(SAT,n^k)$, it witnesses also $\neg LB(SAT,n^k)$: for any $y,a$ of size $|a|<|y|=n$ for a big enough $n$, $C(y;w)=0\rightarrow \neg SAT(y,a)$ and if $C(y;w)=1$ then 
 by $\Sigma^b_0(PV)$-induction (as in the first implication) $C(y(b);w)=1$ for a full assignment $b$ of $y$ for which $SAT(y,b)$ holds.  
\smallskip

Finally, in $PV_1$, if $C$ witnesses $\neg LB(SAT,n^k)$, it witnesses $\neg DCE(SAT,n^k)$: for any $y,a$ of size $|a|<|y|=n$ for a big enough $n$, $(C(y;w)=0\rightarrow \neg SAT(y,a))$, $C(y(a);w)=1\wedge FA(a,y)\rightarrow SAT(y,a)$ and if $C(y(a);w)=1\wedge PA(a,y)$ then for some $b$ extending $a$ $SAT(y,b)$ and thus $C(y(a1);w)=1\vee C(y(a0);w)=1$.    \qed

\subsection{Witnessing errors of p-size circuits}

Using $LB(SAT,n^k), SCE(SAT,n^k)$ and $DCE(SAT,n^k)$ we can define several types of error witnessing of p-size circuits claiming to solve SAT. 
\smallskip

We say somewhat informally that $LB(SAT,n^k)\in$ P if there is a p-time algorithm $A$ which for any sufficiently large $n$ and $n^k$-size circuit $C$ with $n$ inputs finds out $y,a$ such that $LB(C,y,a)$: $C(y)=1\wedge SAT(y,a)$ or $C(y)=0\wedge \forall z \neg SAT(y,z)$. Intuitively, $A$ witnesses the important existential quantifiers in $LB(SAT,n^k)$.
\medskip

We say that $LB(SAT,n^k)$ has an S-T protocol with $l$ rounds if there is a p-time algorithm S such that for any function T and any sufficiently large $n$, whenever S is given $n^k$-size circuit $C$, S outputs $y_1,a_1$ such that either $LB(C,y_1,a_1)$ or otherwise T sends to S $w_1,z_1$ certifying $\neg LB(C,y_1,a_1)$. Then S uses $C, w_1,z_1$ to produce $y_2,a_2$ and the protocol continues in the same way, S possibly using all counter-examples T sent in earlier rounds. But after at most $l$ rounds S outputs $y,a$ such that $LB(C,y,a)$.
\medskip

Analogously, $DCE(SAT,n^k)\in$ P if there is a p-time algorithm $A$ which for any $n^k$-size circuit $C$ with $n$ inputs finds out $y,a$ such that $DCE(C,y,a)$: 

$C(y)=0\wedge SAT(y,a)$ or $C(y(a))=1\wedge FA(a,y)\wedge \neg SAT(y,a)$ or

$C(y(a))=1\wedge PA(a,y)\wedge (C(y(a0))=0\vee C(y(a1))=0)$
\smallskip

$SCE(SAT,n^k)\in$ P if there is a p-time algorithm $A$ which for any $n^k$-size circuit $C$ with $n$ inputs and $n$ outputs finds out $y,a$ such that $SAT(y,a)\wedge\neg SAT(y,C(y))$. 
\smallskip

The phrase that $DCE(SAT,n^k)$ resp. $SCE(SAT,n^k)$ has an S-T protocol with $l$ rounds could be defined similarly but notice that in this case T's advise would consist only of computations $w$ of given circuit $C$ which can be produced by S itself as it has $C$ as input.
\smallskip

In practice, if we want to witness that no small circuit solves SAT, it does not seem sufficient to have a p-time algorithm for $LB(SAT,n^k)$ because such an algorithm could output a tautology but we would not have an apriori way to certify that it is indeed a tautology and hence a correctly witnessed error. Therefore, it seems that practically more appropriate error witnessing is defined by $DCE(SAT,n^k)$ or $SCE(SAT,n^k)$ in which we actually force given circuits to claim inconsistent statements. We discuss it in more details in the next section.

\section{Circuit Lower Bounds in $S^1_2(bit)$}

The provability of circuit lower bounds in $S^1_2(bit)$ gives us an efficient witnessing errors of p-size circuits for SAT described in the previous section. 

\begin{prop} If $S^1_2(bit)\vdash LB(SAT,n^k)$, then $LB(SAT,n^k)$ has an S-T protocol with $poly(n)$ rounds. If $S^1_2(bit)\vdash SCE(SAT,n^k)$, then $SCE(SAT,n^k)\in P$. If $S^1_2(bit)\vdash DCE(SAT,n^k)$, then $DCE(SAT,n^k)\in P$.
\end{prop}

Proof: $LB(SAT,n^k), DCE(SAT,n^k)$ and $SCE(SAT,n^k)$ are universal closures of $\Sigma^b_2(bit)$-formulas so the first implication follows directly from Kraj\'{i}\v{c}ek's witnessing theorem. In case of $SCE(SAT,n^k)$ and $DCE(SAT,n^k)$ T's advise in the resulting S-T protocol consist just of computations of given circuit $C$. This can be, however, produced by S itself as it has $C$ as input. \qed
\smallskip

An efficient witnessing errors of p-time SAT algorithms follows also from instance checkers for SAT, see \cite[chap. 8]{ba}. If we want to check only $n^k$-time algorithms, the instance checker is p-time itself:

\begin{thm} There is a p-time algorithm that given any $n^k$-time algorithm $M$ and a formula $y$ of size $n$ accepts if $M$ solves SAT on all instances, and rejects with probability $\geq 1-1/2^n$ if $M$ does not decide satisfiability of $y$ correctly. 
\end{thm}

Therefore, any $n^k$-time algorithm $M$ claiming to solve SAT can be tested by checking it on formula $F_M(y,a)$ encoding the statement ``$a$ satisfies formula $y$ but $M$ fails to find a satisfying assignment of $y$ (in the same way as $C$ fails to find it in $DCE(SAT,n^k)$)". If $M(F_m(y,a))=1$, by self-reducibility $M$ will be forced to find a satisfying assignment of $F_M$ which is an error of $M$ or it will end up in a local inconsistency which is also error. If $F_M(y,a)$ is unsatisfiable, the checker will use an interactive protocol with $M$ as a Prover to verify that.

In practice, we can test whether a given algorithm $M$ proves theorems efficiently also by taking a statement we consider hard to prove and hard to refute instead of $F_M$ . 
\medskip

Furthermore, if $f$ is one-way function, we can also secretly produce $a\in\{0,1\}^n$ and ask the algorithm whether the statement $f(a)=f(x)$ encoded as a $poly(|a|)$-size formula with free variables $x=x_1,...,x_n$ is satisfiable, see \cite{c2}. In this case, we do not need to use interactive protocols because the algorithm is forced to say that the formula is satisfiable and by the choice of $f$, with high probability it will not find its satisfying assignment. Atserias (private communication) suggested to derandomize this construction and Kraj\'{i}\v{c}ek made the following observation.

\begin{prop}\label{AK} If there exists one-way permutation $f$ secure against p-size circuits, i.e. for any p-size circuits $C_n$ there is a function $\epsilon(n)=n^{-\omega(1)}$ such that for large enough $n$, $$P_{x\in\{0,1\}^n}[C_n(f(x))=x]\leq \epsilon(n)$$ and if there exists $h\in E$ hard on average for subexponential circuits, i.e. there is $\delta>0$ such that for all circuits $C_n$ of size $\leq 2^{\delta n}$ and large enough $n$, $$P_{x\in\{0,1\}^n}[C_n(x)=h(x)]\leq 1/2+1/2^{\delta n}$$ then for each $k$, $SCE(SAT,n^k)\in P$.
\end{prop}

Proof: If there is $h\in E$ hard on average for subexponential circuits, by \cite{nis} for each $l$ there is $c$ and NW-generator $g:\{0,1\}^{c\log n}\mapsto \{0,1\}^n$ such that $g$ is $poly(n)$-time computable and for any $n^l$-size circuits $D_n$, $$|P_{x\in\{0,1\}^{c\log n}}[D_n(g(x))=1]-P_{x\in\{0,1\}^n}[D_n(x)=1]|\leq 1/n$$ This generator allows us to derandomize the construction above: Let $f$ be one-way permutation secure against p-size circuits. Take $l$ such that for each $n^k$-size circuits $C_n$ predicate $C_n(f(x))\neq x$ for $x\in\{0,1\}^n$ can be computed by $n^l$-size circuits. Now, for any $n^k$-size circuit $C_n$ with sufficiently big $n$, for each $x\in \{0,1\}^{c\log n}$ find out whether $C_n(f(g(x)))\neq g(x)$ holds. This can be done in $poly(n)$-time. If we did not succeed at least once, $P_{x\in\{0,1\}^{c\log n}}[C_n(f(g(x)))=g(x)]=1$, and that would break $g$. \qed
\smallskip

In \cite{k2} Kraj\'{i}\v{c}ek also observed that in order to show $SCE(SAT,n^k)\in P/poly$, it suffices to assume that SAT $\notin SIZE(n^{2k})$. It uses a well known combinatorial principle: Let $E\subseteq X\times Y$ be a bipartite graph, $|X|=2^{n^k}, |Y|=2^n$. Then $$\forall x_1,...,x_n\in X\ \exists y\in Y\ \bigwedge_{i=1,...,n} E(x_i,y)\Rightarrow$$ $$\exists y_1,...,y_{n^k}\in Y\ \forall x\in X \ \bigvee_{i=1,...,n^k}E(x,y_i)$$ Now take as $X$ the set of all $n^{k/2}$-size circuits and interpret $E(x,y)$ as ``if $y$ is a satisfiable formula of size $n$, circuit $x$ does not find a satisfying assignment of $y$". If SAT restricted to instances of size $n$ does not have $n^k$-size circuits then for every $n$ circuits $C_1,...,C_n$ of size $n^{k/2}$ there is $y$ such that $\bigwedge_{i=1,...,n} E(C_i,y)$. Else, for any satisfiable $y$ at least one of the $n$ fixed circuits would find a satisfying assignment of $y$. By the principle above, there are then $y_1,...,y_{n^k}$ such that for each $n^{k/2}$-size circuit $C$, $\bigvee_{i=1,...,n^k} E(C,y_i)$. Therefore there is an $n^{2k}$-size circuit which for each $x\in X$ finds $y$ such that $E(x,y)$ by trying $E(x,y_i)$ for $i=1,...,n^k$ and thus using additional satisfying assignments $a_1,...,a_{n^k}$ of respective $y$'s as advice solves $SCE(SAT,n^k)$.       
\smallskip

Similarly, it works for $DCE(SAT,n^k)$ because checking $E(x,y)$, i.e. whether circuit $x$ (with one output) can be used to find the satisfying assignment, is efficient. For $LB(SAT,n^k)$ it could however happen that the search for the satisfying assignment ends in a local inconsistency. 

\begin{prop}[Kraj\'{i}\v{c}ek \cite{k2}] If SAT $\notin SIZE(n^{2k})$, then $SCE(SAT,n^k)$ and $DCE(SAT,n^k)$ are in P/poly.
\end{prop}

Proposition \ref{AK} seems to imply that for proving $S^1_2(bit)\not\vdash SCE(SAT,n^k)$ we need to use other properties than $SCE(SAT,n^k)\in$ P. Moreover, assumptions of Proposition \ref{AK} give us an S-T protocol for $LB(SAT,n^k)$ too.  Informally, any $n^k$-size circuit $C$ claiming to decide SAT can be used to search for satisfying assignments of propositional formulas. Using the algorithm from Proposition \ref{AK}, S can produce $y,a$, such that $SAT(y,a)$ but $C$ will not find any satisfying assignment of $y$. This means that either $C$ claims that $y$ is unsatisfiable or the assignment it finds does not satisfy $y$ or while searching for a satisfying assignment it gets into a local inconsistency which is the only case when S needs to ask for an advice of T, a satisfying assignment of $y$ extending the partial assignment found by $C$. 

\begin{prop}\label{ll} If the same hardness assumption as in Proposition \ref{AK} holds, then $LB(SAT,n^k)$ has an S-T protocol with $poly(n)$ rounds where S is in uniform $AC^0$, and it has also an S-T protocol with 1 round (i.e. 1 advice of T) where S is a p-time algorithm.
\end{prop}

Proof: 

By Proposition \ref{AK} we have a p-time algorithm $A$ solving $SCE(SAT,n^{2k})$. Firstly, we show that $LB(SAT,n^k)$ has an S-T protocol with 1 round and p-time S.
\smallskip
   
For each $n^k$-size circuit $C$ with one output bit, there is a circuit $sC$ of size $\leq n^{2k}$ searching for satisfying assignments of given formulas: For each formula $y$, let $a$ be a partial assignment of $y$ produced by $sC$ so far (empty at the beginning) and denote by $y(a)$ the formula $y$ assigned by $a$. If $C(y(a))=0$, $sC$ outputs an assignment of $y$ full of zeros. If $C(y(a))=1$, it assigns $y_a^1$, the first free variable in $y(a)$, firstly by 1 and then by 0. Denote the resulting formula $y(a1)$ resp. $y(a0)$. If $C(y(a1))=C(y(a0))=1$, $sC$ sets $y_a^1=1$.  If $C(y(a1))=C(y(a0))=0$, $sC$ outputs an assignment of $y$ full of zeros. If $C(y(a1))=1$ and $C(y(a0))=0$, $sC$ sets $y_a^1=1$. If $C(y(a1))=0$ and $C(y(a0))=1$, it sets $y_a^1=0$. In this way $sC$ sets all variables in $y$.                   
\medskip

Given $C$, S can produce $sC$ in p-time and use $A$ to find $y,a_1$ such that $SAT(y,a_1)$ but $\neg SAT(y, sC(y))$.
\smallskip

If $C(y)=0$, S outputs $y,a_1$. Else, S simulates $sC$. If it never happens that $C(y(a1))=C(y(a0))=0$ for any partial assignment $a$ produced by $sC$, S outputs $y, sC(y)$. Otherwise, for some partial assignment $a$ of $y$, $C(y(a))=1$ and $C(y(a0))=C(y(a0))=0$. In such case S outputs $y,a_2$ where $a_2$ is a full assignment of $y$ extending $a$ with all zeros. If this is not a correct answer, T replies with $a_3$ extending $a$ and satisfying $y$. Then S outputs $y(ab),a_3$ where $b\in\{0,1\}$ such that $ab$ is consistent with $a_3$. 
\smallskip

In all cases S succeeds after asking for at most 1 advice of T. 
\medskip

To get S in uniform $AC^0$ note that $A$ actually produces a set $B$ of $\leq n^c$ elements such that each $n^{2k}$-size circuit fails on at least one of them. It suffices to use instead of $A$ the set $B$, i.e. to try all elements from $B$ in place of $a_1$. Moreover, whenever S needs to simulate circuit $C$ on input $y$ it can output $y$ with an arbitrary assignment $c$ of $y$. If this is not a correct answer, T will reply either with a satisfying assignment $d$ of $y$ or with the computation of $C$ on $y$ which can be verified by a uniform constant-depth formula. In the former case S outputs $y,d$ and this time it gets what it wants.      \qed 

\smallskip

This also shows that if $SCE(SAT,n^k)\in$ P, then $DCE(SAT,n^k)\in$ P. All in all, Buss's witnessing does not seem to help us to obtain the unprovability of $LB(SAT,n^k)$ in $PV_1$ or $S^1_2(bit)$. Maybe it could work for intuitionistic $S^1_2$ where the witnessing holds for arbitrarily complex formulas, cf. \cite{b2}. The situation is different in case of weaker theories where we have more efficient witnessing. This will allow us to reduce to some hardness assumptions.

\section{Theories weaker than $PV_1$}

In this section we consider some theories weaker than $PV_1$ like $T_{NC^1}$ for which we will show the unprovability of circuit lower bounds. We could however similarly define a general theory $T_C$ corresponding to a standard complexity class $C$ and our results would work analogously.       

\begin{defi} $T_{NC^1}$ is the first-order theory of all universal $L_{NC^1}$ statements true in the standard model of natural numbers where $L_{NC^1}$ is the language containing names for all uniform $NC^1$ algorithms. Analogously, $T_{PV}$ resp. $T_{AC^0}$ is the true universal theory in the language $L_{PV}$ resp. $L_{AC^0}$ containing names for all p-time algorithms resp. uniform families of $AC^0$ circuits.
\end{defi}

These theories are universal so they admit the KPT theorem from \cite{k4}: 

\begin{thm}[\cite{k4}] If $T_{NC^1}\vdash \exists y A(x,y)$ for open formula $A$, then there is a function $f$ in uniform $NC^1$ such that $A(x,f(x))$ holds for any $x$.

If $T_{NC^1}\vdash\exists y\forall z A(x,y,z)$ for open formula $A$, there are finitely many functions $f_1,...,f_k$ in uniform $NC^1$ such that $$T_{NC^1}\vdash A(x,f_1(x),z_1)\vee A(x,f_2(x,z_1),z_2)\vee...\vee A(x,f(x,z_1,...,z_{k-1}),z_k)$$

Analogously for $T_{AC^0}$ and $T_{PV}$ for which the resulting functions are in uniform $AC^0$ resp. in P. 
\end{thm}

In the field of Bounded Arithmetic there are also standard theories corresponding to uniform $AC^0, NC^1$ and other complexity classes, cf. \cite{c}. Typically, they are presented as two-sorted theories having one sort of variables representing numbers and the second sort of variables representing bounded sets of numbers. The first-sort (number) variables are denoted by lower case letters $x,y,z,...$ and the second-sort (set) variables by capital letters $X,Y,Z,...$ The underlying language includes the symbols $+,\cdot, =,\leq, 0,1$ of first-order arithmetic. In addition it contains symbol $=_2$ interpreted as equality between bounded sets of numbers, $|X|$ for the function mapping an element $X$ of the set sort to the largest number in $X$ plus one, and $\in$ for the relation which holds for a number $n$ and set $X$ if and only if $n$ is an element of $X$. 
\smallskip

Bounded quantifiers for sets have the form $\exists X\leq t\ \phi$ which stands for $\exists X \ (|X|\leq t\wedge \phi)$ or $\forall X\leq t\ \phi$ for $\forall X\ (|X|\leq t\rightarrow \phi)$. Here $t$ is number term which does not involve $X$. $\Sigma^B_0$ formulas are formulas without bounded quantifiers for sets but may have bounded number quantifiers. Each bounded set $X\leq t$ can be seen also as a finite binary string of size $\leq t$ which has 1 in the $i$-th position iff $i\in X$. When we say that a function $f(x,X)$ mapping bounded sets and numbers to bounded sets is in $AC^0$ or $NC^1$ we mean that the corresponding function on finite binary strings and unary representation of $x$  is in $AC^0$ or $NC^1$.
\medskip

The base theory we will consider is $V^0$ consisting of a set of basic axioms capturing the properties of symbols in the two-sorted language and a comprehension axiom schema for $\Sigma^B_0$-formulas stating that for any $\Sigma^B_0$ formula there exists a set containing exactly the elements that satisfy the formula, cf. \cite{c}. Further, Cook and Nguyen define theory $VNC^1$ as $V^0$ extended by the axiom that every monotone formula has an evaluation, see \cite{c}. 

\begin{thm}[Cook-Nguyen \cite{c}] If $VNC^1\vdash \forall x\forall X\exists Y A(x,X,Y)$ for $\Sigma^B_0$-formula $A$, there is a function $f$ in uniform $NC^1$ such that $A(x,X,f(x,X))$ holds for any $x,X$.

If $VNC^1\vdash\forall x\forall X\exists Y\forall Z A(x,X,Y,Z)$ for $\Sigma^B_0$-formula $A$, there are finitely many functions $f_1,...,f_k$ in uniform $NC^1$ such that \\

$A(x,X,f_1(x,X),Z_1)\vee A(x,X,f_2(x,X,Z_1),Z_2)\vee ...$

\ \ \ \ \ \ \qquad\qquad\qquad\qquad\qquad\qquad $... \vee A(x,X,,f(x,X,Z_1,...,Z_{k-1}),Z_k)$\\

Analogously for $V^0$ with the resulting functions in uniform $AC^0$. 

\end{thm}  

$LB(SAT,n^k)$ translates to the two-sorted language as follows\\

$\forall n>n_0\forall C\exists Y\leq n \exists A\leq n\ \forall W\leq n^k\forall Z\leq n [Comp(C,Y,W)\rightarrow$

\ \ \ $(C(Y;W)=1\wedge \neg SAT(Y,Z))\vee (C(Y;W)=0\wedge \ SAT(Y,A))]$\\ 

where $k,n_0$ are constants as before and $Comp(C,Y,W), C(Y;W)=0/1$, $SAT(Y,Z)$ are defined as their first-order counterparts but function $x_i$ is replaced by $i\in X$. 
\smallskip

Similarly, we obtain the two-sorted $SCE(SAT,n^k), DCE(SAT,n^k)$.  
\medskip

Let us also specify the formalization of $LB(SAT,n^k)$ in $T_{AC^0}$. $L_{AC^0}$ contains symbols for $SAT(y,z), Comp(C,y,w)$ and all the predicates we explicitly defined as $\Sigma^b_0(bit)$-formulas because they are not just p-time but in fact constant-depth formulas. Moreover, even if multiplication is not in $L_{AC^0}$ (but in $L_{NC^1}$) we may assume that the $L_{AC^0}$ functions $Comp(C,y,w), C(y;w)=1/0$ contain the bound $|w|\leq |y|^k$.  For simplicity, whenever we speak about $LB(SAT,n^k)$ in $T_{AC^0}$ we mean its formalization where instead of the $\Sigma^b_0(bit)$-formulas we have the respective symbols of $L_{AC^0}$. Similarly for $SCE(SAT,n^k)$, $DCE(SAT,n^k)$ and $T_{NC^1}$. Therefore, $LB(SAT,n^k)$, $SCE(SAT,n^k)$ and $DCE(SAT,n^k)$ in $T_{AC^0}$ and $T_{NC^1}$ have the form $\exists y\forall z \ A(x,y,z)$ for an open formula $A$ (i.e. $A$ has no quantifiers).   
\smallskip

The situation with the provability of polynomial circuit lower bounds in weak theories like $T_{NC^1}, VNC^1, T_{AC^0}...$ is less natural because they cannot fully reason about p-time concepts. In particular, there is a formula $LB_2(SAT,n^k)$ which is equivalent to $LB(SAT,n^k)$ in $S^1_2(bit)$ but not necessarily in $T_{NC^1}$. $LB_2(SAT,n^k)$ is like $LB(SAT,n^k)$ but with $LB(C,y,a)$ expressed positively:\\

$LB_2(SAT,n^k):$
\smallskip

$\forall 1^n>n_0\forall C\ \exists y,a,w\ |a|<|y|=n, |w|\leq n^k\forall z, |z|<|y| \ [\neg Circ(C,y,w)\vee$

\ \ \ \ \ \ $(C(y;w)=0\wedge SAT(y,a))\vee (C(y;w)=1\wedge \neg SAT(y,z))$\\ 


Analogously define $DCE_2(SAT,n^k)$, $SCE_2(SAT,n^k)$ and their two-sorted and $L_{AC^0}$ formulations. 
\smallskip

By the witnessing theorem above, if $T_{NC^1}$ proves $LB(SAT,n^k)$, $LB(SAT,n^k)$ has an $NC^1$ S-T protocol with $O(1)$ rounds which is S-T protocol with $O(1)$ rounds and uniform $NC^1$ S. If $T_{NC^1}\vdash LB_2(SAT,n^k)$, $LB_2(SAT,n^k)$ has an $NC^1$ S-T protocol with $O(1)$ rounds which is defined analogously as for $LB(SAT,n^k)$ but with S producing also computations $w$ of given circuits. As $DCE_2(SAT,n^k)$ has the form $\exists y A(x,y)$ for an open $A$ in $L_{AC^0}$, its provability in $T_{NC^1}$  implies $DCE_2(SAT,n^k)\in NC^1$. Here again, $DCE_2(SAT,n^k)\in NC^1$ is defined as $DCE(SAT,n^k)\in NC^1$ but with the witnessing algorithm producing also computations $w$ of given circuits. Analogously for theories $T_{AC^0}, V^0, VNC^1$ .

\section{Unprovability of circuit lower bounds in subtheories of $PV_1$}

To prove that $VNC^1$ or $T_{NC^1}$ do not prove $LB(SAT,n^k)$ it suffices to show that $LB(SAT,n^k)$ has no S-T protocol with $O(1)$ rounds where S is in uniform $NC^1$. For the unprovability of $LB_2(SAT,n^k)$ it however suffices to refute the existence of S-T protocols with $O(1)$ rounds where S $\in NC^1$ produces $w$'s (computations of given circuits) itself. This is quite simple:   

\begin{prop}\label{lb2} $LB(SAT,n^{k+1})\notin NC^1$, $DCE_2(SAT,n^{k+1})\notin NC^1$ and $LB_2(SAT,n^{k+1})$ has no $NC^1$ S-T protocol with $poly(n)$ rounds unless \\ $SIZE(n^k)\subseteq NC^1$. Unconditionally, for any sufficiently big $k$, $LB(SAT,n^{k})\notin AC^0$, $DCE_2(SAT,n^{k})\notin AC^0$ and $LB_2(SAT,n^{k})$ has no $AC^0$ S-T protocol with $poly(n)$ rounds.  
\end{prop}

Proof: Assume first that $LB(SAT,n^{k+1})\in NC^1$, i.e. there are $NC^1$ circuits $D_m(x)$ such that for sufficiently big $n$ whenever $x\in\{0,1\}^m$ for $m=poly(n)$ encodes an $n^{k+1}$-size circuit $C_n$ with $n$ inputs, $D_m(x)$ outputs $y,a$ such that $$C_n(y)=0\wedge SAT(y,a) \ \ \ or \ \ \ C_n(x)=1\wedge \forall z\neg SAT(y,z)$$ Now any $n^k$-size circuits $B_n$ with $n$ inputs can be simulated by $NC^1$ circuits: For $b\in \{0,1\}^n$ and $z=(z_1,...,z_n)$ denote $R[B_n,b,z]$ the circuit with $n$ inputs $z$ but computing as $B_n$ on $b$, i.e. it does not use inputs $z$ at all. The size of $R[B_n,b,z]$ is $(n^k+n)$. Let $E_n(b)$ be an $AC^0$ circuit which uses description of $B_n$'s as advice and maps $b\in\{0,1\}^n$ to $x\in\{0,1\}^m$ encoding $R[B_n,b,z]$.
\smallskip

For each $b\in\{0,1\}^n$, use $D_m(E_n(b))$ to find $y,a$ and output 0 iff $SAT(y,a)$. 
\smallskip

Deciding $SAT(y,a)$ is by our formalization doable by constant-depth formulas. Therefore, for each $b$, we predict $B_n(b)$ with an $NC^1$ circuit. 
\medskip

If $LB(SAT,n^k)\in AC^0$, we would obtain $AC^0$ circuits for PARITY, which is impossible. 
\medskip

This construction works analogously for $DCE_2(SAT,n^k)$ and as well for $LB_2(SAT,n^k)$ because if there was some $NC^1$ S-T protocol for $LB_2(SAT, n^k)$ S would be forced to produce computations $w$ of given circuits.            \qed 
\smallskip


\begin{corollary} $T_{NC^1}\not\vdash DCE_2(SAT,n^{k+1})$ and $T_{NC^1}\not\vdash LB_2(SAT,n^k)$ unless $SIZE(n^k)\subseteq NC^1$. For any sufficiently big $k$, $V^0\not\vdash DCE_2(SAT,n^k)$ and $V^0\not\vdash LB_2(SAT,n^k)$.  
\end{corollary}

This simple observation does not work if we want to refute that $LB(SAT,n^k)$ has $NC^1$ S-T protocols because T can send to S a computation of the artificially attached circuit. Indeed by Proposition \ref{ll} $LB(SAT,n^k)$ has a uniform $AC^0$ S-T protocol with $poly(n)$ rounds under a plausible assumption. 
\smallskip

We can however show that $LB(SAT,n^k)$ has no $NC^1$ S-T protocols with $O(1)$ rounds under a hardness assumption. To show this we will use an interpretation of suitable NW-generators as p-size circuits which is due to Razborov \cite{r2} and Kraj\'{i}\v{c}ek's proof of a hardness of certain NW-generators for $T_{PV}$ \cite{k3}. It actually seems to be a relatively straightforward modification of the previous simple observation.

\begin{thm}\label{T} If there is $f\in SIZE(n^k)$ such that for all formulas $F_n$ of size $2^{O(n^{2/c})}$, $P_{x\in\{0,1\}^n}[F_n(x)=f(x)]<1/2+1/2^{O(n^{2/c})}$ for infinitely many $n$'s, then $LB(SAT,n^{2kc})$ has no $NC^1$ S-T protocol with $O(1)$ rounds.
\end{thm}

To prove the theorem we will use Nisan-Wigdewrson (NW) generators with specific design properties. Let $A=\{a_{i,j}\}^{i=1,...,m}_{j=1,...,n}$ be an $m\times n$ 0-1 matrix with $l$ ones per row. $J_i(A):=\{j\in \{1,...,n\}; a_{i,j}=1\}$ and $f:\{0,1\}^l\mapsto \{0,1\}$. Then define NW-generator based on $f$ and $A$, $NW_{f,A}:\{0,1\}^n\mapsto\{0,1\}^m$ as $$(NW_{f,A}(x))_i=f(x|J_i(A))$$ where $x|J_i(A)$ are $x_j$'s such that $j\in J_i(A)$.
\smallskip

For any $c\geq 4$, Nisan and Wigderson \cite{nis} constructed $2^n\times n^c$ 0-1 matrix $A$ with $n^{c/2}$ ones per row which is also $(n,n^{c/2})$-design meaning that for each $i\neq j$, $|J_i(A)\cap J_j(A)|\leq n$. Moreover, the matrix $A$ has such a property that there are $n^c$-size circuits which given $i\in \{0,1\}^n$ compute the set $J_i(A)$. Therefore, as it was observed by Razborov \cite{r2}, if $f$ is in addition computable by $n^k$-size circuits, for any $x\in \{0,1\}^{n^c}$, $(NW_{f,A}(x))_x$ is a function on $n$ inputs $y$ computable by circuits of size $n^{2kc}$.
\medskip

Proof(of Theorem \ref{T}): Let $f\in SIZE(n^k)$ and $A$ be a $2^n\times n^c$ $(n,n^{c/2})$-design defined above so for any $x$, $(NW_{f,A}(x))_y$ can be computed from $y$ by an $n^{2kc}$-size circuit. Assume that $LB(SAT,n^{2kc})$ has an $NC^1$ S-T protocol with $O(1)$ rounds. In particular, for each $n^{2kc}$-size circuit $C(y)$ computing $(NW_{f,A}(x))_y$ S either finds out the value of $C(y_1)$ by deciding (in $AC^0$) $SAT(y_1,a_1)$ for $y_1,a_1$ it produced itself or T will send to S $w_1,b_1$ such that $$(C(y_1;w_1)=0\vee \neg SAT(y_1,a_1))\vee (C(y_1;w_1)=1\vee SAT(y_1,b_1))$$ In the later case, S continues with its second try $y_2,a_2$. After at most $t\leq l$ rounds for some fixed constant $l$, S will successfully predict $C(y_t)$.
\smallskip

Let $E_{n^c}(x)$ be $AC^0$ circuits mapping $x\in\{0,1\}^{n^c}$ to a description of an $n^{2kc}$-size circuit with $n$ inputs $y$ computing the function $(NW_{f,A}(x))_y$. We will consider our S-T protocol only on inputs of the form $E_{n^c}(x)$.    
\smallskip

Kraj\'{i}\v{c}ek \cite{k3} showed that if $f$ is in NP$\cap$coNP with unique witnesses such S-T protocol allows us to approximate $f$ by a p-size circuit. We will inspect that his proof works also for $f$ in $P/poly$ and $NC^1$ S-T protocols. In addition we will assume that T in our S-T protocol operates as follows: whenever S outputs $y$ with some $a$, T answers with the lexicographically first satisfying assignment $b$ to $y$ and the unique computation $w$ of given circuit $y$. If there is no such $b$, T replies with a string of zeros. This should replace the uniqueness property assumed in \cite{k3}.
\smallskip

For $u\in \{0,1\}^{n^{c/2}}$ and $v\in\{0,1\}^{n^c-n^{c/2}}$ define $r_y(u,v)\in \{0,1\}^{n^c}$ by putting bits of $u$ into positions $J_y(A)$ and filling the remaining bits by $v$ (in the natural order). For each $x$ there is a trace $tr(x)=y_1,a_1,...,y_t,a_t, t\leq l$ of the S-T communication. 

\begin{clm} There is a trace $Tr=y_1,a_1,...,y_t,a_t, t\leq l$ and $a\in \{0,1\}^{n^c-n^{c/2}}$ such that $Tr=tr(r_{y_t}(u,a))$ for at least a fraction of $2/(3(2^{2n}))^t$ of all $u$'s. 
\end{clm}

$Tr$ and $a$ can be constructed inductively. There are at most $2^{2n}$ tuples $y_j,a_i$, hence there is $y_1,a_1$ such that at least $1/2^{2n}$ traces begin with it. Either there is $a\in \{0,1\}^{n^c-n^{c/2}}$ such that $y_1,a_1=tr(r_{y_1}(u,a))$ for at least $2/(3(2^{2n}))$ of all $u$'s or we can find $y_2,a_2$ such that at least $1/(3(2^{2n}))^2$ traces begin with $y_1,a_1,y_2,a_2$. For the induction step assume we have a trace $y_1,a_1,...,y_i,a_i$ such that at least $1/(3^{i-1}(2^{2n})^i)$ traces begin with it. Either there is $a\in \{0,1\}^{n^c-n^{c/2}}$ such that $y_1,a_1,...,y_i,a_i=tr(r_{y_i}(u,a))$ for at least $2/(3^i(2^{2n})^i)$ of all $u$'s or we can find $y_{i+1},a_{i+1}$ such that at least $1/(3^i(2^{2n})^{i+1})$ traces begin with $y_1,a_1,...,y_{i+1},a_{i+1}$. This proves the claim.
\smallskip

Fix now $Tr$ and $a$ from the previous claim. 
\smallskip

Because $A$ is $(n,n^{c/2})$-design, for any row $y\neq y_t$ at most $n$ $x_j$'s with $j\in J_y(A)$ are not set by $a$. Hence there are at most $2^n$ assignments $z$ to $x_j$'s with $j\in J_y(A)$ not set by $a$. For each such $z$ let $w_z,b_z$ be the T's advice after S outputs $y, a_i$ on any $x$ containing the assignment given by $z$ and $a$. By our choice of T, $b_z$ depends only on $y$ and $w_z$ is uniquely determined by $z$ (and $a$ which is fixed). Let $Y_y, y\neq y_t$ be the set of all these witnesses for all possible $z$'s. The size of each such $Y_y$ is $2^{O(n)}$.   
\medskip

Now we define a formula $F$ that attempts to compute $f$ and uses as advice $Tr, a$ and some $t$ sets $Y_y$. For each $u\in\{0,1\}^{n^{c/2}}$ produce $r_{y_t}(u,a)$ (this is in $AC^0$). Let $V$ be the set of those inputs $u$ for which $tr(r_{y_t}(u,a))$ either is $Tr$ or starts as $Tr$ and let $U$ be the complement of $V$. Define $d_0$ to be the majority value of $f$ on $U$. Then use S to produce $y_1',a_1'$. If $y_1',a_1'$ is different from $Tr$ output $d_0$. Otherwise, find the unique T's advice in $Y_{y_1}$. Again, this is doable by a constant depth formula of size $2^n$ which has $poly(n)$ output bits. It has the form $\bigvee_{z\in\{0,1\}^n}(z=r_{y_t}(u,a)|(J_{y_1}(A)\cap J_{y_t}(A))\rightarrow output=w_z\in Y_{y_1})$. In the same manner continue until S produces $y_t', a_t'$. If $y_t',a_t'$ differs from $Tr$ output $d_0$. Otherwise, output 0 iff $SAT(y_t,a_t)$.         
\medskip

$F$ is a formula with $n^{c/2}$ inputs and size $2^{O(n)}$ because producing $r_{y_t}(u,a)$ is in $AC^0$, searching for T's advice in $Y_i$'s is doable by constant-depth $2^{O(n)}$-size formulas, S is in $NC^1$ and the structure of S-T protocol can be described by a constant-depth formula of size $n^{O(1)}$: 
\smallskip

$(S(x)\notin Tr\rightarrow output=d_0)\wedge (S(x)\in Tr\rightarrow$ 

\ \ \  $((S(x,w_z,b_z)\notin Tr\rightarrow output=d_0)\wedge (...$

\ \ \ \ \ \ $(S(x,w_1,b_1,...,w_t,b_t)\notin Tr\rightarrow output=d_0)\wedge$

\ \ \ \ \ \ $(S(x,w_1,b_1,...,w_t,b_t)\in Tr\rightarrow (output=0\leftrightarrow SAT(y_t,b_t)))...)))$ 
\smallskip

By the choice of $Tr$, for at least a fraction $2/(3(2^n))^t$ of all $u\in\{0,1\}^{n^{c/2}}$ $F$ will successfully predict $f(u)$. Moreover, at most $1/(3(2^n))^t$ of all traces $tr(r_{y_t}(u,a))$ extend $Tr$. Because $d_0$ is the correct value on at least half of $u\in U$, $P_u[F(u)=f(u)]\geq 1/2+1/(3^t2^{nt+1})$   \qed

\begin{corollary} $T_{NC^1}\not\vdash LB(SAT,n^{2kc})$ and $VNC^1\not\vdash LB(SAT,n^{2kc})$ for $k\geq 1, c\geq 4$ unless for each $f\in SIZE(n^k)$ there are formulas $F_n$ of size $2^{O(n^{2/c})}$ such that for sufficiently big $n$'s, $P_{x\in\{0,1\}^n}[F_n(x)=f(x)]\geq 1/2+1/2^{O(n^{2/c})}$.
\end{corollary}

To obtain an unconditional unprovability of circuit lower bounds we can use Hastad's lower bound for constant depth circuits computing the parity function.    

\begin{thm}[Hastad \cite{has}] For any depth $d$ circuits $C_n$ of size $2^{n^{1/(d+1)}}$ and large enough $n$, $P_{x\in\{0,1\}^n}[C_n(x)=PARITY(x)]\leq 1/2+1/2^{n^{1/(d+1)}}$
\end{thm}

If $V^0\vdash LB(SAT,n^k)$, $LB(SAT,n^k)$ has an $AC^0$ S-T protocol with $O(1)$ rounds so the resulting formula $F$ in the proof of Theorem \ref{T} would be actually a constant-depth circuit and PARITY could be approximated by constant depth circuits of size $2^{O(n^{2/c})}$ with advantage $1/2^{O(n^{2/c})}$. This is not enough for the contradiction with Hastad's theorem. Nevertheless, it is sufficient if we replace polynomial circuit lower bounds $LB(SAT,n^k)$ by quasi polynomial lower bounds $LB(SAT,n^{\log n})$:\\

$\forall m>n_0 \forall C\exists y,a\ |a|<|y|=n \forall w, |w|\leq n^{\log n}=m [Comp(C,y,w)\rightarrow$

\ \ \  $(C(y;w)=0\wedge SAT(y,a))\vee (C(y;w)=1\wedge \forall z\neg SAT(y,z))]$\\

where $n$ is the number of inputs to $C$ and $m$ represents $n^{\log n}$ (or simply $|m|=|n|^2$). Analogously, define the two-sorted and $L_{AC^0}$ version of $LB(SAT,n^{\log n})$.

\begin{corollary} $T_{AC^0}\not\vdash LB(SAT,n^{\log n})$. $V^0\not\vdash LB(SAT,n^{\log n})$
\end{corollary}



\section{Acknowledgement}

I would like to thank Jan Kraj\'{i}\v{c}ek and Albert Atserias for useful discussions. This research was supported by grant GAUK 5732/2012 and partially by grants IAA100190902 of GA AV \v{C}R and SVV-2012-267317. A part of this research was done while I 
was a visiting fellow at the Isaac Newton Institute in Cambridge in Spring 2012 supported by grant N-SPP 2011/2012. 

}



\section{Introduction}

The aim of this paper is to show that a lot of complexity theory can be formalized in low fragments of arithmetic like Cook's theory $PV_1$.
\medskip

Our motivation is to demonstrate the power of bounded arithmetic as a counterpart to the unprovability results we already have or want to obtain, and generally to find out how complexity theory behaves in different worlds of bounded arithmetic.
\medskip

Concerning the unprovability results, Pich \cite{clba} proves that under certain hardness assumptions the theory $T_{NC^1}$, the true universal first-order theory in the language containing names for all uniform $NC^1$ algorithms, cannot prove polynomial circuit lower bounds on SAT formalized naturally by a sentence $LB(SAT,n^k)$. In fact, that result generalizes basically to any theory weaker than $PV_1$ in terms of provably total functions. The question whether $PV_1$ proves $LB(SAT,n^k)$ remains open even if we allow standard complexity-theoretic hardness assumptions, 
 see the discussion in Section \ref{lll}. 



Generally, it would be interesting to arrive at a complexity-theoretic statement, not necessarily circuit lower bounds, whose provability in $PV_1$ unexpectedly contradicts some other natural hypothesis. 
To understand better what are plausible candidates for such statements it might help us to investigate the theorems which are provable in low fragments of arithmetic.
\smallskip

In the present paper we will describe the formalization of just a few results; however, this should suffice to illustrate the power of the respective theories. Actually, many classical theorems from complexity theory have been already formalized in bounded arithmetic. In the table closing this section we list some representative examples. It should be understood that any of the formalized results is accompanied by a lot of other theorems that are formalizable in a similar fashion. In fact, some of the formalizations are so evident that they are used without a proof as a folklore. This is the case of Cook-Levin's theorem whose formalization we nevertheless describe for expository reasons in Section \ref{pv} as it gives us the opportunity to introduce some notions. For more details concerning the list see Section \ref{pre}. 
\medskip

The main original contribution of this paper is a formalization of the exponential PCP theorem in the theory $APC_1$ and the PCP theorem in the theory $PV_1$. Perhaps the most challenging part here was to formalize properties of the $(n,d,\lambda)$-graphs needed to derive the PCP theorem. These are usually obtained using algebraic techniques involving norms over real vector spaces coming all the way down to the fundamental theorem of algebra etc. In order to avoid formalization of this machinery (and it is not clear whether this could be done) we employ certain approximations to derive slightly weaker properties of the $(n,d,\lambda)$-graphs in the theory $PV_1$ which, however, suffice to derive the PCP theorem in $PV_1$. 

As the exponential PCP theorem follows trivially from the PCP theorem, the exponential version is actually also provable in $PV_1$. The $PV_1$ proof of the PCP theorem uses (among many other tools) the exponential PCP theorem but scaled down to constant size instances so that to prove the scaled down version we need to reason only about sets of constant size. On the other hand, in $APC_1$ we perform the standard proof of the exponential PCP theorem directly by formalizing a reasoning with p-time definable sets. Hence, the $APC_1$ proof shows different techniques to be available in low fragments of arithmetic.
\smallskip

The paper is organized as follows. In Section \ref{lll} we describe general properties of our formalizations and define theories of bounded arithmetic in which these formalizations take place. In Section \ref{pre} we discuss theorems that have been already formalized in bounded arithmetic as well as the new ones obtained in this paper. Section \ref{pv} illustrates a formalization of the Cook-Levin theorem in $PV_1$. In Section \ref{epcp} we prove the exponential PCP theorem in $APC_1$. Section \ref{rest} formalizes pseudorandom constructions in $PV_1$ which are then used in Section \ref{npcp} to formalize the PCP theorem in $PV_1$.

\smallskip
\begin{tabular}{l l c}
 
Theory & Theorem & Reference \\
\hline\hline

$PV_1$ & Cook-Levin's theorem & Section \ref{pv} \\

 & $(n,d,\lambda)$-graphs & Section \ref{rest} \\

 & the PCP theorem & Section \ref{npcp} \\

$PV_1+WPHP(PV_1)$ & PARITY $\notin AC^0$ & \cite{k} \\

$APC_1$ & BPP, ZPP, AM,... & \cite{e} \\
 & Goldreich-Levin's theorem & \cite{dai} \\
 & the exponential PCP theorem & Section \ref{epcp} \\

$HARD_{\epsilon}$ & Impagliazzo-Wigderson's derandom. & \cite{jerp} \\
$HARD^A$ & Nisan-Wigderson's derandomization  & \cite{jer} \\


$T^1_2+rWPHP(PV_2)$ & $S^P_2\subseteq ZPP^{NP}$ & \cite{hash}  \\

$APC_2$ & Graph isomorphism in coAM & \cite{hash} \\
 
$APC^{\oplus_p P}_2$ & Toda's theorem & \cite{b2} \\
\hline

\end{tabular}
\smallskip

The theories are listed from the weakest to the strongest one.

\section{Formalizations in bounded arithmetic: initial notes}\label{lll}

The usual language of arithmetic contains well known symbols: $0,S,+,\cdot, =,\leq$. To encode reasoning about computations it is helpful to consider also symbols $\lfloor \frac{x}{2}\rfloor, |x|$ and $\#$ with the intended meaning ``the whole part of $\frac{x}{2}$", ``the length of the binary representation of $x$", and $x\# y=2^{|x|\cdot |y|}$. The language $L$ containing all these symbols was used by Buss \cite{b} to define the theory $S^1_2$ (see below).

All theories we will work with, a subset of theories collectively known as bounded arithmetic, contain $L$ as a part of their language.
\medskip

The defining properties of symbols from $L$ are captured by a set of basic axioms denoted as BASIC which we will not spell out, cf. Kraj\'{i}\v{c}ek \cite{k}.
\medskip

A quantifier is sharply bounded if it has the form $\exists x, x\leq |t|$ or $\forall x, x\leq |t|$ where $t$ is a term not containing $x$. 
A quantifier is bounded if it is existential bounded: $\exists x, x\leq t$ for $x$ not occuring in $t$, or universal bounded: $\forall x, x\leq t$ for $x$ not occuring in $t$. By $\Sigma^b_0$ (=$\Pi^b_0=\Delta^b_0$) we denote the set of all formulas in the language $L$ with all quantifiers sharply bounded. For $i\geq 0$, the sets $\Sigma^b_{i+1}$ and $\Pi^b_{i+1}$ are the smallest sets satisfying
\smallskip

\noindent\quad $\(a$ $\Sigma^b_i\cup\Pi^b_i\subseteq \Sigma^b_{i+1}\cap\Pi^b_{i+1}$

\noindent\quad $\(b$ $\Sigma^b_{i+1}$ and $\Pi^b_{i+1}$ are closed under $\wedge,\vee$ and sharply bounded quantification

\noindent\quad $\(c$ $\Sigma^b_{i+1}$ is closed under bounded existential quantification

\noindent\quad $\(d$ $\Pi^b_{i+1}$ is closed under bounded universal quantification
 
\noindent\quad $\(e$ the negation of a $\Sigma^b_{i+1}$-formula is $\Pi^b_{i+1}$

\noindent\quad $\(f$ the negation of a $\Pi^b_{i+1}$-formula is $\Sigma^b_{i+1}$. 
\smallskip

\noindent In words, the complexity of bounded formulas in language $L$ (formulas with all quantifiers bounded) is defined by counting the number of alternations of bounded quantifiers, ignoring the sharply bounded ones. For $i>0$, $\Delta^b_i$ denotes $\Sigma^b_i\cap\Pi^b_i$.

\medskip

An example of a bounded arithmetic theory is the theory $S^1_2$ introduced by Buss \cite{b}. The language of $S^1_2$ is $L$ and its axioms consist of BASIC and $\Sigma^b_1$-PIND scheme which is the following kind of polynomial induction for $\Sigma^b_1$-formulas $A$: $$A(0)\wedge \forall x, (A(\lfloor x/2 \rfloor)\rightarrow A(x))\rightarrow \forall x A(x)$$ Buss \cite{b} showed that whenever $S^1_2$ proves a formula of the form $\exists y, A(x,y)$ for $\Sigma^b_1$-formula $A$, then there is a p-time (i.e. polynomial time) function $f$ such that $A(x,f(x))$ holds for all $x$.
\smallskip

Theories of bounded arithmetic generally cannot prove the totality of functions with superpolynomial growth of length. This follows from a theorem of Parikh \cite{pa}. In particular, $\forall k \ \exists x, |x|=k$ is unprovable. Consequently, if we want to prove in bounded arithmetic a statement of the form ``for all $k,n$, there is an $n^k$-size circuit (encoded by a binary string of some number, i.e. $\exists x, |x|=n^k$) s.t. ..." we need to quantify the exponent $k$ outside of the respective theory. That is, in such cases instead of proving
\[\mbox{$T\vdash$ ``for all $k,n$, there is an $n^k$-size circuit s.t. ..."}\]
 we prove
\[\mbox{``for all $k$, $T\vdash$ for all $m,n$ s.t. $|m|=n$, there is an $n^k$-size circuit s.t. ..."}\]

\noindent Informally speaking, only the ``feasible part" of the theorem is provable inside the theory. 
\smallskip

In our formalizations numbers encode binary strings in a natural way. 
We then follow the convention that inputs of circuits, algorithms or functions are represented by binary strings. For example, when talking about $n^k$-size circuit lower bounds the number of inputs of $n^k$-size circuits is the length of some number, i.e $\exists x,\ n=|x|$. However, it does not necessarily follow that $n$ is smaller, say, $\exists x, \ n=||x||$. 
To indicate sizes of objects inside our theories we employ the shorthand notation $x\in Log\leftrightarrow \exists y, x=|y|$ and $x\in LogLog\leftrightarrow \exists y, x=||y||$. 

On the contrary, for example Razborov \cite{r3} considered (second-order) formalizations of circuit lower bounds (corresponding in first-order logic to the formalization) where p-size (i.e. polynomial size) circuits with $n$ inputs were required to satisfy $n\in Log Log$. Thus, in his formalization, truth tables of functions computed by p-size circuits are encoded by binary strings. The respective theory is much stronger with respect to such formalization; it is as if it could manipulate with exponentially big objects. Formalizing known theorems is then easier and proving unprovability results is on the other hand formally much harder.

Similarly, in propositional proof complexity there are candidate hard tautologies for strong proof systems like Extended Frege which express circuit lower bounds on SAT (and other functions), see formulas $\neg Circuit_t(f)$ in Razborov \cite{r2} or $\tau(tt_{s,k})_f$ in Kraj\'{i}\v{c}ek \cite{k2}. Using a standard translation into first-order logic they again correspond to the formalization where truth tables of SAT are encoded by binary strings. Therefore, by the known relation between propositional proof systems and bounded arithmetics, the hardness of such formulas for Extended Frege would imply a conditional unprovability of superpolynomial circuit lower bounds on SAT in $PV_1$ formalized in such a way that the theory $PV_1$ would be exponentially stronger than it is with respect to the formalization of circuit lower bounds $LB(SAT,n^k)$ considered in Pich \cite{clba}. The formalization $LB(SAT,n^k)$ follows the convention of our current paper. 

However, the fact advocated here, that a lot of complexity theory is formalizable in theories like $PV_1$, suggests that it might be also hard to obtain the unprovability of $LB(SAT,n^k)$ in $PV_1$. Actually, the unprovability of $LB(SAT,n^k)$ in $PV_1$ would imply that there is no provable witnessing of errors of p-time algorithms claiming to solve SAT which is itself (interesting and) a reason to expect hardness of such unprovability result, see Pich \cite{clba}.

\subsection{Theory \texorpdfstring{$PV_1$}{PV1}: formalized p-time reasoning}\label{pvdef}

$PV_1$ introduced in Kraj\'{i}\v{c}ek-Pudl\'{a}k-Takeuti \cite{k4} is a conservative extension of an equational theory $PV$ introduced by Cook \cite{c3}.

\medskip

The language of $PV$ and $PV_1$ consists of symbols for all p-time algorithms given by Cobham's characterization of p-time functions, cf. \cite{cob}. In particular, it contains $L$. By a slight abuse of the notation we denote the language of $PV_1$ and $PV$ also $PV$. A $PV$-formula is a first-order formula in the language $PV$. The hierarchy of $\Sigma^b_i(PV)$- and $\Pi^b_i(PV)$-formulas is defined similarly to $\Sigma^b_i$ and $\Pi^b_i$ (in first-order logic with equality) but in the language of $PV$.

\medskip

In $PV$ we can define p-time concepts and prove their basic properties. More precisely, every p-time function can be straightforwardly defined as a $PV$-function. Therefore, in the theory $PV_1$, which is a universal first-order theory, we can reason about p-time concepts. 
We can interpret provability in $PV_1$ as capturing the idea of what can be demonstrated when our reasoning is restricted to manipulation of p-time objects. However, strictly speaking, this description would also fit the theory $S^1_2$ which in addition uses NP-concepts in induction. Anyway, it is a natural question which properties of p-time concepts are provable using only such p-time reasoning.

\def\banal{
\begin{defi} A function $f$ is defined from functions $g,h_0,h_1$ and $l$ by limited recursion on notation if:
\smallskip

1. $f(x,0)=g(x)$

2. $f(x,s_i(y))=h_i(x,y,f(x,y))$, for $i=0,1$

3. $f(x,y)\leq l(x,y)$

\smallskip

\noindent where $s_0(y), s_1(y)$ are functions adding 0, resp. 1, to the right of the binary representation of $y$, i.e. $s_0(y):=2y$, $s_1(y)=2y+1$.
\end{defi}

\begin{thm}[Cobham \cite{cob}] The set of polynomial time functions is the smallest set of functions containing constant $0$, functions $s_0(y),s_1(y), x\#y$, and closed under:

1. permutation and renaming of variables,

2. composition of functions,

3. limited recursion on notation.
\end{thm}

\begin{defi}[Cook \cite{c3}]\label{newp} We simultaneously define function symbols of rank $k$ and $PV$-derivations of rank $k, k=0,1,...$. The language of $PV$ will then consist of all function symbols of any rank, and a $PV$-derivation will be a $PV$-derivation of any rank.
\medskip

\noindent 1. Function symbols of rank 0 are constant 0, unary $s_0(y),s_1(y)$ and $Tr(x)$; and binary $x_{\frown} y, x\# y$, and $Less(x,y)$

\medskip

\noindent 2. Defining equations of rank $0$ are:
\smallskip

$Tr(0)=0$

$Tr(s_i(x))=x$, $i=0,1$
\smallskip

$x_{\frown} 0=x$

$x_{\frown} (s_i(y))=s_i(x_{\frown} y)$, $i=0,1$
\smallskip

$x\# 0=0$

$x\# s_i(y)=x_{\frown} (x\# y)$, $i=0,1$
\smallskip

$Less(x,0)=x$

$Less(x,s_i(y))=Tr(Less(x,y))$, $i=0,1$
\smallskip

\noindent $Tr(x)$ deletes the rightmost bit, $x_{\frown} y$ is the concatenation, $x\# y$ is $|y|$ concatenated copies of $x$, and $Less(x,y)$ is $x$ with $|y|$ right bits deleted

\medskip

\noindent 3. PV rules are as follows. Let $t,u,v,t_1,u_1,...,t_k,u_k,f,f_1,f_2$ be function symbols.

R1. from $t=u$ derive $u=t$

R2. from $t=u, u=v$ derive $t=v$ 

R3. from $t_1=u_1,...,t_k=u_k$ derive $f(t_1,...,t_k)=f(u_1,...,u_k)$

R4. from $t=u$ derive $t(x/v)=u(x/v)$

R5. Let $E_1,...,E_6$ be the equations (1-3) from the definition of the limited recursion on notation: three for $f_1$ and three for $f_2$ in place of $f$. Then from $E_1,...,E_6$ $PV$ can derive $f_1(x,y)=f_2(x,y)$

\medskip

\noindent 4. PV-derivations of rank $k$ are sequences of equations $E_1,...,E_t$ in which every function symbol is of rank $\leq k$ and every $E_i$ is either a defining equation of rank $\leq k$ or derived from some earlier equations by one of the PV-rules

\medskip

\noindent 5. Let $t$ be a term consisting of function symbols of rank $\leq k$. Then $f_t$ is a function symbol of rank $k+1$ and $f_t=t$ is a defining equation of rank $k+1$.

\medskip

\noindent 6. Other function symbols of rank $k+1$ are obtained as follows. Whenever $g,h_0,h_1,l_0,l_1$ are function symbols of rank $\leq k$ and $\pi_0,\pi_1$ are PV-derivations of rank $k$ of equality $Less(h_i(x,y,z),z_{\frown} l_i(x,y))=0$ then $f=f_{(g,h_0,h_1,l_0,l_1,\pi_0,\pi_1)}$ is a function symbol of rank $k+1$, and the equations defining $f$ from $g,h_i$ by limited recursion on notation are defining equations of rank $k+1$.

\end{defi}

$PV_1$ is a theory in the first-order predicate logic which consists of all equations provable in $PV$ but has also a form of induction axiom: For an open formula $\psi(x)$ define a function $h(b,y)$ by
\smallskip

(a) $h(b,0)=(0,b)$

(b) if $h(b,\lfloor u/2\rfloor)=(x,y)$ and $u>0$ then 
\smallskip

\ \ \ $h(b,u):=(\lceil x+y/2\rceil,y)$ if $\lceil x+y/2\rceil<y\wedge \psi(\lceil x+y/2\rceil)$

\ \ \ \ \ \ \ \ \ \ \ \ $:=(x,\lceil x+y/2\rceil)$ if $x<\lceil x+y/2\rceil\wedge \neg\psi(\lceil x+y/2\rceil)$

\ \ \ \ \ \ \ \ \ \ \ \ $:=(x,y)$ otherwise

\smallskip

Then $PV_1$ contains the universal axiom $$(\psi(0)\wedge\neg \psi(b)\wedge h(b,b)=(x,y))\rightarrow (x+1=y\wedge\psi(x)\wedge \neg\psi(y))$$

Note that $PV_1$ is a universal theory.
}
\smallskip


It can be shown that $PV_1$ proves $\Sigma^b_0(PV)$-induction, cf. Kraj\'{i}\v{c}ek \cite{k}. That is, for any $\Sigma^b_0(PV)$-formula $A$, $PV_1$ proves $$A(0)\wedge\forall x (A(x)\rightarrow A(x+1))\rightarrow \forall x A(x)$$

In $PV$ we can speak about formulas, circuits, Turing machines and other similar notions which can be encoded using finite sequences of numbers. These are encodable in $PV$ in a well-behaved way so that basic operations on sequences like concatenation are definable by terms, i.e. by functions in the language. For more details see Kraj\'{i}\v{c}ek \cite{k} where the function $(w)_i$ which extracts the $i$th element from a sequence $w$ is shown to be $\Delta^b_1$-definable in $S^1_2$ but the definition is given by a p-time predicate so it can be written as an open $PV$-formula.
\smallskip

All $PV$-functions have well-behaved $\Delta^b_1$-definitions in $S^1_2$. Hence, $S^1_2$ can be seen as an extension of $PV_1$, cf. Buss \cite{b}. Moreover, Buss's witnessing theorem \cite{b} implies that $S^1_2$ is $\forall\Sigma^b_1$-conservative over $PV_1$. This means that when proving a $\forall\Sigma^b_1$ statement in $PV_1$ we can actually use $S^1_2$. In particular, we will use an induction scheme denoted as $\Pi^b_1$-LLIND which is provable in $S^1_2$ and says that for any $\Pi^b_1(PV)$-formula $A$ the following holds, $$A(0)\wedge\forall x\leq ||a||\ (A(x)\rightarrow A(x+1))\rightarrow A(||a||)$$

In Proposition \ref{cst}, we will also use an induction scheme which we denote $\Pi^b_1$-LPIND. It is a weaker form of $\Pi^b_1$-PIND, cf. Kraj\'{i}\v{c}ek \cite{k}, so it is derivable in $S^1_2$. $\Pi^b_1$-LPIND says that for any $\Pi^b_1(PV)$-formula $A$ the following implication holds: $$A(a)\wedge A(a^2)\wedge [\forall l\leq ||b||, (A(a^{\lfloor (l-1)/2\rfloor})\wedge A(a^{\lceil (l-1)/2\rceil})\rightarrow A(a^l))]\rightarrow A(a^{||b||})$$

\subsection{Theory \texorpdfstring{$APC_1$}{APC1}: formalized probabilistic p-time reasoning}

To reason about probabilistic p-time concepts we will use an extension of $PV_1$ in which Je\v{r}\'{a}bek \cite{e} developed a well-behaved notion of probability based on an approximate counting. 

In this section, we recall a part of his work which we will use to formalize the exponential PCP theorem.
\smallskip

The dual (or surjective) pigeonhole principle for $f$, written as $dWPHP(f)$, is the universal closure of the formula $$x>0\rightarrow\exists v<x(|y|+1)\forall u<x|y| f(u)\neq v$$ For a set of functions $\Gamma$, $dWPHP(\Gamma):=\{dWPHP(f)|f\in\Gamma\}$.
\medskip

The theory $APC_1$ is defined as $PV_1+dWPHP(PV)$ where $PV$ stands for the set of $PV$-functions.
\smallskip

When a number $a$ is used in a context which asks for a set it is assumed to represent the integer interval $[0,a)$, e.g. $X\subseteq a$ means that all elements of $X$ are less than $a$. If $X\subseteq a$, $Y\subseteq b$, then $X\times Y:=\{bx+y|x\in X, y\in Y\}\subseteq ab$ and $X\dot{\cup} Y:=X\cup\{y+a|y\in Y\}\subseteq a+b$.

We will often work with rational numbers which are assumed to be represented by pairs of integers in the natural way. By a definable set we mean a collection of numbers satisfying some formula, possibly with parameters.
\medskip

Let $n,m\in Log$, $C: 2^n\rightarrow 2^m$ be a circuit and $X\subseteq 2^n, Y\subseteq 2^m$ definable sets.
We write $C:X\twoheadrightarrow Y$ if $Y\subseteq C[X]$, i.e. $\forall y\in Y\exists x\in X,\ C(x)=y$. 
The following definitions are taken from Je\v{r}\'{a}bek \cite{e}.



\begin{defi}[in $APC_1$] Let $X,Y\subseteq 2^n$ be definable sets, and $\epsilon\leq 1$. We say that the size of $X$ is approximately less than the size of $Y$ with error $\epsilon$, written as $X\preceq_{\epsilon} Y$, if there exists a circuit $G$, and $v\neq 0$ such that $$G: v\times (Y\dot{\cup}\epsilon 2^n)\twoheadrightarrow v\times X$$ The sets $X$ and $Y$ have approximately the same size with error $\epsilon$, written as $X\approx_{\epsilon}Y$, if $X\preceq_{\epsilon} Y$ and $Y\preceq_{\epsilon} X$.
\end{defi}

A number $s$ identified with the interval $[0,s)$, so $X\preceq_{\epsilon} s$ means that the size of $X$ is at most $s$ with error $\epsilon$.

\begin{defi}[in $APC_1$] Let $X\subseteq 2^{|t|}$ be a definable set and $0\leq \epsilon, p\leq 1$. We define $$Pr_{x<t}[x\in X]\preceq_{\epsilon}p\ \ \ iff\ \ \ X\cap t\preceq_{\epsilon} pt$$ and similarly for $\approx$.
\end{defi}

The definition of $\preceq_{\epsilon}$ is an unbounded $\exists \Pi^b_2$-formula so it cannot be used freely in bounded induction. This problem was solved by Je\v{r}\'{a}bek \cite{e} by working in a suitable conservative extension of $APC_1$.

\begin{defi}[in $PV_1$] Let $f:2^k\mapsto 2$ be a truth-table of a Boolean function with $k$ inputs ($f$ is encoded as a string of $2^k$ bits, hence $k\in LogLog$). We say that $f$ is (worst-case) $\epsilon$-hard, written as $Hard_{\epsilon}(f)$ if no circuit $C$ of size $2^{\epsilon k}$ computes $f$. The function $f$ is average-case $\epsilon$-hard, written as $Hard^A_{\epsilon}(f)$, if for no circuit $C$ of size $\leq 2^{\epsilon k}$: $$|\{u<2^k|C(u)=f(u)\}|\geq (1/2+2^{-\epsilon k})2^k$$
\end{defi}

\begin{prop}[Je\v{r}\'{a}bek \cite{jer}] For every constant $\epsilon <1/3$ there exists a constant $c$ such that $APC_1$ proves: for every $k\in LogLog$ such that $k\geq c$, there exist average-case $\epsilon$-hard functions $f:2^k\mapsto 2$.  
\end{prop}

$PV_1$ can be relativized to $PV_1(\alpha)$. The new function symbol $\alpha$ is then allowed in the inductive clauses for introduction of new function symbols. 
This means that the language of $PV_1(\alpha)$, denoted also $PV(\alpha)$, contains symbols for all p-time oracle algorithms.

\begin{defi}[Je\v{r}\'{a}bek \cite{jer}] The theory $HARD^A$ is an extension of the theory $PV_1(\alpha)+dWPHP(PV(\alpha))$ by the axioms

 \begin{enumerate}[label=\arabic*.]
 \item $\alpha(x)$ is a truth-table of a Boolean function in $||x||$ variables 
 \item $x\geq c\rightarrow Hard^A_{1/4}(\alpha(x))$
 \item $||x||=||y||\rightarrow \alpha(x)=\alpha(y)$
 \end{enumerate}
   where $c$ is the constant from the previous lemma. 
\end{defi}

\begin{thm}[Je\v{r}\'{a}bek \cite{jer,e}]\label{sec} $HARD^A$ is a conservative extension of $APC_1$. Moreover, there is a $PV(\alpha)$-function $Size$ such that $HARD^A$ proves: if $X\subseteq 2^n$ is definable by a circuit $C$, then $$X\approx_{\epsilon} Size(C,2^n,e)$$ where $\epsilon=|e|^{-1}$
\end{thm} 

We will abuse the notation and write $Size(X,\epsilon)$ instead of $Size(C,2^n,e)$. 

\begin{defi}[in $APC_1$] If $X\subseteq 2^{|t|}$ is defined by a circuit and $\epsilon^{-1}\in Log$, we put $$Pr_{x<t}[x\in X]_{\epsilon}:=\frac{1}{t} Size(X\cap t,\epsilon)$$
\end{defi}

Je\v{r}\'{a}bek \cite{e} showed that these definitions are well-behaved:

\begin{prop}(in $PV_1$)\label{lem} Let $X,X',Y,Y',Z\subseteq 2^n$ be definable sets and $\epsilon, \delta<1$. Then
\begin{enumerate}[label=\roman*)]

\item $X\subseteq Y\Rightarrow X\preceq_{0} Y$

\item $X\preceq_{\epsilon} Y \wedge Y\preceq_{\delta} Z\Rightarrow X\preceq_{\epsilon+\delta} Z$

\item $X\preceq_{\epsilon} X'\wedge Y\preceq_{\delta} Y' \Rightarrow X\times Y\preceq_{\epsilon+\delta+\epsilon\delta} X'\times Y'$
\end{enumerate}
\end{prop}

\begin{prop}(in $APC_1$)\label{really} \hfill
\begin{enumerate}[label=\arabic*.]
\item Let $X,Y\subseteq 2^n$ be definable by circuits, $s,t,u\leq 2^n$, $\epsilon, \delta,\theta, \gamma \leq 1, \gamma^{-1}\in Log $. Then
\begin{enumerate}[label=\roman*)]
\item $X\preceq_{\epsilon} Y\Rightarrow 2^n- Y\preceq_{\epsilon+\delta} 2^n- X $

\item $X\approx_{\epsilon} s\wedge Y\approx_{\delta} t\wedge X\cap Y\approx_{\theta} u\Rightarrow X\cup Y\approx_{\epsilon+\delta+\theta+\gamma} s+t-u$
\end{enumerate}

\item Let $X\subseteq 2^n\times 2^m$ and $Y\subseteq 2^m$ be definable by circuits, $t\preceq_{\epsilon} Y$ and $s\preceq_{\delta} X_y$ for every $y\in Y$, where $X_y:=\{x| \left<x,y\right>\in X\}$. Then for any $\gamma^{-1}\in Log$ $$st\preceq_{\epsilon+\delta+\epsilon\delta+\gamma} X\cap (2^n\times Y)$$

\item (Chernoff's bound) Let $X\subseteq 2^n, m\in Log, 0\leq \epsilon,\delta,p\leq 1$ and $X\succeq_{\epsilon} p2^n$. Then $$\{w\in (2^n)^m|\ |\{i<m|w_i\in X\}|\leq m(p-\delta)\}\preceq_0 c4^{m(c\epsilon-\delta^2)}2^{nm}$$ for some constant $c$, where $w$ is treated as a sequence of $m$ numbers less than $2^n$ and $w_i$ is its $i$-th member.

\end{enumerate}
\end{prop}





\section{Previous formalizations of complexity theory and our contribution}\label{pre}

Many classical theorems from complexity theory have been already formalized in bounded arithmetic. In the following sections we present some representative examples from different areas of complexity theory. The last section describes the formalizations that are obtained in this paper.

\subsection{NP-completeness}

Actually, formalization of some theorems is a folklore used without a proof. For example, Cook-Kraj\'{i}\v{c}ek \cite{ck} mention that NP-completeness of SAT can be formalized in $PV_1$.

\begin{thm}[Cook-Levin's theorem in $PV_1$]\hfill
\begin{enumerate}[label=(\alph*)]
\item For every $\Sigma^b_1$-formula $\phi(x)$, there is a $PV$-function $f(x)$ such that $$PV_1\vdash \phi(x)\leftrightarrow \exists y SAT(f(x),y)$$ where $SAT(z,y)$ is an open $PV$-formula which holds iff truth assignment $y$ satisfies propositional formula $z$. 

\item For each $k$ we have a $PV$-function $f$ such that $PV_1$
  proves: for any $M,x$, $$\exists w, z; |z|, |w|\leq |x|^k,
  M(x,z,w)=1\leftrightarrow \exists y, |y|\leq 3|M||x|^{2k},
  SAT(f(M,x),y)$$ where $M(x,z,w)=1$ is an open $PV$-formula which
  holds iff $w$ is an accepting computation of Turing machine $M$ on
  input $x,z$ (so we are slightly abusing the notation as $M$ is
  actually a free variable in the formula $M(x,z,w)=1$) and $|M|$ is
  the length of $M$'s code.
\end{enumerate}
\end{thm}

\noindent Note that formulations $(a)$ and $(b)$ are essentially equivalent since the formula $\exists w,z; |z|, |w|\leq |x|^k, M(x,z,w)=1$ is $\Sigma^b_1$ and any $\Sigma^b_1$-formula $\phi(x)$ is equivalent in $PV_1$ to a formula $\exists w, z; |z|, |w|\leq |x|^k, M(x,z,w)=1$ for some $k$ and $M$. In $(b)$ we have in addition also an explicit bound on $y$.

For expository reasons we present a proof of $(b)$ in Section \ref{pv}.



\subsection{Randomized computation}

The main application of approximate counting in $APC_1$ is in the formalization of probabilistic algorithms in $APC_1$ and complexity classes like BPP and AM. Je\v{r}\'{a}bek's formalizations involve many other results we will not state explicitly like ``promise BPP $\subseteq$ P/poly'' (Lemma 3.10 in Je\v{r}\'{a}bek \cite {e}), Rabin-Miller algorithm (Example 3.2.10 in Je\v{r}\'{a}bek \cite{jerp}) but also principles like Stirling's bound on binomial coefficients.

\begin{defi}[Je\v{r}\'{a}bek \cite{e}](in $APC_1$) A $PV$-function $r$ and a $PV$-predicate $A$ define a BPP language if for each $x$ either $ Pr_{w<r(x)}[\neg A(x,w)]\preceq_0 1/4$ or $Pr_{w<r(x)}[A(x,w)]\preceq_0 1/4$.
\end{defi}

\begin{thm}[Je\v{r}\'{a}bek \cite{e}] Let $A$ be a $PV$-predicate and $r$ a $PV$-function. There are $\Sigma^b_2$-formulas $\sigma^+(x),\sigma^-(x)$ and $\Pi^b_2$-formulas $\pi^+(x),\pi^-(x)$ such that $APC_1$ proves $$Pr_{w<r(x)}[\neg A(x,w)]\preceq_0 1/4\Rightarrow \pi^+(x)\Rightarrow \sigma^+(x)\Rightarrow Pr_{w<r(x)}[\neg A(x,w)]\preceq_0 1/3$$ $$Pr_{w<r(x)}[A(x,w)]\preceq_0 1/4\Rightarrow \pi^-(x)\Rightarrow \sigma^-(x)\Rightarrow Pr_{w<r(x)}[A(x,w)]\preceq_0 1/3$$
In particular, any definable $BPP$ language is in $\Sigma^b_2\cap\Pi^b_2$.
\end{thm}

\smallskip

In \cite{hash} Je\v{r}\'{a}bek formalized Cai's \cite{cai} result stating that $S^P_2\subseteq ZPP^{NP}$ in the theory $T^1_2+rWPHP(PV_2)$. The complexity class $S^P_2$ consists of languages for which there exists a p-time predicate $R$ such that $$x\in L\Rightarrow \exists y\forall z R(x,y,z)$$ $$x\notin L\Rightarrow \exists z\forall y \neg R(x,y,z)$$ where $|y|,|z|$ are implicitly bounded by a polynomial in $|x|$.
\smallskip

The theory $T^1_2$ is defined as $S^1_2$ but with induction for $\Sigma^b_1$-formulas, $PV_2$ denotes functions computable in polynomial time relative to NP, and  $rWPHP(PV_2)$ is a set of axioms $$x>0\rightarrow \exists y<x(|y|+1) (g(y)\geq x|y|\vee f(g(y))\neq y)$$ for $PV_2$-functions $f,g$. 

Note that $rWPHP(f,g)$ follows from $dWPHP(f)$.

\begin{thm}[Je\v{r}\'{a}bek \cite{hash}](in $T^1_2+rWPHP(PV_2)$) The complexity class $S^P_2$ is contained in $ZPP^{NP}$. That is, for each p-time relation $R$ defining a language $L\in S^P_2$, there exists $ZPP^{NP}$-predicate $P$ definable in $T^1_2+rWPHP(PV_2)$ such that the same theory proves $x\in L\Leftrightarrow P(x)$. 
 \end{thm}

\subsection{Circuit lower bounds}\label{clbs}

In \cite[Section 15.2]{k} Kraj\'{i}\v{c}ek proves PARITY $\notin AC^0$ in the theory $PV_1+WPHP(PV_1)$. By $WPHP(PV_1)$ he denotes the set of axioms $$a>0\rightarrow \exists y\leq 2a\forall x\leq a, f(x)\neq y$$ for every $PV_1$-function symbol $f(x)$ where $f$ may have other arguments besides $x$ and they are treated as parameters in the axioms. 

It is known that $WPHP(PV_1)$ and $dWPHP(PV)$ are equivalent over $S^1_2$. Further, the theory $PV_1+dWPHP(PV)$ is $\forall\Sigma^b_1$-conservative over $$PV_1+\{\exists y< a\#a\ \forall x< a,\ f(x)\neq y| \mbox{ for PV-functions }f\}$$ (noted in Je\v{r}\'{a}bek \cite{jerr} as a corollary of earlier results).

\begin{thm}[Kraj\'{i}\v{c}ek \cite{k}, Section 15.2] Let $d, k$ be arbitrary constants. Then the theory $PV_1+WPHP(PV_1)$ proves that for any sufficiently large $n\in Log$ there are no depth $d$ circuits of size $\leq kn^k$ computing $PARITY(x_1,...,x_n)$.
\end{thm}

In \cite{r3} Razborov developes a logical formalism supporting his feeling that $S^1_2$ is the right theory to capture that part of reasoning in Boolean complexity which led to actual lower bounds for explicitly given Boolean functions. He formalizes lower bounds for constant-depth circuits over the standard basis, lower bounds for monotone circuits, lower bounds for constant-depth circuits with MOD-$q$ gates, and lower bounds for monotone formulas based on communication complexity.

Importantly, his formalizations presented in second-order logic correspond in first-order logic to the formalization where the number of inputs of circuits in the respective theorems is in $LogLog$. This makes it more suitable for encoding into the propositional setting but it also makes the formalization results formally weaker.

\subsection{Interactive proofs}

Je\v{r}\'{a}bek \cite{hash} formalized the equivalence of public-coin and private-coin interactive protocols in the theory $APC_2:=T^1_2+dWPHP(PV_2)$. This is illustrated on the example of the isomorphism problem: given two structures $G_0$ and $G_1$ (as tables) of the same signature, determine whether $G_0\simeq G_1$.

\begin{defi}[Je\v{r}\'{a}bek \cite{e}](in $APC_2$) A pair $\left<\phi,r\right>$ where $\phi(x,w)$ is a $\Sigma^b_1$-formula, and $r$ is a $PV$-function, defines an $AM$ language if for each $x$ either $Pr_{w<r(x)}[\neg \phi(x,w)]\preceq^1_0 1/4$ or $Pr_{w<r(x)}[\phi(x,w)]\preceq^1_0 1/4$ where $\preceq^1_0$ denotes $\preceq_0$ relativized with a $\Sigma^b_1$-complete oracle.
\end{defi}

\begin{thm}[Je\v{r}\'{a}bek \cite{hash}](in $APC_2$) Graph nonisomorphism is in AM.
\end{thm}

\subsection{Cryptography}

Recently, Dai Tri Man Le \cite{dai} formalized Goldreich-Levin's theorem in $APC_1$.

\begin{thm}[Dai Tri Man Le \cite{dai}](in $APC_1$) Let $f:\{0,1\}^n\rightarrow \{0,1\}^n$ be a function computed by a circuit of size $t$, and suppose that there exists a circuit $C$ of size $s$ such that $$Pr_{(x,r)\in\{0,1\}^{2n}}[C(f(x),r)=\bigoplus^n_{i=1}x_ir_i]_{\epsilon}\geq \frac{1}{2}+\frac{1}{p(n)}$$ If $\epsilon=\frac{1}{poly(n)}$ is sufficiently small, then there is a circuit $C'$ of size at most $(s+t)poly(n,1/\epsilon)$ and $q=poly(n)$ such that $$Pr_{(x,r')\in \{0,1\}^n\times \{0,1\}^q}[f(C'(f(x),r'))=f(x)]_{\epsilon}\geq \frac{1}{4p(n)}-\frac{15\epsilon}{2}$$
\end{thm}

\subsection{Complexity of counting}

In \cite{b2}, Buss, Ko\l odziejczyk and Zdanowski derived Toda's theorem in an extension of the theory $APC_2$.
\medskip

For a fixed prime $p\geq 2$, they denote by $C^k_p$ for $k\in [p]$ quantifiers counting mod $p$. The intended meaning of $C^k_p x\leq t A(x)$ is that the number of values $x\leq t$ for which $A$ is true is congruent to $k$ mod $p$. See \cite{b2} for the explicit list of axioms defining $C^k_p$.

A $\oplus_p P$ formula is a formula which is either atomic, or of the form $C^k_p x\leq t A(x)$ where $A$ is sharply bounded. $\Sigma_0^{b,\oplus_p P}=\Pi^{b,\oplus_p P}_0$ is the set of formulas obtained as the closure of $\oplus_p P$ formulas under Boolean connectives $\vee,\wedge, \neg$ and under sharply bounded quantifiers. For $i\geq 1$, the strict formula sets $\hat{\Sigma}^{b,\oplus_p P_i}$ are defined in the usual way by counting the number of alternations of bounded quantifiers.

$T^{1,\oplus_p P}_2$ is the theory axiomatized by the axioms for $PV_1$ symbols, the $C^k_p$ axioms for sharply bounded formulas $A(x)$, and $\hat{\Sigma}^{b,\oplus_p P}_1$-IND.

$APC_2^{\oplus_p P}:=T_2^{1,\oplus_p P}+dWPHP(PV_2^{\oplus_p P})$ where $PV_2^{\oplus_p P}$ means functions that can be computed in polynomial time relative to $NP^{\oplus_p P}$.

$\Sigma^b_{\infty}(\oplus_p)$ denotes formulas formed from bounded existential, universal, and $C_p$ quantifiers.

In $APC_2^{\oplus_p P}$, we say that a language is in $BP\cdot \oplus_p P$ if there exists $PV_1$ functions $f$ and $u$ such that for all $x$, $$x\in L\Leftrightarrow Pr_{r<u(x)}[f(x,r)\notin \oplus^1_p SAT]\preceq_0 1/4$$ $$x\notin L\Leftrightarrow Pr_{r<u(x)}[f(x,r)\notin \oplus^0_p SAT]\preceq_0 1/4$$ where $\oplus^i_p SAT$ is the set of propositional formulas $\phi$ such that the number of satisfying assignments of $\phi$ is congruent to $i$ mod $p$ for some prime $p$.

\begin{thm}[Buss, Ko\l odziejczyk, Zdanowski \cite{b2}] $APC_2^{\oplus_p P}$ proves that any $\Sigma^b_{\infty}(\oplus_p)$ formula defines a property in BP$\cdot\oplus_p$P.
\end{thm}

\subsection{Derandomization}

The approximate counting developed in $APC_1$ relies on a formalization of the derandomization result by Nisan and Wigderson \cite{nis}.

\begin{defi}[Je\v{r}\'{a}bek \cite{e}](in $APC_1$) A definable randomized algorithm is given by a pair of $PV$-functions $f,r$ such that
\medskip

$\exists w<r(x) \ f(x,w)\neq *\rightarrow Pr_{w<r(x)}[f(x,w)=*]\preceq_0 1/2$
\medskip

\noindent where $*$ is a special symbol signaling a rejecting computation.
\end{defi}

The special symbol $*$ could be avoided but it is useful for denoting a ``failure-state'' of probabilistic algorithms. It can be used when the input random string does not encode the expected structure, say a graph or a formula.

\begin{thm}[Je\v{r}\'{a}bek \cite{jer}]  Let $F$ be a randomized
  algorithm that is definable in $S^1_2+dWPHP(PV)$. Then there are $PV$-functions $h$ and $g$ such that $HARD^A$ proves $$\exists y\ y=F(x)\leftrightarrow h(x,\alpha(g(x)))\neq *$$ $$\exists y\ y=F(x)\rightarrow h(x,\alpha(g(x)))=F(x)$$
\end{thm}

Je\v{r}\'{a}bek \cite{jerp} formalized also Impagliazzo-Wigderson's \cite{imp} derandomization which draws the same conclusion assuming only worst-case hardness. This turned out to be much harder than the Nisan-Wigderson construction mainly because list decoding of error-correcting codes used in the construction requires several algebraic tools concerning finite fields.

\begin{thm}[Je\v{r}\'{a}bek \cite{jerp}] Let $F$ be a randomized
  algorithm that is definable in $S^1_2+dWPHP(PV)$, and let $\epsilon>0$. Then there are $PV$-functions $h$ and $g$ such that $HARD_{\epsilon}$ proves $$\exists y\ y=F(x)\leftrightarrow h(x,\alpha(g(x)))\neq *$$ $$\exists y\ y=F(x)\rightarrow h(x,\alpha(g(x)))=F(x)$$ Here, $HARD_{\epsilon}$ is defined as an extension of $S^1_2(\alpha)$, i.e. relativized $S^1_2$, by the following axioms: 
\begin{enumerate}[label=\arabic*.]
\item[1.] $\alpha(x):2^{||x||}\rightarrow 2$
\item[2.] $x\geq c\rightarrow Hard_{\epsilon}(\alpha(x))$
\end{enumerate} for a standard constant $c$.
\end{thm}

\subsection{Contribution of our paper: the PCP theorem and the \texorpdfstring{$(n,d,\lambda)$}{(n,d,lambda)}-graphs}

We add to the list of formalized results mentioned in previous sections formalizations of the exponential PCP theorem, the PCP theorem, and certain pseudorandom constructions involving the so called $(n,d,\lambda)$-graphs which are needed in the proof of the PCP theorem. The exponential PCP theorem was proved in Arora-Safra \cite{san}, and the PCP theorem is originally from Arora-Safra \cite{san} and Arora et.al. \cite{sansud}. In \cite{din} Dinur gave a simpler proof of the PCP theorem which we will formalize.

\begin{defi}\label{dal} (in $APC_1$)
Let $k,k',d$ be constants, $x\in \{0,1\}^n$ for $n\in Log$. Further, let $w\in \{0,1\}^{kn^k}$ (represent random bits), $\pi$ be a $k'n^{k'}$-size circuit with $m$ inputs where $m$ might differ from $n$, and $D$ be a $kn^k$-time algorithm.

Denote by $D^{\pi,w}(x)$ the output of $D$ on input $x$ and with access to $\pi$ specified by (random bits) $w$ as follows. $D$ computes $\pi$ on at most $d$ different inputs: first, it produces strings $\hat{w}_1,...,\hat{w}_d$ where each $\hat{w}_i\in\{0,1\}^m$, then it computes $\pi(\hat{w}_1),...,\pi(\hat{w}_d)$ and finally computes its output which is either 1 or 0. 
\end{defi}

We formulate the exponential PCP theorem in $APC_1$ as follows. For an explanation and a discussion concerning the choice of the formulation see Section \ref{epcp}.

\begin{theorems:exp}[The exponential PCP theorem in $APC_1$]\label{theorems:exp} There are constants $d,k,k'$ and a $kn^k$-time algorithm $D$ (given as a $PV$-function) computing as in Definition \ref{dal} such that $APC_1$ proves that for any $x\in\{0,1\}^n$, $n\in Log$: $$\exists y SAT(x,y)\rightarrow \exists k'n^{k'}size\ circuit\ \pi\ \forall w<2^{kn^k}, D^{\pi,w}(x)=1$$ $$\forall y \neg SAT(x,y)\rightarrow \forall k'n^{k'}size\ circuit\ \pi, Pr_{w<2^{kn^k}}[D^{\pi,w}(x)=1]\preceq_0 1/2$$
\end{theorems:exp}

We also formalize pseudorandom constructions involving the $(n,d,\lambda)$-graphs in $PV_1$ but leave the presentation of these results to Section \ref{rest} as it would require introducing too many definitions now.
\medskip

In order to formalize the PCP theorem we use the notion of probability $Pr$ on spaces of polynomial size $poly(n)$ for $n\in Log$ which is assumed to be defined in a natural way using an exact counting of sets of polynomial size which is also assumed to be defined in $PV_1$ in a standard way. The notion of probability $Pr$ should not be confused with the definition of $Pr$ in $APC_1$. We formulate (the more important implication of) the PCP theorem in $PV_1$ as follows.

\begin{defi}\label{dall}(in $PV_1$)
Let $k,c,d$ be constants, $x\in \{0,1\}^n, n\in Log, w\in \{0,1\}^{c\log n}$, $\pi\in \{0,1\}^{dn^c}$, and be $D$ be a $kn^k$-time algorithm.

Denote by $D^{\pi,w}(x)$ the output of $D$ on input $x$ and with access to $\pi$ specified by $w$ as follows. $D$ uses at most $c\log n$ random bits $w$ and makes at most $d$ nonadaptive queries to locations of $\pi$, i.e. $D$ can read bits $\pi_{i_1},...,\pi_{i_d}$ for $i_1,...,i_d$ produced by $D$. Then it computes its outputs, 1 or 0. 
\end{defi}

In Definition \ref{dall} we abuse the notation and use the shortcut $D^{\pi,w}(x)$ in different meaning than in Definition \ref{dal}. This should not lead into confusion.

\begin{theoremss:pcptt}[The PCP theorem in $PV_1$] There are constants $d,k,c$ and a $kn^k$-time algorithm $D$ (given as a $PV$-function) computing as in Definition \ref{dall} such that $PV_1$ proves that for any $x\in\{0,1\}^n, n\in Log$: $$\exists y SAT(x,y)\rightarrow \exists \pi\in\{0,1\}^{dn^c} \ \forall w<n^c, D^{\pi,w}(x)=1$$ $$\forall y \neg SAT(x,y)\rightarrow \forall \pi\in\{0,1\}^{dn^c}, Pr_{w<n^c}[D^{\pi,w}(x)=1]\leq 1/2$$
\end{theoremss:pcptt}

Note that the exponential PCP theorem follows from the PCP theorem. Hence, the exponential version is also provable in $PV_1$. The $PV_1$ proof of the PCP theorem uses (among many other tools) the exponential PCP theorem but scaled down to constant size instances so that to prove the scaled down version we need to reason only about sets of constant size. On the other hand, in $APC_1$ we perform a reasoning with p-time definable sets. Hence, the $APC_1$ proof shows different tools to be available in low fragments of arithmetic.

\section{The Cook-Levin theorem in \texorpdfstring{$PV_1$}{PV1}}\label{pv}

This section serves mainly as an illustration of some techniques available in $PV_1$ which we later use freely in our arguments.

\begin{thm}(The Cook-Levin theorem in $PV_1$) For each $k$, we have a $PV$-function $f$ such that $PV_1$ proves: for any $M,x$, $$\exists w,z; |z|, |w|\leq |x|^k, M(x,z,w)=1\leftrightarrow \exists y, |y|\leq 3|M||x|^{2k}, SAT(f(M,x),y))$$ where $M(x,z,w)=1$ is an open $PV$-formula which holds iff $w$ is an accepting computation of Turing machine $M$ on input $x,z$, and $|M|$ is the length of $M$'s code. 
\end{thm}

\proof 

\def\netreba{
Note firstly that $PV$ can introduce functions using conditional definitions:
\medskip

$f(x):=g_i(x)$ if $P_i(x), i=0,1$
\medskip

\noindent where $g_0,g_1$ are functions already defined in $PV$ and $P_0, P_1$ are disjoint and exhaustive open $PV$-formulas. This is because such $P_0,P_1$ define p-time truth functions which can be introduced as $PV$-functions $P'_0,P'_1$ and $f$ can be then defined as 

$f(x):=P'_0(x)g_0(x)+P'_1(x)g_1(x)$. 
\smallskip
}

First, we show that for some $PV$-function $f$, $PV_1$ proves ($*$): $$\forall M,x,z,w; |z|, |w|\leq |x|^k \exists y; |y|\leq 3|M||x|^{2k}\  (M(x,z,w)=1\rightarrow SAT(f(M,x),y))$$
The Turing machine $M$ is represented as a binary string encoding a tuple $(Q, \Sigma, b, F, \rho)$ where $Q$ is the set of states, $\Sigma$ is the set of tape symbols, $b\in Q$ is the initial state, $F\subseteq Q$ is the set of accepting states, and $\rho\subseteq ((Q-F)\times \Sigma)\times (Q\times \Sigma \times \{-1,1\})$ is the transition function.

We assume that the open $PV$-formulas $M(x,z,w)=1$ and $SAT(x,y)$ are already constructed in a well-behaved way.

The propositional formula $f(M,x)$ will be built from atoms
$T_{i,j,s}$ with intended interpretation ``tape cell $i$ of $M$
contains symbol $j$ at step $s$", atoms $H_{i,s}$ for ``$M$'s head is
at tape cell $i$ at step $s$", and atoms $Q_{q,s}$ for ``$M$ is in
state $q$ at step $s$". These atoms are assumed to be encoded in a
standard way.

Given $M,x$ we define $f(M,x)$ gradually by introducing more and more complex functions. This is supposed to illustrate the way in which $PV_1$ introduces new functions.

Let us start with a definition of function $f_{input}(x,y)$ mapping $x,y$ to a conjunction of $|y|$ atoms representing first $|y|$ bits of binary string $x$:
\[\eqalign{
  f_{input}(x,0)&:=0\cr
  f_{input}(x,s_i(y))&:={}'f_{input}(x,y)\wedge T_{|y|,i,0}\ '\
  \mbox{if}\ 
  |y|\leq |x|\wedge x_{|y|}=i\ \mbox{,}\ i= 0,1
  }
\]

\noindent 
where $'A\wedge B'$ 
is a code of the conjunction of propositional formulas encoded in $A$ and $B$. 

Next, put $f_{ins}(M,x):='f_{input}(x,x)\wedge Q_{b,0}\ '$.
\medskip

Then, define $f_{symb}(M,x,[t,l,m])='f_{ins}(M,x)\wedge G'$ where $G$ is a conjunction of formulas $(T_{t',l',m'}\rightarrow \neg T_{t',l'',m'})$ for all $l'\neq l''$ and $t',m'$ such that $[t',l',m'],[t',l'',m']\leq [t,l,m]$. This guarantees that cell $t'\leq t$ contains only one symbol at step $m'\leq m$.
\[\eqalign{
  f_{symb}(M,x,0)&:=f_{ins}(M,x)\cr
  f_{symb}(M,x,s_i([t,l,m]))&:={}'f_{symb}(M,x,[t,l,m])\wedge
  (T_{t',l',m'}\rightarrow \neg T_{t',l'',m'})\ '\cr
&\qquad\mbox{if}\ l'\neq l''\wedge [t',l',m'],[t',l'',m']\leq [t,l,m]\
\mbox{,}\ i\in \{0,1\}
  }
\]

Similarly, define $f_{state}(M,x,[t,l,m])$ by extending $f_{symb}(M,x,[t,l,m])$ with 
\begin{enumerate}[label=\arabic*.]
\item $Q_{t',m'}\rightarrow \neg Q_{t'',m'}$ for $t'\neq t''$ ($M$ cannot be in two different states at step $m'$)
\item $H_{t',m'}\rightarrow \neg H_{t'',m'}$ for $t'\neq t''$ (Head cannot be in two different positions at step $m'$)
\item $T_{t',l'',m'}\wedge T_{t',l',m'+1}\rightarrow H_{t',m'}$ for $l'\neq l''$ and $t',t''\leq t; l',l''\leq l;m'\leq m$
\end{enumerate}

Further, in this way introduce function $f_{trans}$ capturing $M$'s transition function $\rho$.
\[\eqalign{
  f_{trans}(M,x,c)&:={}'f_{state}(M,x,[|x|^k,|x|^k,|x|^k])\wedge\cr
&\qquad(H_{j,c}\wedge Q_{q,c}\wedge T_{j,\sigma, c}\rightarrow
\bigvee_{(q,\sigma,q',\sigma',d)\in \rho} (H_{j+d,c+1}\wedge
Q_{q',c+1}\wedge T_{j,\sigma',c+1})'
  }
\]
Finally, $f(M,x):='f_{trans}(M,x,|x|^k)\wedge \bigvee_{r\in F, t\leq |x|^k} Q_{r,t}\ '$.

This defines a $PV$-function $f$. To see that ($*$) holds, given $M,x,w$, we define $y$ assigning 0 or 1 to atoms of the formula $f(M,x)$ as follows:
\begin{enumerate}[label=\arabic*.]
\item $y(T_{j,i,0})=1$ iff $x_j=i$ for $i=0,1$ and $j<|x|$.
\item[\cW] $y(T_{j,i,t})=1$ iff $w$ says that tape cell $j$ of $M$ at step $t$ contains $i$


\item[2.] $y(H_{j,c})=1$ iff $w$ says that at step $c$ head is in position $j$


\item[3.] $y(Q_{r,t})=1$ iff $w$ contains $M$ in state $r$ at step $t$
\end{enumerate}

\smallskip

Informally, if $w$ indeed encodes an accepting computation of Turing machine $M$ on input $x,z$, then the previous definition produces $y$ which satisfies all conjuncts in formula $f(M,x)$ because these are copying the conditions from the definition of $M(x,z,w)=1$. Therefore, we can conclude that $M(x,z,w)=1\rightarrow SAT(f(M,x),y)$ in the theory $PV_1$.
\smallskip

Analogously, $PV_1\vdash \forall M,x,y,\exists w,z (SAT(f(M,x),y)\rightarrow M(x,z,w)=1)$.
 \qed

\section{The exponential PCP theorem in \texorpdfstring{$APC_1$}{APC1}}\label{epcp}


The exponential PCP theorem was proved in Arora-Safra \cite{san}. We formalize it in the theory $APC_1$ basically following the presentation in Arora-Barak \cite{ba}. However, there is a crucial change: we cannot use the Fourier transformation to derive the linearity test because it would require manipulations with exponentially big objects and it is not clear whether this could be done (for example, using a representation by circuits). Instead, we formalize the so called majority correction argument as it is presented in Moshkovitz \cite{mos}. Other parts of the proof work without much change. It is essential that all sets used to express probabilities are definable by p-size circuits so that $APC_1$ can work with them and the proof itself does not use more than basic operations on these sets which are available in $APC_1$. 







\smallskip

Recall Definition \ref{dal} introducing the predicate $D^{\pi,w}(x)$. The algorithm $D$ will represent the so called verifier of probabilistically checkable proofs $\pi$. The verifier is usually defined so that $\pi$ is allowed to be any string of arbitrary length and $D$ has an oracular access to $\pi$, it can ask for any bit of $\pi$. Then, for a language $L$, $L\in PCP(poly(n),1)$ standardly means that there is a p-time algorithm $D$ such that: 
\begin{itemize}
\item[1.] If $x\in L$, then there is a string $\pi$ (proof) such that $D$ with input $x$ of length $n$ and $poly(n)$ random bits asks for at most $O(1)$ bits of $\pi$ and accepts (with probability 1); 
\item[2.] If $x\notin L$, then for any $\pi$, $D$ with input $x$ of length $n$ and $poly(n)$ random bits asks for at most $O(1)$ bits of $\pi$ and accepts with probability $\leq 1/2$.
\end{itemize}

\noindent The exponential PCP theorem says that $NP\subseteq PCP(poly(n),1)$. As the verifier uses $poly(n)$ random bits, the proof $\pi$ can be seen as a string of size $2^{poly(n)}$. In our formalization, $n\in Log$ so bounded arithmetic cannot encode the exponentially big proofs by binary strings. In order to be able to speak about them we represent such proofs by p-size circuits. More precisely, for a $k'n^{k'}$-size circuit $\pi$ with $m$ inputs and $x\in \{0,1\}^m$, $\pi(x)$ is the $x$-th bit of the proof represented by $\pi$. Hence, the condition 1.) in our formulation of the exponential PCP theorem will look formally stronger but it follows trivially from the standard proof. In condition 2.) our $D$ will recognize errors only in proofs that are represented by $k'n^{k'}$-size circuits. We can interpret it as if the proofs that are not represented by such circuits were automatically rejected. Alternatively, we could also represent proofs by oracles which would maybe better reflect the nature of the exponential PCP theorem. However, then we would need to perform the formalization in the theory $APC_1$ extended by such oracles.
\smallskip

As the NP-completeness of SAT is provable in $PV_1$ it is sufficient to show in $APC_1$ that SAT $\in PCP(poly(n),1)$. This should justify Theorem \ref{theorems:exp} as the right formulation of the exponential PCP theorem in $APC_1$.


\proof (of Theorem \ref{theorems:exp})
For any $x\in\{0,1\}^n$, the algorithm $D$ firstly reduces SAT instance $x$ to a set of quadratic equations: It obtains 3SAT formula equivalent to $x$ by introducing new variable for each gate of the formula encoded in $x$ and clauses representing the gate. For each clause of the form $x_1\vee x_2\vee x_3$ it produces two equations $(1-x_1)y=0$ and $y-(1-x_2)(1-x_3)=0$ where $y$ is a new variable. Analogously for other possible clauses, if some $x_i$ occurs in the clause negatively, $1-x_i$ in the resulting equations is replaced by $x_i$. In this way $D$ produces a set of quadratic equations which is solvable in $F_2$ if and only if $x$ is satisfiable. More precisely, there is $k_0$ such that if $x$ encodes a propositional formula with $n_0$ variables it can be efficiently mapped to a set of $m\leq |x|^{k_0}$ quadratic equations on $n_1\leq |x|^{k_0}$ variables $u_1,...,u_{n_1}$ (w.l.o.g. $u_1=1$). The set of equations can be represented by an $m\times n^2_1$ matrix $A$ and a string $b\in\{0,1\}^m$ satisfying: $$\exists y\  SAT(x,y)\rightarrow \exists u\ Au\otimes u=b$$ $$\forall y\  \neg SAT(x,y)\rightarrow \forall u\ Au\otimes u\neq b$$ where $u\in\{0,1\}^{n_1}$ and $u\otimes u$ is a vector of bits $u_iu_j, i,j\in [n_1]$ ordered lexicographically. 
\smallskip

The algorithm $D$ will interpret $k'n^{k'}$-size circuits $\pi$ with $n_1^2+n_1+1$ inputs $b,z,z',$ where $b\in\{0,1\}, z\in\{0,1\}^{n_1}, z'\in\{0,1\}^{n^2_1}$, as circuits allowing us to access functions $f_{\pi}=WH(u)$ and $g_{\pi}=WH(u\otimes u)$ for some $u\in\{0,1\}^{n_1}$ in the following way, $\pi(0,z,z')=WH(u)(z)$ and $\pi(1,z,z')=WH(u\otimes u)(z')$. Here, $WH(u)(z):=\Sigma^{n_1}_{i=1}u_iz_i\  mod\ 2$. Similarly for $WH(u\otimes u)(z')$. $WH$ stands for ``Walsh-Hadamard''.
\smallskip

For any $x\in\{0,1\}^n $, the algorithm $D$ with $\leq kn^k$ random bits $w=r^l_1,...,r^l_7$ for $l=1,...,m_0$, where $m_0$ is a constant, $r^l_1,r^l_2,r^l_3\in \{0,1\}^{n_1}$, $r^l_4,r^l_5,r^l_6\in \{0,1\}^{n^2_1}, r^l_7\in \{0,1\}^m$  and with access to an $k'n^{k'}$-size circuit $\pi$ accepts if and only if for each $l=1,...,m_0$, $\pi$ passes the following tests

\begin{itemize}

\item ``linearity'': $f(r^l_1+r^l_2)=f(r^l_1)+f(r^l_2)$ and $g(r^l_4+r^l_5)=g(r^l_4)+g(r^l_5)$

\item ``$g_{\pi}$ encodes $u\otimes u$'': $g'(r^l_1\otimes r^l_2)=f'(r^l_1)f'(r^l_2)$

\item ``$g_{\pi}$ encodes a satisfying assignment'': $g'(z)=\Sigma^{m}_{i=1} (r^l_7)_ib_i$ for $z$ representing the sum $\Sigma^{m}_{i=1} (r^l_7)_i(A_iu\otimes u)$ where $A_iu\otimes u$ is the lefthand-side of the $i$-th equation in $Au\otimes u=b$
\end{itemize}

\noindent Here, $f=f_{\pi}, g=g_{\pi}$, $f'(r^l_1)=f(r^l_1+r^l_3)+f(r^l_3), f'(r^l_2)=f(r^l_2+r^l_3)+f(r^l_3)$ and similarly $ g'(r^l_1\otimes r^l_2)=g(r^l_1\otimes r^l_2+r^l_6)+g(r^l_6), g'(z)=g(z+r^l_6)+g(r^l_6)$.
\smallskip

For any $x\in\{0,1\}^n$, if $\exists y SAT(x,y)$ then there is $u\in\{0,1\}^{n_1}$ solving the corresponding equations $Au\otimes u=b$. Thus there is a $k'n^{k'}$-size circuit $\pi$ with $n^2_1+n_1+1$ inputs given by $\pi(0,z,z'):=WH(u)(z)$ and $\pi(1,z,z'):=WH(u\otimes u)(z')$ which passes all the tests: for any $w$, the linearity is clearly satisfied by the definition. Further:
\smallskip

\noindent $g'(r^l_1\otimes r^l_2)=g(r^l_1\otimes r^l_2+r^l_6)+g(r^l_6)=g(r^l_1\otimes r^l_2)=\Sigma^{n_1}_{i,j=1} u_iu_j(r^l_1)_i(r^l_2)_j$

\quad\quad\quad $= \Sigma^{n_1}_{i=1} u_i(r^l_1)_i\Sigma^{n_1}_{j=1} u_j(r^l_2)_j= f(r)f(r')=f'(r)f'(r')$
\smallskip

\noindent and as $Au\otimes u=b$ also $g'(z)=\Sigma^m_{i=1}(r^l_7)_ib_i$.
\smallskip

Now we will show that the algorithm $D$ recognizes incorrect proofs with high probability. The argument relies on the Test of linearity which we prove in Section \ref{ct}.

\begin{prop}[Test of linearity in $APC_1$]\label{lin}
Let $\epsilon$ be sufficiently small, $\epsilon^{-1}\in Log$ and let $f$ be a function on $n_1\in Log$ inputs represented by a circuit such that for each linear function $g$ with $n_1$ inputs,$$Pr_{x\in\{0,1\}^{n_1}}[f(x)=g(x)]_{\epsilon}< p$$ Then $Pr_{x,y}[f(x+y)=f(x)+f(y)]_{\epsilon}\preceq_{11\epsilon+13\epsilon^2+2\epsilon^3 } max\{29/32,1/2+p/2\}$. 
\smallskip

\noindent(We abuse the notation and use $f$ also in place of circuits representing $f$. Note that $g$ is represented by $n_1$ coefficients.)
\end{prop}

\begin{clm}[Local decoding in $APC_1$] Let $s<1/4, \epsilon\leq 1$ and $f$ be a function on $n_1\in Log$ inputs represented by a circuit such that there is a linear function $f_l$ which satisfies $Pr_{x<2^{n_1}}[f(x)=f_l(x)]_{\epsilon}\geq 1-s$. Then for each $x<2^{n_1}$,

 $Pr_{r<2^{n_1}}[f_l(x)=f(x+r)+f(r)]_{\epsilon}\succeq_{6\epsilon} 1-2s$.
\end{clm}

\proof (of the claim) By the assumption and Proposition \ref{really} 1.i), for $x<2^{n_1}$, \\ $\{r| f(r)\neq f_l(r)\}\cap 2^{n_1}\preceq_{2\epsilon} s2^{n_1}$ and $\{r|f(x+r)\neq f_l(x+r)\}\cap 2^{n_1}\preceq_{2\epsilon} s2^{n_1}$ which implies $\{r|f(r)\neq f_l(r)\vee f(x+r)\neq f_l(x+r)\}\cap 2^{n_1}\preceq_{4\epsilon} 2s2^{n_1}$. By linearity of $f_l$, for any $x<2^{n_1}$, $\{r| f_l(x)\neq f(x+r)+f(r)\}\subseteq \{r| f_l(r)\neq f(r)\vee f_l(x+r)\neq f(x+r)\}$. 

Thus, $Pr_{r}[f_l(x)=f(x+r)+f(r)]_{\epsilon}\succeq_{6\epsilon} 1-2s$, which proves the claim.\qed

Assume that $\forall y\neg SAT(x,y)$, so $\forall u, Au\otimes u\neq b$ and let $\pi$ be arbitrary circuit of size $k'n^{k'}$. Further, let $\epsilon$ be sufficiently small, $\epsilon^{-1}\in Log$ and denote by $D^{\pi,w}_1(x)$, $D^{\pi,w}(x)$ with $m_0=1$, i.e. $D$ performing only one round of testing.

If for each linear function $g_l$, $Pr_{x\in \{0,1\}^{n^2_1}}[g(x)=g_l(x)]_{\epsilon}<31/32$ or for each linear function $f_l$, $Pr_{x\in \{0,1\}^{n_1}}[f(x)=f_l(x)]_{\epsilon}<31/32$, then by the test of linearity, we have $Pr_{w}[D^{\pi,w}_1(x)=1]_{\epsilon}\preceq_{13\epsilon+13\epsilon^2+2\epsilon^3} 63/64$. Otherwise, there are linear functions $g_l$, $f_l$ such that by local decoding, for each $x\in \{0,1\}^{n^2_1}$, it holds $Pr_r[g_l(x)=g'(x)]_{\epsilon}\succeq_{6\epsilon} 15/16$ where $g'(x)=g(x+r)+g(r)$ and for each $x\in\{0,1\}^{n_1}$, $Pr_r[f_l(x)=f'(x)]_{\epsilon}\succeq_{6\epsilon} 15/16$ where $f'(x)=f(x+r)+f(r)$.

We need to show that even in the latter situation verifier $D$ accepts with small probabilty. For this, we distinguish two cases: 1. $g_l\neq WH(u\otimes u)$, i.e. $\exists x,y,\ g_l(x\otimes y)\neq f_l(x)f_l(y)$; 2. $g_l=WH(u\otimes u)$. Here, by the linearity of $f_l$, we have $f_l=WH(u)$ for some $u$ and $f_lf_l=WH(u\otimes u)$.
\medskip

\begin{clm}\label{guu} If $g_l\neq WH(u\otimes u)$, then  $Pr_{r_1,r_2}[g_l(r_1\otimes r_2)\neq f_l(r_1)f_l(r_2)]\succeq_{2\epsilon} 1/4$
\end{clm}

\proof Let $U,W$ be matrices such that $g_l(x\otimes y)=xUy$ and $f_l(x)f_l(y)=xWy$.

If $U\neq W$, then $\{r_2\in 2^{n_1}|Ur_2\neq Wr_2\}\succeq_{0} 2^{n_1}/2$ as witnessed by the following circuit: Let $(i,j)$ be a position where $U$ and $W$ differ. Consider the circuit mapping $r_2$ from $\{r_2\in 2^{n_1}| Ur_2\neq Wr_2\}$ to $\hat{r_2}$ where $\hat{r_2}<2^{n_1}/2$ is obtained from $r_2$ by erasing its $j$th bit $(r_2)_j$. For each $r_2<2^{n_1}/2$, let $r^0_2<2^n$ be such that $r_2=\hat{r^0_2}$ and $(r^0_2)_j=0$ and let $r^1_2<2^{n_1}$ be such that $r_2=\hat{r^1_2}$ and $(r^1_2)_j=1$. Then, for each $r_2<2^n/2$, $r^0_2$ or $r^1_2$ is in $\{r_2\in 2^{n_1}|Ur_2\neq Wr_2\}$.

Furthermore, if $U\neq W$, we similarly observe that $\{r_1\in 2^{n_1}| r_1Ur_2\neq r_1Wr_2\}\succeq_0 2^n/2$ for each $r_2<2^{n_1}$. Hence, by Proposition \ref{really} $2.$, $\{\left<r_1,r_2\right>|g_l(r_1\otimes r_2)\neq f_l(r_1)f_l(r_2) \}\succeq_{\epsilon} 2^{2n}/4$. This proves the claim.\qed

Suppose now that $g_l\neq WH(u\otimes u)$. As $\{\left<r_1,r_2\right>|g'(r_1\otimes r_2)=f'(r_1)f'(r_2)\}$ is a subset of 
\[\eqalign{
\bigl\{\langle r_1,r_2\rangle\,|\,g'(r_1\otimes r_2)=g_l(r_1\otimes r_2)&\wedge
g_l(r_1\otimes r_2)=f_l(r_1)f_l(r_2)\cr
&\wedge  f'(r_1)=f_l(r_1) \wedge f'(r_2)=f_l(r_2)\bigr\}\cr
\cup\,\bigl\{\langle r_1,r_2\rangle\,|\,g'(r_1\otimes r_2)\neq g_l(r_1\otimes r_2) &\vee f'(r_1)\neq f_l(r_1) \vee f'(r_2)\neq f_l(r_2)\bigr\}
}\]
 which is $\preceq_{28\epsilon} 15/16 (2^{2n_1})$ by Claim \ref{guu}, we can conclude that
\[Pr_w[D^{\pi,w}_1(x)=1]_{\epsilon}\preceq_{2\epsilon} Pr_{r_1,r_2}[g'(r_1\otimes r_2)=f'(r_1)f'(r_2)]_{\epsilon}\preceq_{28\epsilon} 15/16.\] 

\noindent It remains to consider the case that $g_l=WH(u\otimes u)$. 
For each $u<2^{2n_1}$, $R=\{r|\Sigma_i r_i(A_iu\otimes u)\neq \Sigma_i r_ib_i\}\cap 2^m\succeq_0 2^m/2$ as it is witnessed by the following circuit. Let $j$ be the first such that $A_ju\otimes u\neq b_j$. The circuit maps each $r\in R$ to $\hat{r}$ where $\hat{r}<2^m/2$ is obtained from $r$ by erasing its $j$th bit $r_j$. For each $r<2^m/2$, let $r^0<2^m$ be such that $r=\hat{r^0}$ and $r^0_j=0$ and let $r^1<2^m$ be such that $r=\hat{r^1}$ and $r^1_j=1$. Then, for each $r<2^m/2$, $r^0\in R$ or otherwise $\Sigma_i r^0_i(A_iu\otimes u)=\Sigma_i r^0_ib_i$ and hence $r^1\in R$.
\smallskip

Furthermore, assuming $g_l=WH(u\otimes u)$, $\{r|g'(z)=\Sigma_i
r_ib_i\}$ is a subset of
\[\{r|\Sigma r_i(A_iu\otimes u)= \Sigma_i r_ib_i \wedge g_l(z)=g'(z)\}\cup \{r|g_l(z)\neq g'(z)\}\]
Thus, $Pr_{w}[D^{\pi,w}_1(x)=1]_{\epsilon}\preceq_{2\epsilon} Pr_r[g'(z)=\Sigma_i r_ib_i]_{\epsilon}\preceq_{10\epsilon} 9/16$.
\smallskip

In all cases, $Pr_w[D^{\pi,w}_1(x)=1]_{\epsilon}\preceq_{28\epsilon} 63/64$ so $$\{w\in 2^{3n_1+n_1^2+m}|D^{\pi,w}_1(x)=0\}\succeq_{30\epsilon} 1/64 (2^{3n_1+n_1^2+m})$$

Therofore, for sufficiently big constant $m_0$, Chernoff's bound from Proposition \ref{really} with $\delta^2:=c30\epsilon+1/100^2$ and sufficiently small $\epsilon$ implies that $Pr_{w<2^{kn^k}}[D^{\pi,w}(x)=1]\preceq_0 1/2$. 
\smallskip

To conclude the proof of the exponential PCP theorem in $APC_1$ it thus remains to derive the Test of linearity.

\subsection{Test of linearity in \texorpdfstring{$APC_1$}{APC1}}\label{ct}

In this section we prove Proposition \ref{lin} in the theory $APC_1$. 
\medskip

We cannot use the Fourier transformation argument directly as in Arora-Barak \cite{ba} which would require to prove the existence of exponentially long Fourier expansions (and it is not clear if this could be managed, for example, using a representation by p-size circuits). Instead we formalize the so called majority correction argument. Our presentation is a minor modification of Moshkovitz \cite{mos}. 
\smallskip

Let $\epsilon>0$ be sufficiently small and $\epsilon^{-1}\in Log$. Define $g_{\epsilon}: 2^n\mapsto 2$ by $$g_{\epsilon}(x)=1 \quad\quad \equiv_{def} \quad\quad Pr_{y<2^n}[f(y)+f(x+y)=1]_{\epsilon}\geq 1/2$$

Therefore, for any $x<2^n$, $P_x:=Pr_{y<2^x}[g_{\epsilon}(x)=f(y)+f(x+y)]_{\epsilon}\geq 1/2$. Hence, $g_{\epsilon}(x)$ is the majority value of the expression $f(y)+f(x+y)$ for possible $y$'s.

We will now derive three claims that can be combined into a proof of Proposition \ref{lin}.

\begin{clm}\label{CL:b-one}
 $Pr_{\left<x,y\right>}[f(x+y)\neq f(x)+f(y)]_{\epsilon}\succeq_{8\epsilon+13\epsilon^2+2\epsilon^3} \frac{1}{2} Pr_{x}[f(x)\neq g_{\epsilon}(x)]_{\epsilon}$ 
\end{clm}

\proof
This holds trivially if $Size(\{x|g_{\epsilon}(x)\neq f(x)\}\cap 2^n,\epsilon)=0$. Otherwise,
define sets 
\[T:=\{\left<x,y\right>| f(x+y)\neq f(x)+f(y)\}\ \mbox{and}\
G:=\{x|g_{\epsilon}(x)\neq f(x)\}.\]
 Then,
\[\eqalign{
  &Pr_{x<2^n,y<2^n}[f(x+y)\neq f(x)+f(y)]_{\epsilon}\cr
\geq{}& Size(T\cap (G\times 2^n)\cap 2^{2n},\epsilon)/2^{2n}\cr
={}&
\frac{Size((G\cap 2^n)\times 2^n,\epsilon)}{2^{2n}}\cdot
\frac{Size(T\cap (G\times 2^n)\cap 2^{2n},\epsilon)}{Size((G\cap
  2^n)\times 2^n,\epsilon)}
 }
\]
By Proposition \ref{lem} iii), $(G\cap 2^n)\times 2^n\approx_{\epsilon} Size(G\cap 2^n,\epsilon) 2^n$, so the first fraction in the expression above is $\approx_{2\epsilon} Pr_{x<2^n}[g_{\epsilon}(x)\neq f(x)]_{\epsilon}$.

Further, for each $x\in G\cap 2^n$, $P_x\geq 1/2$ and in particular, $2^n/2\preceq_{\epsilon} T_x=\{y| \left<x,y\right>\in T\}$.
Hence, by Proposition \ref{really} 2., $Size(G,\epsilon) 2^n/2\preceq_{3\epsilon+\epsilon^2} T\cap (G\times 2^n)$, and $$\frac{Size(T\cap (G\times 2^n)\cap 2^{2n},\epsilon)}{Size((G\cap 2^n)\times 2^n,\epsilon)}\succeq_{4\epsilon+\epsilon^2} \frac{Size(G,\epsilon)2^n}{2Size((G\cap 2^n)\times 2^n,\epsilon)} \succeq_{2\epsilon} 1/2$$ Applying now Proposition \ref{lem} iii) we obtain Claim \ref{CL:b-one}.\qed

\begin{clm}\label{CL:b-two}
If $Pr_{\left<x,y\right>}[f(x+y)\neq f(x)+f(y)]_{\epsilon}<\frac{3}{32}$, then $\forall x<2^n$, $P_x>\frac{3}{4}$.
\end{clm}

\proof
Fix $x<2^n$ and define 
\[\eqalign{
A&:=\{\left<y,z\right> | g_{\epsilon}(x)=f(y)+f(x+y)\wedge g_{\epsilon}(x)=f(x+z)+f(z)\}\cr
B&:=\{\left<y,z\right> | g_{\epsilon}(x)\neq f(y)+f(x+y)\wedge
g_{\epsilon}(x)\neq f(x+z)+f(z)\}
}\]
Then, $Pr_{y,z}[f(y)+f(x+y)=f(z)+f(x+z)]_{\epsilon}= Pr_{y,z}[\left<y,z\right>\in A\cup B]_{\epsilon}$.

By \ref{really} 1.ii), $(A\cup B)\cap 2^{2n}=(A\cap 2^{2n})\cup (B\cap 2^{2n})\approx_{3\epsilon} Size(A\cap 2^{2n},\epsilon)+Size(B\cap 2^{2n},\epsilon)$.
Thus, $Pr_{y,z}[\left<y,z\right>\in A\cup B]_{\epsilon}\approx_{4\epsilon} Pr_{y,z}[\left<y,z\right>\in A]+Pr_{y,z}[\left<y,z\right>\in B]$.

Next, let $A':=\{y|g_{\epsilon}(x)=f(x+y)+f(x)\}$. Using Proposition \ref{lem} iii) twice, $A\cap 2^{2n}$ is $(A'\cap 2^n)\times (A'\cap 2^n)\approx_{2\epsilon} Size(A'\cap 2^n,\epsilon) Size(A'\cap 2^n,\epsilon)$. Therefore, $Pr_{y,z}[\left< y,z\right>\in A]\approx_{3\epsilon} P_x^2$.  

As by Proposition  \ref{really} 1.i), $\{y|g_{\epsilon}(x)\neq f(x+y)+f(y)\}\cap 2^n=2^n-A'\cap 2^n$ is $\approx_{2\epsilon} 2^n-Size(A'\cap 2^n,\epsilon)$, we analogously obtain $Pr_{y,z}[\left<y,z\right>\in B]\approx_{9\epsilon} (1-P_x)^2$.
Therefore, $Pr_{y,z}[f(y)+f(y+x)=f(z)+f(x+z)]\approx_{17\epsilon} P_x^2+(1-P_x)^2$.

Define now, 
\[\eqalign{
C&:=\{\left< y,z\right> |f(y+z)\neq f(y)+f(z)\}\cr
D&:=\{\left< y,z\right> |f(y+z)\neq f(x+y)+f(x+z)\}
}\]
Then, $2^{2n}-(C\cap 2^{2n})\cup (D\cap 2^{2n}) \subseteq (A\cup B)\cap 2^{2n}$ and by Proposition \ref{lem} i) we have $2^{2n}-(C\cap 2^{2n})\cup (D\cap 2^{2n})\preceq_0 (A\cup B)\cap 2^{2n}$. 

By Proposition \ref{really} 1.ii), $(C\cap 2^{2n})\cup (D\cap 2^{2n})\preceq_{3\epsilon} Size(C\cap 2^{2n},\epsilon)+Size(D\cap 2^{2n},\epsilon)$, so $2^{2n}- Size(C\cap 2^{2n},\epsilon)- Size(D\cap 2^{2n},\epsilon)\preceq_{4\epsilon}2^{2n}-(C\cap 2^{2n})\cup (D\cap 2^{2n})$.

Moreover, by the assumption, $Pr_{y,z}[f(y)+f(z)\neq f(y+z)]_{\epsilon}<3/32$ and similarly, $Pr_{y,z}[f(y+z)\neq f(x+y)+f(x+z)]_{\epsilon}< 3/32$. Therefore, $$Pr_{y,z}[f(y)+f(x+y)=f(z)+f(x+z)]_{\epsilon}\succeq_{5\epsilon}13/16$$
This shows that $P_x^2+(1-P_x)^2\succeq_{22\epsilon}\frac{13}{16}$ and $2(P_x-\frac{1}{4})(P_x-\frac{3}{4})+\frac{10}{16}\succeq_{22\epsilon} \frac{13}{16}$. As $P_x\geq 1/2$, $P_x< 3/4$ would imply $\frac{10}{16}2^n\succeq_{22\epsilon}\frac{13}{16}2^n$ contradicting dual weak pigeonhole principle. Hence, Claim \ref{CL:b-two} follows.\qed

\begin{clm}\label{CL:b-three}
If $Pr_{x,y}[f(x+y)\neq f(x)+f(y)]_{\epsilon}<3/32$, then $g_{\epsilon}$ is linear.
\end{clm}

\proof
By Claim \ref{CL:b-two}, $\forall x,y<2^n$, 
\[\eqalign{
Pr_z[g_{\epsilon}(x)\neq f(x+z)+f(z)]_{\epsilon}&\preceq_{3\epsilon} 1/4\cr
Pr_z[g_{\epsilon}(y)\neq f(y+z)+f(z)]_{\epsilon}&\preceq_{3\epsilon} 1/4\cr
Pr_z[g_{\epsilon}(x+y)\neq
f(y+z)+f(z+x)]_{\epsilon}&\preceq_{3\epsilon} 1/4
}\]
Therefore, 
\[\eqalign{Pr_z[g_{\epsilon}(x)=f(x+z)+f(z)&\wedge g_{\epsilon}(y)=
  f(y+z)+f(z)\cr
& \wedge g_{\epsilon}(x+y)=
f(y+z)+f(z+x)]_{\epsilon}\succeq_{16\epsilon} 1/4 
}\]
The last estimation implies that if $\epsilon$ is sufficiently small, there exists $z_0$ (and we can efficiently find it) such that
\[\eqalign{
g_{\epsilon}(x)&=f(x+z_0)+f(z_0)\cr
g_{\epsilon}(y)&=f(y+z_0)+f(z_0)\cr 
g_{\epsilon}(x+y)&=f(y+z_0)+f(z_0+x)
}\]
which shows that $g_{\epsilon}(x)+g_{\epsilon}(y)=g_{\epsilon}(x+y)$ and proves Claim \ref{CL:b-three}.\qed

We can now derive Proposition \ref{lin}. Assume that for each linear function $g$ we have $Pr_x[g(x)=f(x)]_{\epsilon}<p$. By Claim \ref{CL:b-three}, $Pr_{x,y}[f(x+y)\neq f(x)+f(y)]_{\epsilon}\geq 3/32$ or $g_{\epsilon}$ is linear. This means that either $Pr_{x,y}[f(x+y)=f(x)+f(y)]_{\epsilon}\preceq_{3\epsilon}29/32$ or $Pr_x[g_{\epsilon}(x)=f(x)]<p$. In the latter case, $Pr_x[g_{\epsilon}(x)\neq f(x)]\succeq_{3\epsilon}1-p$ and by Claim \ref{CL:b-one},\\ $Pr_{x,y}[f(x+y)=f(x)+f(y)]_{\epsilon}\preceq_{11\epsilon+13\epsilon^2+2\epsilon^3} 1/2+p/2$.

\section{Pseudorandom constructions in \texorpdfstring{$PV_1$}{PV1}}\label{rest}

In order to derive the PCP theorem in $PV_1$ we will need to prove in the theory $PV_1$ the existence and some properties of the $(n,d,\lambda)$-graphs (see their definition below). While the construction itself is very combinatorial, its analysis uses algebraic techniques, e.g. properties of eigenvectors, which we do not know how to formalizable in $PV_1$. 

Using an equivalent combinatorial definition of the $(n,d,\lambda)$-graphs it is possible to derive their existence and main properties by only combinatorial tools. However, we need it for the algebraic equivalent and the implication producing the algebraic $(n,d,\lambda)$-graphs from the combinatorial $(n,d,\lambda)$-graphs is one of those which seem to require the algebraic techniques we are trying to avoid.

Therefore, we will employ an approximation of some algebraic tools which will allows us to derive slightly weaker results about the algebraic $(n,d,\lambda)$-graphs that are, however, sufficient to derive the PCP theorem. 

For the history of the field leading to the results presented in this section see Arora-Barak \cite[Chapter 21]{ba}.

\subsection{Definition and some properties of the \texorpdfstring{$(n,d,\lambda)$}{(n,d,lambda)}-graphs}

In $PV_1$ we say that a graph $G$ is $d$-regular if each vertex appears in exactly $d$ edges. We allow $G$ to have multiple edges and self-loops. The random-walk $n\times n$ matrix $A$ of a $d$-regular graph $G$ with $n$ vertices consists of elements $A_{i,j}$ being the number of edges between the $i$-th and the $j$-th vertex in $G$ divided by $d$. All our graphs will be undirected, hence, their random-walk matrices will be symmetric. For any $k$ and a graph $G$ with $n$ vertices, we denote by $G^k$ the graph with $n$ vertices which has an edge between the $i$th and the $j$th vertex for each $k$ step path between the $i$th and the $j$th vertex in $G$.
\smallskip

We would like to define now the second largest eigenvalue of $G$ denoted as $\lambda(G)$. The parameter $\lambda(G)$ corresponds to a certain expansion property of $G$ (see Proposition \ref{expans}) and normally it is defined as the maximum value of $||Ax||$ over all vectors $x$ in $n$-dimensional real vector space such that $||x||=1$ and $\Sigma_i x_i=0$. Here, $||y||=(\Sigma_i y_i^2)^{1/2}$ and $A$ is the random-walk matrix of graph $G$ with $n$ vertices. In $PV_1$ we will approximate this definition using a sufficiently dense net of rational numbers. 
\smallskip

The theory $PV_1$ proves that each $x$ is the value of an expression of the form $\Sigma^{|x|}_{i=0} 2^iy_i$ for $y_i\in \{0,1\}$ which is encoded in a natural way. In $PV_1$ we write that $x\in Q^n/m$ if $x=(x_1,...,x_n)$ and each $x_i$ is $\frac{a}{b}$ or $-\frac{a}{b}$ for $a\in [m]\cup \{0\}, b\in [m]=\{1,...,m\}$ where $a,b$ are represented by products of such expressions $\Sigma_i 2^iy_i, y_i\in \{0,1\}$. These products are also encoded in a natural way. In such cases we might write $a=c\cdot d$ to specify that $a$ is represented by a product of $c$ and $d$ where $c,d$ might be products of other expressions of the form $\Sigma_i 2^i y_i$.  

Let $L$ be a sufficiently big constant, then $SQRT$ is a function which given nonnegative $r\in Q/m$, $m>1$, produces $SQRT(r)\in Q/(Lm)^7$ such that $$0\leq (SQRT(r))^2-r\leq \frac{1}{L}$$ where we ignore the difference between $SQRT(r)$ and the value of the expression it encodes. Moreover, $SQRT$ satisfies the following: If input $r$ is a fraction of the form $\frac{c\cdot c\cdot e}{d\cdot d\cdot f}\in Q/m$ where $c,d$ are sums $\Sigma_i 2^iy_i$ with $y_i\in \{0,1\}$ (and $e,f$ might be products of such sums), then $$SQRT(\frac{c\cdot c\cdot e}{d \cdot d\cdot f})=\frac{c}{d}\cdot SQRT(\frac{e}{f})\ \ \ \ \ \ \ \ \ (*)$$ which is illustrating the representation of the number encoded in $SQRT(r)$. The representation of $\frac{c^2e}{d^2f}$ guarantees that $SQRT$ does not need to perform factorization.

The function $SQRT$ is essentially the usual algorithm approximating square root by a digit-by-digit search. We will assume that $SQRT$ works as follows: given $r\in Q/m$, it first finds out maximal $e,f\in [m]$ such that the current representation of $r$ is $\frac{e\cdot e}{f\cdot f}\frac{p}{q}$ for some $p,q\in [m]$, and then by a digit-by-digit search it finds the first $c\in [L^7m^6]$ such that $SQRT(r)$ which is $\frac{ec}{2fLqm^4}\in Q/(Lm)^7$ satisfies $0\leq (\frac{ec}{2fLqm^4})^2-r\leq \frac{1}{L}$. To get such $c$ we want to satisfy $c^2-4pqL^2m^8\leq 4Lm^6$. Thus $c\leq 2\sqrt{pq}Lm^4+2\sqrt{L}m^3\leq 7m^6$. The value $c$ is then produced by a p-time algorithm approximating $2\sqrt{pq}Lm^4$ so it is unique and its existence is provable in $PV_1$.
\medskip

For $x\in Q^n/m$, put $||x||:=SQRT(\Sigma_i x_i^2)$ where the input 
$\Sigma_i x_i^2\in Q/(nm^{2n})$ 
is computed so that if each $x_i=\pm\frac{a_ic}{b_id}$ for some common $c,d$, then $\Sigma_i x^2_i$ is 
represented as $\frac{e\cdot c\cdot c}{f\cdot d\cdot d}$ for some $e,f$. 

By the definition, if $x\in Q^n/m$, $x\neq 0$, then $\frac{x}{||x||}\in Q^n/((Lnm^{2n})^7m)$ and using $(*)$, $||\frac{x}{||x||}||=1$. Note that $||x||$ might be a fraction so we assume that $\frac{x}{||x||}$ is rearranged appropriately. 

However, by $||x||^2$ we always mean $\left<x,x\right>$ where $\left<x,y\right>:=\Sigma_i x_iy_i$ for $x,y\in Q/m$. The $n$-dimensional unite vector is defined as $\mathbf{1}:=(1/n,...,1/n)$. 
\smallskip

The parameter $\lambda(G)$ is defined as the maximum value of $||Ax||$ over all possible vectors $x\in Q^n/(Ln)^{(Ln)^L}$ such that $||x||=1$ and $\left<x,\mathbf{1}\right>=0$. Here again, 
the vector $Ax\in Q^n/(n(d(Ln)^{(Ln)^L})^n)$ (with elements of length $poly(n)$) is computed so that if each $x_i=\pm\frac{a_ic}{b_id}$ for some common $c,d$, then $(Ax)_j=\pm\frac{c\cdot e_j}{d\cdot f_j}$ for some $e_j,f_j$. 


\smallskip

We will not need to prove $\exists y,\ y=\lambda(G)$ in $PV_1$ but we will work with formulas of the form $\lambda(G)\leq y$ which are $\Pi^b_1$. To see this note that in $\lambda(G)\leq y$ we universally quantify over all $x$'s in $Q^n/(Ln)^{(Ln)^L}$. For each $j$, there are $\leq m^j$ ways how to represent $b\in [m]$ as a product of $j$ numbers 
so this is a universal quantification over $\leq 2^{n^{O(1)}}$ $x$'s. For each such $x$, predicates $||x||=1$ and $||Ax||\leq y$ are computable in time $n^{O(1)}$.

\begin{defi}
A $d$-regular graph $G$ with $n$ vertices is $(n,d,\lambda)$-graph if $\lambda(G)\leq \lambda<1$.
\end{defi}

We will often use Cauchy-Schwarz inequality in $PV_1$ which can be obtained in the standard way.

\begin{prop}(Cauchy-Schwarz inequality in $PV_1$) For every $n,m$ and $x,y\in Q^n/m$, $\left<x,y\right>^2\leq ||x||^2\cdot ||y||^2$ and therefore, if $n\in Log$ (and thus $||x||$ exists), also $\left<x,y\right>\leq ||x||\cdot ||y||$.  
\end{prop}

\proof If $y=0$, the inequality holds. Otherwise, let $z:=x-\frac{\left<x,y\right>}{\left<y,y\right>}y$. Then, $\left<z,y\right>=\left<x,y\right>-\frac{\left<x,y\right>}{\left<y,y\right>}\left<y,y\right>=0$. Therefore, $||x||^2=(\frac{\left<x,y\right>}{\left<y,y\right>})^2||y||^2+||z||^2=\frac{\left<x,y\right>^2}{||y||^2}+||z||^2\geq \frac{\left<x,y\right>^2}{||y||^2}$. \qed 

\smallskip

In Peano Arithmetic, regular graphs $G$ satisfy $\lambda(G)\leq 1$ but in $PV_1$ we will have just $\lambda(G)\leq 1+\epsilon+1/L$ for any rational $\epsilon>0$. Fortunately, this is enough to derive the PCP theorem in $PV_1$.

\begin{prop}\label{eps} For any $d$ and any rational $\epsilon >0$, $PV_1$ proves that for any $d$-regular graph $G$ with $n\in Log$ vertices, $\lambda (G)< 1+\epsilon+1/L$. 
\end{prop}

\proof As the statement we want to prove is $\forall\Sigma^b_1$, by $\forall\Sigma^b_1$-conservativity of $S^1_2$ over $PV_1$, we can work in the theory $S^1_2$. 

Let $A$ be the random-walk matrix of $G$. We want to show that $\lambda (G)<1+\epsilon+1/L$. Using Cauchy-Schwarz inequality, for every $x\in Q^n/(Ln)^{(Ln)^L}$ such that $||x||=1$, $$||Ax||^2=\Sigma_i (\Sigma_j A_{i,j}x_j)^2\leq \Sigma_i (\Sigma_j A_{i,j}^2\Sigma_j x_j^2)\leq\Sigma_i \Sigma_j A_{i,j}^2\leq \Sigma_i\Sigma_j A_{i,j}= \Sigma_i 1=n$$ As $A_{i,j}=A_{j,i}$, we have $\left<x,Ay\right>=\Sigma_i (x_i\Sigma_j A_{i,j}y_j)=\Sigma_j(y_j\Sigma_i(x_iA_{j,i}))=\left<Ax,y\right>$ and\\ $||Ax||^4=\left<Ax,Ax\right>^2=\left<A^2x,x\right>^2\leq ||A^2x||^2$ where $A^2$ is the random-walk matrix of $G^2$, so also $||A^2x||^2\leq n$ and $||Ax||^4\leq n$. This shows that $$\forall k\leq  K\log \log n\ (\forall A, ||Ax||^2\leq n^{1/(2^k)}\rightarrow \forall A, ||Ax||^2\leq n^{1/(2^{k+1})})$$ where $K$ is a sufficiently big constant depending only on $\epsilon$ and the universal quantifier before $A$ goes only over random-walk matrices of $d$-regular graphs with $n$ vertices. Note also that $n^{1/(2^k)}$ might be irrational but we can assume that it is approximated with a sufficiently small constant error so that the predicate $||Ax||^2\leq n^{1/(2^k)}$ is $\Pi^b_1$. 

Then, by $\Pi^b_1$-LLIND (available in $S^1_2$), we have $\forall A, ||Ax||^2\leq n^{1/(\log n)^K}$ which is $<(1+\epsilon)^2$ by the choice of $K$ and therefore $||Ax||\leq 1+\epsilon +1/L$.
 \qed 

We can now prove that the $(n,d,\lambda)$-graphs satisfy a useful expansion property. The term $\frac{\lambda d}{Ln^2}$ occuring in its formulation is an error resulting from our approximations in $PV_1$.

\begin{prop}(in $PV_1$)\label{expans} If $G$ is $(n,d,\lambda)$-graph with $n\in Log$ vertices $V$ and edges $E$, then for every $S\subseteq V, |S|\leq n/2$, $$|E(S,V-S)|\geq \frac{d|S|(1-\lambda)}{2}-\frac{\lambda d}{Ln^2}$$ where $E(S,T)$ denotes the set of edges $(i,j)\in E$ with $i\in S, j\in T$. 
\end{prop}

\proof It suffices to show: $$|E(S,V-S)|\geq (1-\lambda)\frac{d|S||V-S|}{n}-\frac{\lambda}{Ln^2}$$ Let $x\in Q^n/n$ be the following vector: $x_i=|V-S|$ if $i\in S$ and $x_i=-|S|$ if $i\in V-S$. Put $Z:=\Sigma_{i,j} A_{i,j}(x_i-x_j)^2$ for the random-walk matrix $A$ of $G$. Then,\\ $Z=\frac{2}{d}|E(S,V-S)|(|S|+|V-S|)^2$. As $A$'s rows and columns sum up to one, we have also $$Z=\Sigma_{i,j} A_{i,j}x^2_i-2\Sigma_{i,j}A_{i,j}x_ix_j+\Sigma_{i,j} A_{i,j}x^2_j=2||x||^2-2\left<x,Ax\right>$$ Further, $\Sigma x_i=0$ and $\frac{x}{||x||}\in Q^n/((Lnn^{2n})^7n)$ so $||Ax||= ||A\frac{x}{||x||}||||x||\leq \lambda ||x||$. By Cauchy-Schwarz inequality, $\left<x,Ax\right>\leq ||x||\cdot ||Ax||$. Therefore, $$\frac{1}{d}|E(S,V-S)|(|S|+|V-S|)^2\geq (1-\lambda)||x||^2-\lambda/L$$ It remains to observe that $||x||^2=|S||V-S|(|S|+|V-S|)$  \qed

In the following proposition we use the notion of probability $Pr$ on sets of polynomial size $poly(n)$ for $n\in Log$. We assume that this is defined in $PV_1$ in a natural way using an exact counting of sets of polynomial size $poly(n), n\in Log$ which is also definable in $PV_1$ in a usual way. This should not be confused with the definition of $Pr$ in $APC_1$.

\begin{prop}\label{ecc} For any $d,l<L$, 
$PV_1$ proves that for each $(n,d,\lambda)$-graph $G$ with $n\in Log$ vertices $V$, for any $S\subseteq V, |S|\leq |V|/2$, $$Pr_{(i,j)\in E(G^l)}[i\in S\wedge j\in S]\leq \frac{|S|}{|V|}(\frac{|S|}{|V|}+2\lambda^l)$$ where 
$E(G^l)$ denotes the set of all edges in $G^l$.
\end{prop}

\proof For empty $S$ the statement holds. Otherwise put $S:=\{i_1,...,i_{|S|}\}$. If $\left<x,\mathbf{1}\right>=0$, then $\left<Ax,\mathbf{1}\right>=0$ for the random-walk matrix $A$ of $G$. As $A^l$ is the random-walk matrix of $d^l$-regular graph $G^l$, $A^{l-1}\in Q^{n\times n}/d^{l-1}$ and 
$\frac{A^{l-1}x}{||A^{l-1}x||}\in Q^n/(Ln((d^{l-1}n)^nn)^{2n})^7(d^{l-1}n)^nn$ for $x\in Q^n/n$.
By the choice of $d,l$, this does not exceed the range $(Ln)^{{Ln}^L}$ and we can apply $\lambda(G)\leq \lambda$ to obtain $||A^lx||\leq \lambda^l||x||$ for any $x\in Q^n/n$ with $\left<x,\mathbf{1}\right>=0$. 
Now, use the inequality from the proof of Proposition \ref{expans}:
$$\frac{|E(S,V-S)|}{d^l}\geq \frac{|S||V-S|(1-\lambda^l)}{|V|}-\frac{\lambda^l}{Ln^2}$$ Then, $Pr_{(i,j)\in E(G^l)}[i\in S\wedge j\in S]=\frac{1}{|V|}\Sigma^{|S|}_{m=1}(1-Pr[j\notin S|i=i_m])$ is $$\frac{|S|}{|V|}(1-\Sigma_{m=1}^{|S|} \frac{|E(i_m,V-S)|}{|S|d^l})=\frac{|S|}{|V|}(1-\frac{|E(S,V-S)|}{|S|d^l})\leq \frac{|S|}{|V|}(\frac{|S|}{|V|}+2\lambda^l)\eqno{\qEd}$$

\subsection{A technical tool}\label{techtool}

Sometimes we will need to use an assumption which has the form ``$||Ax||\leq \lambda$ for $x\in Q^n/(Ln)^{(Ln)^L}$'' even for $x$'s exceeding the range fixed by $(Ln)^{(Ln)^L}$. We will now prove a simple approximation lemma which allows this in some cases. It illustrates a type of approximation which we use more often. The matrix $A$ in its formulation will not need to represent a random-walk matrix. In our applications $A$ will be a result of certain operations on random-walk matrices.

\begin{prop}\label{flow}(in $PV_1$) Let $A$ be an $n\times n$ matrix of elements from $Q/(2L^2n^5d)$, for $n\in Log$. Further, let $s\in Log$. If $||Ax||^2\leq y(||x||^2+1/L)$ for any $x\in Q^n/(Ln)^{(Ln)^L}$, then for any $x\in Q^n/m$, $$||Ax||^2\leq (y(1+\frac{1}{L})+\frac{1}{L})(||x||^2+\frac{1}{Ls})$$
\end{prop}

\proof 

For $x\in Q^n/m$ and $s\in Log$, define $||x||'$ in the same way as $||x||$ but with $SQRT$ redefined so that $0\leq (SQRT(||x||^2))^2-||x||^2\leq 1/(Ls)$. 


It suffices now to approximate $\frac{x}{||x||'}, x\neq 0$ by $c\in Q^n/(Ln)^{(Ln)^L}$ with $||c||^2\leq 1$ such that $|||A\frac{x}{||x||'}||^2-||Ac||^2|\leq \frac{1}{L}$. Then,
\medskip
 
   $||Ax||^2\leq ||A\frac{x}{||x||'}||^2(||x||^2+\frac{1}{Ls})\leq (y(||c||^2+1/L)+\frac{1}{L})(||x||^2+\frac{1}{Ls})\leq$

   $\leq (y(1+1/L)+\frac{1}{L})(||x||^2+\frac{1}{Ls})$
\medskip

The approximation: for each $i$, $|\frac{x_i}{||x||'}|\leq 1$ so we can find $c_i$ (i.e. $PV_1$ can prove its existence) such that $0\leq \frac{x_i}{||x||'}-c_i\leq 1/(18L^5n^{13}d^2)$. Then $||c||^2\leq ||\frac{x}{||x||'}||^2\leq 1$ and for each $l$, $|A_{l,i}\frac{x_i}{||x||'}-A_{l,i}c_i|\leq 1/(9L^3n^8d)$. Hence, $|(A\frac{x}{||x||'})_l-(Ac)_l|\leq 1/(9L^3n^7d)$. As $(A\frac{x}{||x||'})_l,(Ac)_l\leq 3L^2n^6d$, we can conclude $|||A\frac{x}{||x||'}||^2-||Ac||^2|\leq 1/L$ \qed

\smallskip

Using a similar approximation, we will derive one more useful lemma.

\smallskip

For any $n\times n$ matrix $A$ with elements from $Q/m$, we say that $||A||\leq 1$ iff for every $x\in Q^n/(Ln)^{(Ln)^L}$, $||Ax||^2\leq (1+2/L)(||x||^2+1/L)$.

\begin{prop}\label{myreorder} For any $\lambda$ and $d<L$, $PV_1$ proves the following. Let $A$ be a random-walk matrix of a $d$-regular graph $G$ with $n\in Log$ vertices such that $\lambda(G)\leq \lambda\in Q/(Ln^2)$. Let $J$ be $n\times n$ matrix such that $J_{i,j}=1/n$ for every $i,j$. Then, $$A=(1-\lambda)J+\lambda C$$ for some $C$ with $||C||\leq 1$
\end{prop}

\proof Define $C:=\frac{1}{\lambda}(A-(1-\lambda)J)$ $\in Q^{n\times n}/(2L^2n^5d)$. We want to prove that for any $x\in Q^n/(Ln)^{(Ln)^{L}}$, $||Cx||^2\leq (||x||^2+1/L)(1+2/L)$. Decompose $x$ as $x=\alpha \mathbf{1}+y$ for some $\alpha\in Q/((Ln)^{(Ln)^L})^{n+1}$ where $\left<\mathbf{1},y\right>=0$. 

Similarly as in Proposition \ref{flow}, approximate $\frac{y}{||y||}$  by vector $c$ with $||c||^2\leq 1$ so that $||A\frac{y}{||y||}||^2\leq ||Ac||^2+\lambda^2/L$ and $\frac{c}{||c||}\in Q^n/(Ln)^{(Ln)^L}$. This time we can do it without the absolute value because all elements of $A$ are positive. Note also that for $d<L$ the range of $\frac{c}{||c||}$ does not exceed $(Ln)^{(Ln)^L}$.

Since $A\mathbf{1}=\mathbf{1}$ and $J\mathbf{1}=\mathbf{1}$, we have $C\alpha \mathbf{1}=\alpha \mathbf{1}$. As $\left<y,\mathbf{1}\right>=0$, $Jy=0$ and $Cy=\frac{1}{\lambda}Ay$. Using $\left<Ay,\alpha \mathbf{1}\right>=0$ and $||Ac||\leq \lambda ||c||$, we obtain,
\medskip

 $||Cx||^2=||\alpha 1+\frac{1}{\lambda}Ay||^2=||\alpha \mathbf{1}||^2+||\frac{1}{\lambda}Ay||^2\leq ||\alpha\mathbf{1}||^2+\frac{1}{\lambda^2}(||Ac||^2+\frac{\lambda^2}{L})(||y||^2+\frac{1}{L})\leq$ 
 
 $||\alpha \mathbf{1}||^2+(1+2/L)(||y||^2+1/L)\leq (1+2/L)(||x||^2+1/L)$ \qed

\subsection{The tensor product}\label{sec:tens}

The explicit construction of the $(n,d,\lambda)$-graphs needs two graph products, the tensor product and the replacement product, which we describe in this and the next section. More details about the tensor product and the replacement product can be found in \cite[Section 21.3.3]{ba} resp. \cite[Section 21.3.4]{ba} .

\begin{defi}(in $PV_1$)
If $A=\{a_{i,j}\}_{i,j=1,...,n}$ is the $n\times n$ random-walk matrix of $d$-degree graph $G$ and $A'=\{a'_{i',j'}\}$ is the $n'\times n'$ random-walk matrix of $d'$-degree graph $G'$, then the random-walk matrix of $G\otimes G'$, denoted as $A\otimes A'$ is the $nn'\times nn'$ matrix that in the $\left<i,i'\right>$th row and the $\left<j,j'\right>$th column has the value $a_{i,j}a'_{i',j'}$.
\end{defi} 

This means that $G\otimes G'$ has a cluster of $n'$ vertices for every vertex in $G$. If $(i,j)$ is an edge in $G$ and $(i',j')$ is an edge in $G'$, then there is an edge between the $i'$-th vertex in the cluster corresponding to $i$ and the $j'$-th vertex in the cluster corresponding to $j$. Therefore, $G\otimes G'$ has degree $d'd$ and $nn'$ vertices. We can see matrix $A\otimes A'$ as consisting of blocks of the form $a_{i,j}A'$, that is, intuitively, $A\otimes A'$ is matrix $A$ with elements multiplied by copies of $A'$.

\smallskip

In Peano Arithmetic, $\lambda(G\otimes G')\leq max\{\lambda(G),\lambda(G')\}$ for regular graphs $G,G'$. The standard derivation of this bound uses the existence of an orthogonal basis of eigenvectors for symmetric matrices which uses the fundamental theorem of algebra (applied to determinant of matrix $A-xI$ consisting of exponentially many terms). We do not know how to formalize this in $PV_1$. Instead, we will derive a weaker bound which is sufficient for our purposes.

\smallskip

Note first a simple consequence of Cauchy-Schwarz inequality.

\begin{prop}(in $PV_1$) For every two $n\times n$ matrices $A,B$ and $x\in Q^n/m$ where $n\in Log$, we have $||(A+B)x||\leq ||Ax||+||Bx||+1/L^{1/2}$.
\end{prop}

\proof $||(A+B)x||^2=\left<(A+B)x,(A+B)x\right>=||Ax||^2+2\left<Ax,Bx\right>+||Bx||^2\leq$

\quad\quad\quad\quad\quad\quad\quad \ $\leq  ||Ax||^2+2||Ax||||Bx||+||Bx||^2\leq (||Ax||+||Bx||)^2$ 
\medskip

\noindent and so $||(A+B)x||\leq ||Ax||+||Bx||+1/L^{1/2}$. \qed

\begin{prop}\label{tens} $PV_1$ proves that if $G$ is a $d$-regular graph with $n\in Log$ vertices and $G'$ is a $d'$-regular graph with $n'\in Log$ vertices such that $d,d'<L$, $\lambda(G)\leq \lambda\in Q/(Ln^2)$ and $\lambda(G')\leq \lambda'\in Q/(L(n')^2)$, then $$\lambda (G\otimes G')\leq ((1+6/L)^2+1/L)(max\{ \lambda+\lambda'-\lambda\lambda',\lambda\lambda',\lambda',\lambda\})+3/L^{1/2}$$ (Note that $PV_1$ does not need to know that $\lambda(G)\leq 1$ or $\lambda(G')\leq 1$.) 
\end{prop}

\proof Let $A$ be the random-walk matrix of $G$ of the form $n\times n$ and $A'$ be the random-walk matrix of $G'$ of the form $n'\times n'$. By Proposition \ref{myreorder} $A=(1-\lambda)J_{n}+\lambda C$ for some $C$ with $||C||\leq 1$ and $n\times n$ all $1/n$ matrix $J_n$. Similarly, $A'=(1-\lambda')J_{n'}+\lambda' C'$ for some $C'$ with $||C'||\leq 1$ and $n'\times n'$ all $1/n'$ matrix $J_{n'}$. 

As tensor product satisfies $(A+B)\otimes C=A\otimes C+B\otimes C$ and $A\otimes (B+C)=A\otimes B+A\otimes C$, for any $x\in Q^{nn'}/(Lnn')^{(Lnn')^L}$ we have ($*$):
\medskip

$||A\otimes A' x||\leq (1-\lambda)||(J_{n}\otimes J_{n'})x||+(1-\lambda)\lambda'||(J_n\otimes C')x||$

$\ \ \ \ \ \ \ \ \ \ \ \ \ \ \ \ \ \ +\lambda(1-\lambda')||(C\otimes J_{n'})x||+\lambda'\lambda||(C\otimes C')x||+3/L^{1/2}$
\medskip

If $\Sigma_i x_i=0$, then $J_n\otimes J_{n'}x=0$. 
If $x\in Q^n/(Ln)^{(Ln)^L}$, $||J_nx||^2=\frac{1}{n}(\Sigma_i x_i)^2\leq ||x||^2$ where we used $\left<x,(1,...,1)\right>^2\leq n||x||^2$ which follows from Cauchy-Schwarz inequality. Therefore, $||J_n||\leq 1$ and similarly $||J_{n'}||\leq 1$. 

If $\lambda>1$ or $\lambda'>1$, we can trivially upper bound the term corresponding to $1-\lambda$ resp. $1-\lambda'$ in $(*)$ by 0. In all cases, to finish the proof it suffices to show that for any $n\times n$ matrix $A\in Q^{n\times n}/(2L^2n^5d)$, $n'\times n'$ matrix $B\in Q^{n'\times n'}/(2L^2(n')^5d)$ such that $||A||\leq 1,||B||\leq 1$, for any $x\in Q^{nn'}/(Lnn')^{(Lnn')^L}$ with $||x||=1$, $||(A\otimes B)x||\leq (1+6/L)^2+1/L$ holds.
\medskip

For any $x\in Q^{nn'}/m'$ and $i\in [n']$ define $x^i\in Q^n/m$ so that for each $j\in [n]$, $$x^i_j=\Sigma_{k\in \{n'(j-1)+1,...,n'j\}} B_{i,(k-n'(j-1))}x_k$$ Then, $||(A\otimes B)x||^2=\Sigma_{i\in [n']}||Ax^i||^2$ and as by Proposition \ref{flow} for each $i$,\\ $||Ax^i||^2\leq (\Sigma_{j\in [n]} (x^i_j)^2+1/(Ln'))((1+1/L)(1+2/L)+1/L)$, we have, $$||(A\otimes B)x||^2\leq (1/L+\Sigma_{i\in [n']}\Sigma_{j\in [n]}(x^i_j)^2)(1+6/L)$$ Since also $||Bx||^2\leq (||x||^2+1/(Ln))((1+1/L)(1+2/L)+1/L)$, for each $j\in [n]$, $$\Sigma_{i\in [n']}(\Sigma_{k\in \{n'(j-1)+1,...,n'j\}}B_{i,(k-n'(j-1))}x_k)^2\leq (\frac{1}{Ln}+\Sigma_{k\in \{n'(j-1)+1,...,n'j\}}(x_k)^2)(1+\frac{6}{L})$$ Therefore, if $||x||=1$, then $||(A\otimes B)||^2\leq (1/L+(1+6/L) (1+1/L))(1+6/L)$, and $||(A\otimes B)x||\leq (1+6/L)^2+1/L$.
 \qed

\subsection{The replacement product}

If $G$ is an $n$-vertex $d$-degree graph, we can give a number from 1 to $d$ to each neighbor of each vertex and then the rotation map $\hat{G}:[n]\times [d]\mapsto [n]\times [d]$ maps a pair $\left<v,i\right>$ to $\left<u,j\right>$ where $u$ is the $i$-th neighbor of $v$ and $v$ is the $j$-th neighbor of $u$. Using this rotation map, we define the replacement product.
\medskip

Let $G,G'$ be graphs such that $G$ has $n$ vertices and degree $D$, and $G'$ has $D$ vertices and degree $d$. Further, let $A,A'$ denote the random-walk matrices of $G$ and $G'$ respectively, and $\hat{A}$ be the permutation matrix corresponding to the rotation map of $G$ which means that $\hat{A}$ is an $nD\times nD$ matrix whose $(i,j)$th column is all zeroes except a single 1 in the $(i',j')$ position where $(i',j')=\hat{G}(i,j)$. Then the replacement product of $G$ and $G'$, denoted $G\oslash G'$, is the graph with the random-walk matrix $$A\oslash A':=1/2\hat{A}+1/2(I_n\otimes A')$$ where $I_n$ is $n\times n$ 0-1 matrix with 1's only on the diagonal.

This means that $G\oslash G'$ has a copy of $G'$ for every vertex in $G$ and if $(i,j)$ is an edge in $G$ then there are $d$ parallel edges between the $i'$-th vertex in the copy of $G'$ corresponding to $i$ and the $j'$ vertex in the copy of $G'$ corresponding to $j$ where $i'$ is the index of $j$ as neighbor of $i$ and $j'$ is the index of $i$ as neighbor of $j$ in $G$. Therefore, $G\oslash G'$ has degree $2d$ and $nD$ vertices.
\medskip

\begin{prop}(in $PV_1$)\label{repl} Let $d,D<L$. Suppose $G$ is a $D$-degree graph with $n\in Log$ vertices and $G'$ is a $d$-degree graph with $D$ vertices. If $\lambda(G)\leq 1-\epsilon\in Q/(Ln^2)$ and $\lambda(H)\leq 1-\delta\in Q/(LD^2)$ for $n\in Log$, rational $\epsilon$ and rational $\delta \in [0,1]$, then $$\lambda((G\oslash H)^3)\leq (1-\epsilon\delta^2/8)(1+8/L^{1/2})^{9}+\delta^2/(2L^{1/2})+2/L^{1/2}$$
\end{prop}

In Proposition \ref{repl}, Peano Arithmetic could prove $\lambda(G\oslash H)\leq 1-\frac{\epsilon\delta^2}{24}$ following the argument in Arora-Barak \cite{ba}. In \cite{ba} this is derived using the equation $\lambda(G^l)=\lambda(G)^l$ which uses the existence of an orthogonal basis of eigenvectors for symmetric matrices. Again, in $PV_1$ we prove just a weaker bound for $(G\oslash H)^3$ (i.e. not for the product $G\oslash H$ but its power) which is sufficient for our purposes.

\proof Let $A$ resp. $B$ be the random-walk matrix of graph $G$ with $n$ vertices resp. graph $H$ with $D$ vertices and $\hat{A}$ be the permutation matrix corresponding to the rotation map of $G$. By definition, $A\oslash B=\frac{1}{2}(\hat{A}+I_n\otimes B)$ and 
\[\eqalign{
  (A\oslash B)^3={}&\frac{1}{8}(\hat{A}^3+\hat{A}(I\otimes B)\hat{A}+(I\otimes B)\hat{A}^2+(I\otimes B)^2\hat{A}+\hat{A}^2(I\otimes B)+\cr
&+\hat{A}(I\otimes B)^2+(I\otimes B)\hat{A}(I\otimes B)+(I\otimes B)^3)
}\]
By Proposition \ref{myreorder}, $B=\delta J+(1-\delta)C$ for some $C$ with $||C||\leq 1$ and $D\times D$ all $1/D$ matrix $J$. Therefore, 
\[\eqalign{
  (I\otimes B)\hat{A}(I\otimes B)
={}&\delta^2(I\otimes J)\hat{A}(I\otimes J)+\delta(1-\delta)(I\otimes J)\hat{A}(I\otimes C)+\cr
&+\delta(1-\delta)(I\otimes C)\hat{A}(I\otimes J)+(1-\delta)^2(I\otimes C)\hat{A}(I\otimes C)
}\]
Since $||C||\leq 1$ and $||I||\leq 1$, for any $x$ with $||x||\leq 1$, we have $||(I\otimes C)x||^2\leq (1+6/L)^4$ as in the proof of Proposition \ref{tens}. Similarly, $||(I\otimes J)x||^2\leq (1+6/L)^4$.

If a matrix $A$ satisfies $||Ax||^2\leq (1+6/L)^4$ for $||x||\leq 1$, then for any $B$ and $x$, 
$||(AB)x||^2=||A\frac{Bx}{||Bx||}||^2(SQRT(||Bx||^2))^2\leq (1+\frac{6}{L})^4(SQRT(||Bx||^2))^2$. Consequently, $||(AB)x||\leq (1+6/L)^2||Bx||+1/L^{1/2}$.

As $||\hat{A}||\leq 1$, this shows that for any $x,||x||\leq 1$ and $\delta \in [0,1]$,  
\[\eqalign{
||((I\otimes B)\hat{A}(I\otimes B))x||
\leq{}& \delta^2||((I\otimes J)\hat{A}(I\otimes J))x||+(1-\delta^2)((1+\frac{6}{L})^8+\cr
&+(1+\frac{6}{L})^4/L^{1/2}+(1+\frac{6}{L})^2/L^{1/2}+\frac{1}{L^{1/2}})+\frac{3}{L^{1/2}}
}\]
Further, for any $x, ||x||=1$ and $\delta \in [0,1]$, 
\[||(I\otimes B)x||\leq \delta||(I\otimes J)x||+(1-\delta)||(I\otimes C)x||+1/L^{1/2}\leq (1+6/L)^2+2/L^{1/2}\]
Hence, $||(I\otimes B)x||^2\leq (1+8/L^{1/2})^4$, and using an analogous argument as above we can bound $||(A\oslash B)^3x||$. For any $x,||x||=1$,
\[||(A\oslash B)^3x||\leq (1-\frac{\delta^2}{8})(1+8/L^{1/2})^{9}+\frac{\delta^2}{8}||((I\otimes J)\hat{A}(I\otimes J))x||+2/L^{1/2}\]
Observe that $(I\otimes J)\hat{A}(I\otimes J)=A\otimes J$ because $(I\otimes J)\hat{A}(I\otimes J)$ is the random-walk matrix of a graph with the number of edges between its nodes $(i,j)$ and $(i',j')$ being the number of $k$'s in $[D]$ for which there is $k'$ such that 
$\hat{G}(i,k)=(i',k')$. That is,
\[((I\otimes J)\hat{A}(I\otimes J))_{(i,j),(i',j')}=\frac{1}{D}a_{i,i'}=(A\otimes J)_{(i,j),(i,j')}\]
Then, by Proposition \ref{tens}, for any $x, ||x||=1$ such that $\Sigma_i, x_i$ (and so $Jx=0$) we have: 
\[||(I\otimes J)\hat{A}(I\otimes J)x||=||(A\otimes J)x||\leq (1-\epsilon)((1+6/L)^2+1/L)+3/L^{1/2}\]
which completes the proof.  \qed

\subsection{The construction of the \texorpdfstring{$(n,d,\lambda)$}{(n,d,lambda)}-graphs}

Finally, we are ready to construct the $(n,d,\lambda)$-graphs in the theory $PV_1$, see Arora-Barak \cite[Chapter 21]{ba} for the history of the result. However, we will do it just for $n$'s of the form $c^k$ where $c$ is a constant and $k\in LogLog$. It is possible to extend the construction to any $n$ (cf. \cite{ba})  but at least a straightforward application of the extension requires algebraic techniques which we are avoiding. More specifically, it uses a converse of Proposition \ref{expans} which in turn uses facts about eigenvectors derived from the fundamental theorem of algebra. Nevertheless, the weaker construction is sufficient to derive the PCP theorem in $PV_1$.

\begin{prop}\label{cst} For any rational $c\in (0,1)$ there are $d,b$ and $L$ (the constant from the definition of $\lambda(G)$) such that $PV_1$ proves that for each $k\in LogLog$ and $n=(2d)^{100k}$ there is a $(2d)^b$-regular graph $G_n$ with $n$ vertices and $\lambda(G_n)<c$.
\end{prop}

\proof For $c\in (0,1)$, let $e$ be such that $1/2^e<c$ and $b>e$ be a sufficiently big constant. Then, define $((2d)^{100k}, (2d)^b,1/2^e)$-graphs in $PV_1$ as follows.
\medskip

\noindent 1. Let $H$ be a $((2d)^{100},d,0.01)$-graph where $d$ is a sufficiently big constant so that such a graph exists. Let $G_1$ be a $((2d)^{100},(2d)^b,\frac{1}{2^b})$-graph and $G_2$ be a $((2d)^{200},(2d)^b,\frac{1}{2^b})$-graph. These graphs can be found by brute force, cf. \cite{ba}. More precisely, as our $H$ take the graph $H$ from the proof of Theorem 21.19 in \cite{ba} and as our $G_1,G_2$ take $G_1^b,G_2^b$ for $G_1,G_2$ from the same proof in \cite{ba}.
\smallskip

\noindent 2. For $(2d)^{100k}$ with $k> 2$, define $G_k:=((G_{\lfloor (k-1)/2\rfloor}\otimes G_{\lceil (k-1)/2\rceil})\oslash H)^b$ 
\medskip

Note that for given $(2d)^{100k}$, $G_k$ is produced by a specific p-time computation which exists provably in $PV_1$.
\medskip

\begin{clm}
For every $(2d)^{100k}$, $G_k$ is a $((2d)^{100k},(2d)^b,1/2^e)$-graph. 
\end{clm}

\proof The claim is proved by $\Pi^b_1(PV)$-LPIND induction. As graphs $G_k$ are constructed by a p-time function, the statement we want to obtain is $\forall\Sigma^b_1$. Hence, by $\forall\Sigma^b_1$-conservativity of $S^1_2$ over $PV_1$, we can work in the theory $S^1_2$ (which proves $\Pi^b_1(PV)$-LPIND).

 For $k=1,2$, $PV_1$ can verify the claim directly. For $(2d)^{100k}$ with $k>2$, let $n_k$ be the number of vertices of $G_k$. If $n_{\lfloor (k-1)/2\rfloor}=(2d)^{100\lfloor (k-1)/2\rfloor}$ and $n_{\lceil (k-1)/2\rceil}=(2d)^{100\lceil (k-1)/2\rceil}$, then $n_k=n_{\lfloor (k-1)/2\rfloor}n_{\lceil (k-1)/2\rceil}(2d)^{100}=(2d)^{100k}$.

Considering the degree, if $G=G_{\lfloor (k-1)/2\rfloor}$ has degree $(2d)^b$, then $(G\otimes G)$ has degree $(2d)^{2b}$, $(G\otimes G)\oslash H$ has degree $2d$ and $G_k$ has degree $(2d)^b$.

The eigenvalue analysis: if $\lambda(G)\leq 1/2^e$ (which is a $\Pi^b_1(PV)$-formula), then assuming $L$ is sufficiently big, $1/2^e\in Q/(Ln^2)$ and  by Proposition \ref{tens} $\lambda(G\otimes G)\leq 2/2^e$. Hence, 
by Proposition \ref{repl}, 
\[\lambda(((G\otimes G)\oslash H)^3)\leq (1-(1-2/2^e)\frac{(0.99)^2}{8})(1+8/L^{1/2})^{9}+\frac{(0.99)^2}{2L^{1/2}}+2/L^{1/2}\]
The conclusion $\lambda(((G\otimes G)\oslash H)^b)\leq 1/2^e$ is a consequence of the fact that the assumption $\lambda(G)\leq \lambda$ implies $\lambda(G^b)\leq \lambda^b(1+4/L)+3d^b/L^{1/2}$ (where $L$ is quantified after $d,b$ so the term $3d^b/L^{1/2}$ can be made arbitrarily small). To see that the implication holds, note that similarly as in the proof of Proposition \ref{ecc}, $\lambda(G)\leq \lambda$ implies that for any $x\in Q^n/((Ln^3)^nn)$ with  $\left<x,\mathbf{1}\right>=0$, we have $||A^bx||\leq \lambda^b||x||$ where $A^b\in Q^{n\times n}/d^b$ is the random-walk matrix of $G^b$. We need a similar bound even for $x\notin Q^n/((Ln^3)^nn)$. Fortunately, if $x\notin Q^n/((Ln^3)^nn), ||x||=1$, $\left<x,1\right>=0$, we can again approximate $x$ by vector $c\in Q^n/((Ln^3)^nn)$: for each $i$, $|x_i|\leq 1$ (otherwise $||x||>1$) so we can find $c_i\in Q/((Ln^3)^nn)$ such that $|x_i-c_i|\leq 1/(Ln^2)$ and $\left<c,1\right>=0$. The values $c_i$ are produced provably in $PV_1$ by a p-time algorithm which choses $i_0$ satisfying $x_{i_0}\geq 1/(Ln^2)$, then finds the smallest $c_i>x_i$ such that $c_i-x_i<1/(Ln^3), c_i\in Q/(Ln^3), i\neq i_0$ and puts $c_{i_0}=\Sigma_{i\neq i_0} c_i\in Q/((Ln^3)^nn)$. The chosen $c$ satisfies $||c||^2\leq 1+3/(Ln)$ and 
$|(A^bx)_j-(A^bc)_j|\leq d^b/(Ln)$. Since also $(A^bx)_j,(A^bc)_j\leq 2d^bn$, we have $|||A^bx||^2-||A^bc||^2|\leq 5d^{2b}/L$ and $$||A^bx||^2\leq\lambda^{2b}(||c||^2+1/L)+5d^{2b}/L\leq \lambda^{2b}(1+4/L)+5d^{2b}/L$$ Thus, $||A^bx||\leq \lambda^b(1+4/L)+3d^b/L^{1/2}$.
 \qed

Note that in the previous proposition, $d$ does not depend on $L$ and $b$ can be chosen arbitrarily big.

\section{The PCP theorem in \texorpdfstring{$PV_1$}{PV1}}\label{npcp}

The PCP theorem obtained in Arora-Safra \cite{san} and Arora et.al. \cite{sansud} (see Arora-Barak \cite[Chapter 22]{ba} for the history of the theorem) is a strengthening of the exponential PCP theorem in which the verifier $D$ uses only $O(\log n)$ random bits. Using these random bits, $D$ asks for at most $O(1)$ bits of  the given proof $\pi$. Hence, $\pi$ can be seen as a string of size $poly(n)$. In particular, it can be represented by a binary string in our formalization.

We will follow Dinur's \cite{din} simplified proof of the PCP theorem as it is presented in Arora-Barak \cite{ba}. This will go rather smoothly (once we have a suitable formalization of the $(n,d,\lambda)$-graphs) because the proof is combinatorial and it needs to count only sets of polynomial size. These are subsets of $\{1,...,poly(n)\}$ where $n\in Log$ for which we assume to have exact counting in $PV_1$ defined in a natural way.






\smallskip

Recall the verifier $D^{\pi,w}(x)$ from Definition \ref{dall}. In the standard definition, $\pi$ would be allowed to be a string of arbitrary length and $D$ would have an oracular access to $\pi$, it could ask for any bit of $\pi$.  Then, for a language $L$, $L\in PCP(\log n,1)$ standardly means that there is a p-time algorithm $D$ such that: 

\begin{enumerate}[label=\arabic*.]
\item If $x\in L$, then there is a string $\pi$ such that $D$ with input $x$ of length $n$ and $O(\log n)$ random bits asks for at most $O(1)$ bits of $\pi$ and accepts (with probability 1); 
\item If $x\notin L$, then for any $\pi$, $D$ with input $x$ of length $n$ and $O(\log n)$ random bits asks for at most $O(1)$ bits of $\pi$ and accepts with probability $\leq 1/2$.
\end{enumerate}

\noindent The PCP theorem says that $NP=PCP(\log n,1)$. In our formalization, proofs $\pi$ will be represented by p-size strings, hence, the statement of the PCP theorem is modified accordingly. As in the case of the exponential PCP theorem, we could alternatively represent proofs $\pi$ by oracles which would maybe better reflect the nature of the PCP theorem but then we would need to formalize the PCP theorem in a theory extended by such oracles.
\smallskip

In this Section we use the notion of probability $Pr$ on spaces of polynomial size $poly(n)$ which is assumed to be defined in a natural way using the exact counting of sets of polynomial size in $PV_1$. This should not be confused with the definition of $Pr$ in $APC_1$.
\smallskip

First we formalize the easier implication of the PCP theorem: $PCP(\log n,1)\subseteq NP$.

\begin{thm} Let $c,d,k$ be arbitrary constants, then $PV_1$ proves that for any $kn^k$-time algorithm $D$ there exists $2kcn^{2kc}$-time algorithm $M$ such that for each $x\in\{0,1\}^n$: $$\exists \pi\in\{0,1\}^{dn^c} \ \forall w<n^c, D^{\pi,w}(x)=1\rightarrow \exists y\in\{0,1\}^{dn^c}\ M(x,y)=1$$ 
$$\forall \pi\in\{0,1\}^{dn^c}\ Pr_{w<n^c}[D^{\pi,w}(x)=1]\leq 1/2\rightarrow \forall y\in\{0,1\}^{dn^c}\ M(x,y)=0$$
\end{thm}

\proof Given a $kn^k$-time algorithm $D$, define the algorithm $M$ as follows. $M$ accepts $x,y$ if and only if $y=(y_0,...,y_{n^c-1})\in\{0,1\}^{dn^c}$ with $y_i$'s in $\{0,1\}^d$ and for all the $y_i$'s the algorithm $D$ on input $x$, random bits $i$ and with access to $\pi$ which results in $d$ bits $y_i$ accepts.
\smallskip

Suppose there is $\pi\in\{0,1\}^{dn^c}$ such that for each $w<n^c$, $D$ on input $x$ with bits $r_w\in\{0,1\}^d$ obtained from $d$-times accessing $\pi$ accepts. Then for $y=(y_0,...,y_{n^c-1})$ with $y_{w}=r_{w}$ we have that for each $y_i\in y$ the algorithm $D$ on input $x$ and with access to $\pi$ which results in $d$ bits $y_i$ accepts. Therefore, $M(x,y)=1$.
\smallskip

Now assume that for any $\pi\in \{0,1\}^{dn^c}$, $Pr_{w<n^c}[D^{\pi,w}(x)=1]\preceq_0 1/2$. 
Then for any $y=(y_0,...,y_{n^c-1})$ with $y_i$'s in $\{0,1\}^d$ there is $y_i$ such that $D$ on $x$, random bits $i$, and with access to $\pi$ resulting in $y_i$ rejects. Otherwise, for some $\pi\in \{0,1\}^{dn^c}$ we have $\{w<n^c|D^{\pi,w}(x)=1\}=n^c$ contradicting the assumption. Hence, $M(x,y)=0$ . \qed

\smallskip

As the NP-completeness of SAT is provable in $PV_1$, the important implication of the PCP theorem, $PCP(\log n,1)\subseteq NP$, can be stated in $PV_1$ as Theorem \ref{pcp}.

\begin{theoremss:pcptt}[The PCP theorem in $PV_1$]\label{pcp} There are constants $d,k,c$ and a $kn^k$-time algorithm $D$ (given as a PV-function) computing as in Definition \ref{dall} such that $PV_1$ proves that for any $n\in Log$ and $x\in\{0,1\}^n$, $n\in Log$:

$$\exists y SAT(x,y)\rightarrow \exists \pi\in\{0,1\}^{dn^c}\ \forall w<n^c\ D^{\pi,w}(x)=1$$ 
$$\forall y \neg SAT(x,y)\rightarrow \forall \pi\in\{0,1\}^{dn^c}\  Pr_{w<n^c}[D^{\pi,w}(x)=1]\leq 1/2$$
\end{theoremss:pcptt}

The proof is summarized at the end of this section. It is a sequence of certain reductions between the so called CSP instances (CSP stands for constraint satisfaction problem) so we need to start with a reformulation of Theorem \ref{pcp} in terms of these reductions.

\begin{defi}[in $PV_1$] Let $q,W$ be constants, and $n,m\in Log$. A $qCSP_W$ instance $\phi$ is a collection of circuits $\phi_1,...,\phi_m$ (called constraints) mapping $[W]^n$ to $\{0,1\}$. Each $\phi_i$ is encoded by a binary string, it has $n$ inputs which are taking values that are bit strings in $\{0,1\}^{\log W}$ but depends on at most $q$ of them: for every $i\in [m]$ there exist $f_1,...,f_q\in [n]$ and $f:\{0,1\}^q\mapsto \{0,1\}$ such that $\phi_i(u)=f(u_{f_1},...,u_{f_q})$ for every $u\in [W]^n$. We say that $q$ is the arity of $\phi$. By $qCSP$ instance we mean a $qCSP$ instance with binary alphabet.

An assignment $u\in [W]^n$ satisfies $\phi_i$ if $\phi_i(u)=1$, and instance $\phi$ is satisfiable if $val(\phi):=max_{u\in [W]^n}\frac{\Sigma^m_{i=1}\phi_i(u)}{m}=1$.
\end{defi}

We will not need to prove the totality of the function $val(\phi)$ in $PV_1$. It will be sufficient for us to work with formulas of the form $val(\phi)\leq y$ which are $\Pi^b_1$.

\begin{defi}[in $PV_1$] Let $q,q',W,W'$ be arbitrary constants. A p-time function $f$ (given as a PV-function) mapping $qCSP_W$ instances to $q'CSP_{W'}$ instances, abbreviated as $f:qCSP_W\rightarrow q'CSP_{W'}$, is a $CL$-reduction (short for complete linear-blowup reduction) if 
for every $qCSP_W$ instance $\phi$:

\begin{itemize}
\item Completeness: If $\phi$ is satisfiable then so is $f(\phi)$.

\item Linear blowup: If there are $m$ constraints in $\phi$, then $f(\phi)$ has at most $Cm$ constraints and alphabet $W'$, where $C$ can depend on $q$ (but not on $m$ or the number of variables in $\phi$).
\end{itemize}

\noindent For a constant $k$, a function $f$ is $CL^k$-reduction if it is a $CL$-reduction computable in time $kn^k$.
\end{defi}

Theorem \ref{pcp} then follows from the following proposition.

\begin{prop}\label{cc} There are constants $q_0\geq 3, \epsilon_0>0$ and a $CL$-reduction $f:q_0CSP\rightarrow q_0CSP$ such that $PV_1$ proves that for every $q_0CSP$ instance $\phi$, every $\epsilon<\epsilon_0$, $$val(\phi)\leq 1-\epsilon \rightarrow val(f(\phi))\leq 1-2\epsilon$$
\end{prop}

\proof (of Theorem \ref {pcp} from Proposition \ref{cc}) The statement
we want to derive is a $\forall\Sigma^b_1$-formula. Hence, we can work
in the theory $S^1_2$. As $q_0\geq 3$, $q_0CSP$ is a generalization of
3SAT and by the NP-completeness of 3SAT (derived similarly as the
NP-completeness of SAT), for some $k'$, there is a $k'n^{k'}$-time
function $h$ mapping propositional formulas to $q_0CSP$ instances such
that for every $n\in Log$ and $x\in \{0,1\}^n$, $\exists y
SAT(x,y)\rightarrow val(h(x))=1$ and $\forall y\neg
SAT(x,y)\rightarrow val(h(x))\leq 1-1/m$ where $m\in Log$ is the
number of constraints in $h(x)$. Applying Proposition \ref{cc} we
obtain a $kn^k$-time function $f^{\log m}\circ h$ for some constant
$k$ such that 
\[\eqalign{
  \exists y SAT(x,y)\rightarrow val(f^{\log m}\circ h(x))&=1\cr
  \forall y\neg SAT(x,y)\rightarrow val(f^{\log m}\circ h(x))&\leq
  1-\epsilon_0
}\]
 Here, we used $\Pi^b_1$-LLIND (available in $S^1_2$) for
 $\Pi^b_1$-formulas $val(f^i(\phi))\leq 1-2^i\epsilon$ where $i\leq
 |m|$.  Therefore, for some constants $d',c'$, and an algorithm $D'$
 which given any formula $x$ and proof $\pi$ accepts if and only if
 $\pi$ encodes a satisfying assignment to randomly chosen constraint
 in $f^{\log m}\circ h(x)$ we have: 
\[\eqalign{
  \exists y SAT(x,y)\rightarrow \exists \pi\in\{0,1\}^{d'n^{c'}}\
  \forall w D'^{\pi,w}(x)^=1\cr
  \forall y\neg SAT(x,y)\rightarrow \forall \pi\in\{0,1\}^{d'n^{c'}}\
  Pr_{w}[D'^{\pi,w}(x)=1]&\leq 1-\epsilon_0
}\]
The gap can be amplified to 1/2 by choosing sufficiently many (but
constant number of) constraints in $f^{\log m}\circ h(x)$ and
accepting if and only if $\pi$ encodes satisfying assignments to all
of them. This requires Chernoff's bound but only over sets of
polynomial size for which we have exact counting in $PV_1$. \qed

\smallskip

Proposition \ref{cc} is an immediate consequence of the following two statements. The first one provides us a $CL$-reduction producing CSP instances which increase the gap between 0 and the minimal number of unsatisfied constraints. However, the alphabet of the resulting instances increases too. The second statement takes it back to binary while losing just a factor of 3 in the gap.

\begin{prop}[Gap amplification in $PV_1$]\label{cc1} 
For every $l,q$ there are $W,\epsilon_0$ and a $CL$-reduction $g_{l,q}: qCSP\rightarrow 2CSP_W$ such that $PV_1$ proves that for every $qCSP$ instance $\phi$ and for every $\epsilon<\epsilon_0$ $$val(\phi)\leq 1-\epsilon \rightarrow val(g_{l,q}(\phi))\leq 1-l\epsilon$$
\end{prop}

\begin{prop}[Alphabet reduction in $PV_1$]\label{cc2} There is $d$ such that for any $W$ there is a $CL$-reduction $h:2CSP_W\rightarrow dCSP$ such that $PV_1$ proves that for every $2CSP_W$ instance $\phi$, and for each $\epsilon$ $$val(\phi)\leq 1-\epsilon\rightarrow val(h(\phi))\leq 1-\epsilon/3$$
\end{prop}

Proposition \ref{cc} can be obtained from previous two propositions by taking $l=6$ in Proposition \ref{cc1} and $q=max\{d,3\}$ for $d$ from Proposition \ref{cc2}.
\smallskip

We firstly derive Proposition \ref{cc2} using the following application of the exponential PCP theorem which is scaled down so that we need to reason only about sets of constant size.

\begin{prop}\label{prox} There are constants $d,k'$ and an algorithm $D$ such that for every $s$, $PV_1$ proves: given any $s$-size circuit $C$ with $2n_1$ inputs, $D$ runs in time $s^{k'}$, examines $\leq d$ bits in the provided strings and
\begin{enumerate}[label=\arabic*.]

\item If $C(u_1, u_2)=1$ for $u_1,u_2 \in\{0,1\}^{n_1}$, there is a string $\pi_3$ of size $2^{s^{k'}}$ such that 
\[\forall w<2^{s^{k'}}\ D^{(WH(u_1),WH(u_2),\pi_3),w}(C)=1.\]

\item For bit strings $\pi_1,\pi_2,\pi_3$ where
  $\pi_1,\pi_2\in\{0,1\}^{2^{n_1}}$, $\pi_3\in\{0,1\}^{2^{s^{k'}}}$,
  if 
\[Pr_{w<2^{s^{k'}}}[D^{(\pi_1,\pi_2,\pi_3),w}(C)=1]\geq 1/2\]
 then 
\[Pr_{w< 2^{n_1}}[(\pi_1)_w=WH(u_1)(w)]\geq 0.99\ \mbox{and}\ 
  Pr_{w<2^{n_1}}[(\pi_2)_w=WH(u_2)(w)]\geq 0.99
\]
for some $u_1,u_2\in\{0,1\}^{n_1}$ such that $C(u_1, u_2)=1$.
\end{enumerate}
\end{prop}

\proof (of Proposition \ref{cc2} from Proposition \ref{prox}) The $CL$-reduction $h$ works as follows. Let $\phi$ be a $2CSP_W$ instance with constraints $\phi_1,\phi_2,...,\phi_m$ on variables $u_1,...,u_n$ which are taking values that are in $\{0,1\}^{\log W}$. Each constraint $\phi_S(u_i,u_j)$ is a circuit applied to the bit strings representing $u_i,u_j$. Without loss of generality $s\leq 2^{4\log W}$ is an upper bound on the size of this circuit. 
\smallskip

Given such $\phi$, $h$ replaces each variable $u_i$ by a sequence $U_i=(U_{i,1},...,U_{i,W})$ of $W$ binary variables ($U_i$ is long enough to represent $WH(u_i)$). Then, for each constraint $\phi_S(u_i,u_j)$ it applies Proposition \ref{prox} where $\phi_S(u_i,u_j)$ is the circuit whose assignment is being verified. The resulting $s^{k'}$-time algorithm $D$ can be represented as a $2^{s^{O(1)}}$-size $dCSP$ instance $\psi_S(U_i,U_j, \Pi_S)$ where $U_i,U_j$ play the role of $\pi_1,\pi_2$ and $2^{s^{k'}}$ new binary variables $\Pi_S$ play the role of $\pi_3$. The arity $d$ of $\psi_S(U_i,U_j,\Pi_S)$ is the number of bits $D$ reads in the proof which is a fixed constant independent of $W$ and $\epsilon$. The instance $\psi_S(U_i,U_j,\Pi_S)$ contains one constraint for each possible random string in $D$, so the fraction of its satisfied constraints is the acceptance probability of $D$. The $CL$-reduction $h$ thus maps $2CSP_W$ instances $\phi$ to $dCSP$ instances $\psi$ where each $\phi_S(u_i,u_j)$ is replaced by a $dCSP$ instance $\psi_S(U_i,U_j,\Pi_S)$. As $2^{s^{O(1)}}$ is a constant independent of $m$ and $n$, linear blowup is preserved.
\smallskip

If $\phi$ is satisfiable, then by property 1 in Proposition \ref{prox} so is $\psi$. We want to show that if some assignment satisfies more than $1-\epsilon/3$ fraction of the constraints in $\psi$, then we can construct an assignment for $\phi$ satisfying more then $1-\epsilon$ fraction of its constraints: For each $i$, if $U_i$ is 0.99-close to some linear function $WH(a_i)$, i.e. $Pr_x[U_{i,x}=WH(a_i)(x)]\geq 0.99$,  then use (the determined) $a_i$ as the assignment for $u_i$, and otherwise use arbitrary string. The algorithm is p-time because the size of each $U_i$ is constant. If the decodings $a_i,a_j$ of $U_i,U_j$ do not satisfy $\phi_S(u_i,u_j)$, then by property 2 in Proposition \ref{prox} at least half of constraints in $\psi_S$ is not satisfied. Hence, the fraction of unsatisfied constraints in $\phi$ is $<2\epsilon/3$. 
\qed 

\proof (of Proposition \ref{prox}) $PV_1$ can prove the statement from Proposition \ref{prox} simply by examining all possible cases of which there is a constant number. Hence, the provability of the statement follows from it being true. Nevertheless, we present also the standard proof itself.
\smallskip

The algorithm $D$ firstly reduces the problem of satisfiability of the given circuit $C$ with $s$ wires (inputs are considered as wires in the circuit) to the question of solvability of a set of quadratic equations with $t=s^{O(1)}$ variables similarly as in the proof of the exponential PCP theorem. $D$ expects $\pi_3$ to contain linear functions $f,g$ which are $WH(z)$ and $WH(z\otimes z)$ respectively for $z\in\{0,1\}^t$ satisfying the set of quadratic equations and checks these functions as in the exponential PCP theorem. Moreover, $D$ checks that $\pi_1$ and $\pi_2$ are 0.99-close to some linear functions. That is, if $D$ accepts $\pi_1,\pi_2,\pi_3$ with probability $\geq 1/2$, it is because the set of quadratic equations is satisfiable and $Pr_w[(\pi_1)_w=WH(u_1)(w)]\geq 0.99$, $Pr_w[(\pi_2)_w=WH(u_2)(w)]\geq 0.99$ for some $u_1,u_2\in\{0,1\}^{n_1}$.

Finally, $D$ checks that $\pi_1,\pi_2$ encode strings whose concatenation is the same as the first $2n_1$ bits of the string encoded by $f$ (without loss of generality the first $2n_1$ bits encode satisfying assignement for $C$) by performing the following concatenation test:
\smallskip

Pick random $x,y\in\{0,1\}^{n_1}$ and denote by $XY\in\{0,1\}^t$ the string whose first $n_1$ bits are $x$, the next $n_1$ bits are $y$ and the remaining bits are all 0. Accept if and only if $f(XY)=\pi_1(x)+\pi_2(y)$.
\smallskip

The algorithm $D$ runs in time $s^{k'}$ and examines $\leq d$ bits in $\pi_1,\pi_2,\pi_3$ for some constants $k',d$. It satisfies the first property from Proposition \ref{prox}. Moreover, assuming that $\pi_1=WH(u), \pi_2=WH(v)$ and $z$ is the string encoded by a linear function $f$, the concatenation test rejects with probability 1/2 if $u, v$ differs from the first $2n_1$ bits of $z$. Hence, if $D$ accepts $\pi_1,\pi_2,\pi_2$ with probability $\geq 1/2$, it is because $\pi_1,\pi_2$ are 0.99-close to linear functions encoding $u_1,u_2$ such that $C(u_1, u_2)=1$.
 \qed

In the rest of this section we derive Proposition \ref{cc1}. To do this, we will need two facts about probability:

\begin{prop}\label{stat}

\noindent 1. Let $t$ be a square and $S_t$ be the binomial distribution over $t$ fair coins, i.e. $Pr[S_t=k]=t!/((t-k)!k!)2^{-t}$. Then for $i\in\{0,1\}$ and any $\delta$ such that $0\leq\delta<1$, $PV_1$ proves: $$\Sigma_k |Pr[S_t=k]-Pr[S_{t+(-1)^i\lfloor\delta\sqrt{t}\rfloor}=k]|\leq 20\delta$$

\noindent 2. For any $k$, $PV_1$ proves that for each $n\in Log$, if $V$ is a nonnegative random variable defined on a sample space of size $n^k$, then $Pr[V>0]\geq E[V]^2/E[V^2]$.
\end{prop}

The first part of Proposition \ref{stat} is an estimation of a so called statistical distance of two binomial distributions which is known to hold (see \cite{ba} page 469) and as all its parameters are quantified outside of the theory $PV_1$, it is trivially provable by an explicit ``brute force'' enumeration.

The second part is obtained from a simple expansion: $$(E[X])^2=(E[X\cdot 1_{X>0}])^2\leq E[X^2] E[(1_{X>0})^2]=E[X^2]Pr[X>0]$$ where we used a form of Cauchy-Schwarz inequality $E[XY]^2\leq E[X^2]E[Y^2]$ which can be derived in the same way as our Cauchy-Schwarz inequality from Section \ref{rest} but with $\left< x,y\right> :=E[XY]$.

The proof of Proposition \ref{cc1} is divided into two parts. The first part shows how to reduce any $qCSP$ instance into a $2CSP_W$ instance which is nice (in a sense defined below) and the second part gives us a CL-reduction from nice instances which amplifies the gap as it is required in Proposition \ref{cc1}.

\begin{defi} (in $PV_1$) \hfill
\begin{enumerate}[label=\arabic*.]
\item Let $\phi$ be a $2CSP_W$ instance mapping $[W]^n$ to $\{0,1\}$. The constraint graph of $\phi$ is the graph $G$ with vertex set $[n]$ where for every constraint $\phi$ depending on the variables $u_i, u_j$, the graph $G$ has the edge $(i,j)$. $G$ is allowed to have parallel edges and self-loops. Then $G$ is $d$-regular for some constant $d$ independent of $W$, and at every node, at least half the edges incident to it are self-loops.

\item A $qCSP_W$ instance $\phi$ is nice if $q=2$ and the constraint
  graph of $\phi$ denoted $G$ satisfies $\lambda(G)\leq 0.9$.
\end{enumerate}
\end{defi}

\noindent The reduction into nice instances which we need is a consequence of the following three Propositions.

\begin{prop}\label{n1} There is a constant $k$ such that for every $q$ there is a $CL^k$-reduction $h:qCSP\rightarrow 2CSP_{2^q}$ such that $PV_1$ proves that for any $qCSP$ instance $\phi$ and any $\epsilon$ $$val(\phi)\leq 1-\epsilon\rightarrow val(h(\phi))\leq 1-\epsilon/q$$
\end{prop}

\proof  The $CL^k$ reduction works as follows. Given $qCSP$ instance $\phi$ over $n$ variables $u_1,...,u_n$ with $m$ constraints, it produces $2CSP_{2^q}$ instance $\psi$ over the variables $u_1,...,u_n$, $y_1,...,y_m$ such that for each $\phi_i$ in $\phi$ depending on the variables $u_{f_1},...,u_{f_q}$, $\psi$ contains $q$ constraints $\psi_{i,j}, j=1,...,q$ where $\psi_{i,j}(y_i,u_{f_j})$ is true iff $y_i$ encodes an assignment to $u_{f_1},...,u_{f_q}$ satisfying $\phi_i$ and $u_{f_j}\in \{0,1\}$ agrees with the assignment $y_i$. 

The number of constraints in $\psi$ is $qm$ and if $\psi$ is satisfiable, then so is $\psi$. Suppose that $val(\phi)\leq 1-\epsilon$ and let $u_1,...,u_n,y_1,...,y_m$ be any assignment to $\psi$. By the assumtion, there is a set $S\subseteq [m]$ of size $\geq \epsilon m$ such that all constraints $\phi_i, i\in S$ are violated by $u_1,...,u_n$. Then, for any $i\in S$ there is $j\in [q]$ such that $\psi_{i,j}$ is violated.   \qed

\begin{prop}\label{n2} There are constants $d,e, k$ such that for every $W$ there is a $CL^k$-reduction $h:2CSP_W\rightarrow 2CSP_W$ such that $PV_1$ proves that for any $2CSP_W$ instance $\phi$, and any $\epsilon$ $$val(\phi)\leq 1-\epsilon\rightarrow val(h(\phi))\leq 1-\epsilon/(100Wed)$$ and the constraint graph of $h(\phi)$ is $d$-regular.
\end{prop}

\proof By Proposition \ref{cst} and Proposition \ref{expans} there are constants $d,e$ such that for each $e^t, t\in LogLog$, there is a $d$-regular graph $G_{e^t}$ which for any $S\subseteq V, |V|=e^t, |S|\leq e^t/2$ satisfies $|E(S,V-S)|\geq d|S|/4-1/8$. In particular, for each $W$ and $S\subseteq V$, $|S|\leq e^t/2$, we have ($*$): $|E(S,V-S)|\geq |S|/(10W)$. 

The $CL^k$-reduction $h$ works as follows.
\smallskip

Let $\phi$ be a $2CSP_W$ instance. First, erase variables in $\phi$ that do not appear in any constraint. Suppose next that $u_l$ is a variable that appears in $c'\geq 1$ constraints. Put $c:=e^t$ for the smallest natural $t$ such that $c'\leq e^t$. Replace $u_l$ by $c$ variables $y^1_l,...,y^c_l$ so that in each constraint $u_l$ originally appeared in we have different $y^f_l$ (different $c$'s might be needed for each $u_l$). Add a constraint requiring that $y^j_l\leftrightarrow y^{j'}_l$ for every edge $(j,j')$ in the graph $G_c$. Do this for every variable in $\phi$ until each variable appears in $d+1$ constraints, $d$ equality constraints and one original constraint resp. a null constraint that always accepts which is added if necessary. Denote the resulting $2CSP_W$ instance as $\psi$ ($=h(\phi)$). 
\smallskip

If $\phi$ has $m$ constraints, $\psi$ has $\leq m+2dem+2em$ constraints ($m$ original constraints, $\leq 2em$ null constraints and $\leq 2dem$ ``$y^j_l\leftrightarrow y^{j'}_l$'' constraints). If $\phi$ is satisfiable, then so is $\psi$. Suppose that $val(\phi)\leq 1-\epsilon$ and let $y$ be any assignment to $\psi$.
Consider then the plurality assignment $u$ to $\phi$'s variables: $u_i$ gets the most likely value that is claimed for it by $y^1_i,...,y^c_i$. Define $t_i$ to be the number of $y^j_i$'s that disagree with the plurality value of $u_i$.
\smallskip

If $\Sigma^n_{i=1} t_i\geq \epsilon m/2$, then by ($*$) there are $\geq \epsilon m/(20W)$ equality constraints violated in $\psi$.

Suppose that $\Sigma^n_{i=1} t_i<\epsilon m/2$. Since $val(\phi)\leq 1-\epsilon$, there are $\geq \epsilon m$ constraints in $\phi$ violated by $u$. All of these constraints are also present in $\psi$.  If more than $\epsilon m/2$ of them were assigned a different value by $y$ than by $u$, then $\Sigma^n_i t_i\geq \epsilon m/2$. Thus $y$ violates $\geq \epsilon m/2$ constraints in $\psi$. 
\smallskip

Note that all the sets we counted had polynomial size so we had exact counting for them in $PV_1$.
\qed

\begin{prop}\label{n3} There are constants $d,e,k$ such that for any $d',W$ there is a $CL^k$-reduction $h:2CSP_W\rightarrow 2CSP_W$ such that $PV_1$ proves that for any $2CSP_W$ instance $\phi$ with $d'$-regular constraint graph for $d\geq d'$ and for any $\epsilon$, $$val(\phi)\leq 1-\epsilon\rightarrow val(h(\phi))\leq1-\epsilon/(10de)$$ Moreover, the constraint graph $G$ of $h(\phi)$ is $4d$-regular with at least half the edges coming out of each vertex being self-loops and $\lambda(G)\leq 0.9$. \end{prop}

\proof By Proposition \ref{cst} there are constants $d,e$ such that for each $e^t$ where $t\in LogLog$, there is a $d$-regular graph $G_{e^t}$ in $PV_1$ with $\lambda(G_{e^t})\leq 0.1$. The $CL^k$-reduction $h$ works as follows. 

Let $\phi$ be a $2CSP_W$-instance with $n$ variables, $m$ constraints, and $d'$-regular constraint graph $G'$ for $d'\leq d$. Without loss of generality $2m\geq n$. Otherwise, $\phi$ contains variables that are not in any constraint so $d'=0$ and $\phi$ is empty. Add new vertices and self-loops to $G'$ so that it becomes $d$-regular with $e^t$ vertices for the smallest $e^t\geq n$. For each of these new vertices add new variables and for the new self-loops add null constraints that always accept. Then add null constraints for every edge in the graph $G_{e^t}$. Finally, add $2d$ null constraints forming self-loops for each vertex in $G_{e^t}$. 

The resulting instance $\psi$(=$h(\phi)$) has $4d$-regular constraint graph with $\leq 2den$ constraints, and at least half the edges coming out of each vertex being self-loops. Assuming $val(\phi)<1-\epsilon$, there are $\geq\epsilon m\geq\epsilon 2den/(4de)$ violated constraints in $\psi$.

Let $G$ be $\psi$'s constraint graph and $A$ its random-walk matrix. Then $A=3/4B+C/4$ for $C$ the random-walk matrix of $G_{e^t}$ and $B$ the random walk matrix of a $3d$-regular graph. In Section \ref{sec:tens}, we observed that for any $x\in Q^n/m$, $||Ax||\leq 3/4||Bx||+1/4||Cx||+1/L^{1/2}$ and by Proposition \ref{eps}, for any $\delta>0$, $\lambda(B)\leq 1+\delta+1/L$. Thus, assuming $\delta$ is sufficiently small and $L$ sufficiently big, $\lambda(G)\leq 3/4(1+\delta +1/L^{1/2})+1/4\lambda(G_{e^t})+1/L\leq 0.9$.
\qed 

\smallskip

Note that the constant $d$ from Proposition \ref{n3} can be chosen so that it is bigger than the constant $d$ from Proposition \ref{n2}. Therefore, Propositions \ref{n1}, \ref{n2} and \ref{n3} show that there are constants $d,e,k$ such that for any $q$ (and $W=2^q$) there is a $CL^k$-reduction $h:qCSP\rightarrow 2CSP_{2^q}$ such that $PV_1$ proves that $h$ maps any $qCSP$ instance into an instance which is nice with the constraint graph being $d$-regular while the fraction of violated constraints is reduced by a factor at most $1/(1000We^2d^2q)$. This shows that to derive Proposition $\ref{cc1}$ it suffices to prove the following powering proposition:

\begin{prop}\label{pprop} There is $k$ such that for any $W>0$ and sufficiently big square $t\geq 1$ there is an algorithm $A$ with properties described below such that $PV_1$ proves that for any nice $2CSP_W$ instance $\psi$ with $n$ variables with $n\in Log$ the algorithm $A$ produces a $2CSP_{W'}$ instance $\psi^t$ such that:
\begin{enumerate}[label=\rm\arabic*.]

\item $W'\leq W^{d^{5t}}$, where $d$ is the degree of $\psi$'s constraint graph. The instance $\psi^t$ has $\leq d^{5t}n$ constraints.

\item If $\psi$ is satisfiable, then so is $\psi^t$.

\item For every $\epsilon<1/(d\sqrt{t})$, $$val(\psi)\leq 1-\epsilon\rightarrow val(\psi^t)\leq 1-\epsilon\sqrt{t}/(10^6dW^5)$$

\item The formula $\psi^t$ is produced from $\psi$ (by $A$) in time
  $(nd)^kW^{kd^{5t}}$.
\end{enumerate} 
\end{prop}

\proof 

(It might be helpful to the reader to consult the proof we present here in conjunction with the exposition from \cite[Lemma 22.9]{ba} where some concepts are explained with additional details.)
\medskip

Let $\psi$ be a $2CSP_W$ instance with $n$ variables $u_1,...,u_n$ and $m\leq nd/2$ constraints and let $G$ denote the constraint graph of $\psi$.
\smallskip

The formula $\psi^t$ will have $n$ variables $y_1,...,y_n$ over an alphabet of size $W'=W^{d^{5t}}$. A value of a variable $y_i$ is a $d^{5t}$-tuple of values in $\{0,...,W-1\}$ and we will think of it as giving a value $y_i(u_j)$ in $\{0,...,W-1\}$ to every variable $u_j$ in $\psi$ where $j$ can be reached from $i$ using a path of $\leq t+\sqrt{t}$ steps in $G$. Since $G$ is $d$-regular the number of such nodes is $\leq d^{t+\sqrt{t}+1}\leq d^{5t}$.
\smallskip

For every path $p=\left<i_1,...,i_{2t+2}\right>$ in $G$ we will have a constraint $C_p$ in $\psi^t$ depending on variables $y_{i_1}$ and $y_{i_{2t+2}}$ which outputs 0 if and only if there is some $j\in [2t+1]$ such that
\smallskip

\noindent 1. $i_j$ can be reached from $i_1$ using a path of $\leq t+\sqrt{t}$ steps in $G$

\noindent 2. $i_{j+1}$ can be reached from $i_{2t+2}$ using a path of $\leq t+ \sqrt{t}$ steps n $G$

\noindent 3. $y_{i_1}(u_{i_j}), y_{i_{2t+2}}(u_{i_{j+1}})$ violate the constraint in $\psi$ depending on $u_{i_j}$ and $u_{i_{j+1}}$
\medskip

The $2CSP_{W'}$ instance $\psi^t$ can be produced in time $(nd)^kW^{kd^{5t}}$ and has $\leq d^{5t}n$ constraints. Any assignment $u_1,...,u_n$ satisfying $\psi$ induces an assignment $y_1,...,y_n$ satisfying $\psi^t$: each $y_i$ encodes values $u_j$ for $j$'s that can be reached from $i$ by $\leq t+\sqrt{t}$ steps in $G$. Therefore, it remains to show that for $\epsilon<1/(d\sqrt{t})$, $val(\psi)\leq 1-\epsilon\rightarrow val(\psi^t)\leq 1-\epsilon\sqrt{t}/(10^6dW^5)$.
\smallskip

Every assignment $y$ for $\psi^t$ induces the so called plurality assignment $u$ for $\psi$: $u_i$ gets the value $\sigma y(u_i)$ which is the most likely value $y_k(u_i)$ for $y_k$'s where $k$ is obtained by taking a $t$-step random walk from $i$ in $G$. If more than one value is most likely, take the lexicographically first one.
  
Suppose that $val(\psi)\leq 1-\epsilon$, then there is a set $F$ of $\epsilon m$ constraints violated by the plurality assignment.
\smallskip

Pick a random path $p=\left<i_1,...,i_{2t+2}\right>$ in $G$. For $j\in\{1,...,2t+1\}$ we say that the edge $(i_j,i_{j+1})$ in $p$ is truthful if $y_{i_1}(u_{i_j})=\sigma y(u_{i_j})$ and $y_{i_{2t+2}}(u_{i_{j+1}})=\sigma y(u_{i_{j+1}})$. Let $\delta=1/(1000W)$ and denote by $V$ the number of edges in $\left<i_t,...,i_{t+\lfloor\delta\sqrt{t}\rfloor+1}\right>$ that are truthful and in $F$. That is, $V$ is a nonnegative random variable defined on a sample space of size $poly(n)$. If there is at least one such edge, the corresponding constraint in $\psi^t$ is unsatisfied so we want to show that $Pr_p[V>0]\geq \epsilon\sqrt{t}/(10^6dW^5)$.
\smallskip

For each edge $e$ of $G$ and each $j\in \{1,2,...,2t+1\}$, $Pr_p[e=(i_j,i_{j+1})]=1/m$, i.e. each edge has the same probability to be the $j$-th edge in $p$.

\begin{clm}
For any edge $e$ of $G$ and any $j\in\{t,...,t+\lfloor\delta\sqrt{t}\rfloor\}$, $$Pr_p[(i_j,i_{j+1}) \mbox{ is truthful }| \ e=(i_j,i_{j+1})]\geq 1/(2W^2)$$
\end{clm}

\proof
To prove the claim, let $i_1$ be the endpoint of a random walk $p_1$ of length $j$ out of $i_j$ and $i_{2t+2}$ be the endpoint of a random walk $p_2$ of length $2t-j$ out of $i_{j+1}$.  We need to show that $$Pr_{p_1}[y_{i_1}(u_{i_j})=\sigma y(u_{i_{j}})]Pr_{p_2}[y_{i_{2t+2}}(u_{i_{j+1}})=\sigma y(u_{i_{j+1}})]\geq 1/(2W^2)$$ Since half of the edges incident to each vertex are self-loops, we can see an $l$-step random walk from a vertex $i$ as follows: 
\begin{enumerate}[label=\arabic*.]
\item throw $l$ fair coins and let $S_l$ denote the number of ``heads"; 
\item take $S_l$ non-self-loop steps along the graph. 
\end{enumerate}

Denote by $l(p)$ the length of a path $p$ not counting self-loops. Then,
\[\eqalign{
Pr_{p_1}[y_{i_1}(u_{i_j})=\sigma y(u_{i_j})]
&=\Sigma_l Pr[S_j=l]Pr_{p_1}[l(p_1)=l\wedge y_{i_1}(u_{i_j})=\sigma
y(u_{i_j})]\cr
&\geq \Sigma_lPr[S_t=l]Pr_{p_1}[l(p_1)=l\wedge y_{i_1}(u_{i_j})=\sigma
y(u_{i_j})]-20\delta\cr
&\geq 1/W-20\delta
}\]
where the first inequality results from Proposition\ \ref{stat}, while
the last inequality follows from the definition of the plurality assignment which implies that for $j=t$, $Pr_{p_1}[y_{i_1}(u_{i_j})=\sigma y(u_{i_j})]\geq 1/W$. Similarly we obtain
\[Pr_{p_2}[y_{i_{2t+2}}(u_{i_{j+1}})=\sigma y(u_{i_{j+1}})]\geq
(1/W-20\delta).\]
 This proves our claim.\qed

The claim implies $Pr_p[(i_j,i_{j+1}) \mbox{ is truthful and in
}F]\geq |F|/(m2W^2)$ for any $j$ from
$\{t,...,t+\lfloor\delta\sqrt{t}\rfloor\}$. Without a loss of
generality, $|\{t,...,t+\lfloor\delta\sqrt{t}\rfloor\}|$ is
$\lceil\delta\sqrt{t}\rceil$. Thus by linearity of expectation, 
\[E[V]\geq \epsilon\lceil\delta\sqrt{t}\rceil /(2W^2)\]
By Proposition \ref{stat} 2., $Pr[V>0]\geq E[V]^2/E[V^2]$, so to conclude the proof it suffices to show that $E[V^2]\leq 50d\epsilon\lceil\delta\sqrt{t}\rceil$.

Denote by $V'$ the number of edges in $\left<i_t,...,i_{t+\lfloor\delta\sqrt{t}\rfloor+1}\right>$ that are in $F$. For any $j$ from $\{t,...,t+\lfloor\delta\sqrt{t}\rfloor\}$ put $I_j:=1$ iff $(i_j,i_{j+1})\in F$. Further, let $S$ be the set of vertices contained in an edge from $F$. Then, assuming that the constant $L$ from our definition of $\lambda(G)$ satisfies $L>d$ and $L>\delta\sqrt{t}$,
\[\eqalign{
  E[V^2]
&\leq E[V'^2]\cr
&=E[\Sigma_{j,j'}I_jI_{j'}]\cr
&=E[\Sigma_j I^2_j]+E[\Sigma_{j\neq j'}I_jI_{j'}]\cr
&=\epsilon\lceil\delta\sqrt{t}\rceil+2\Sigma_{j<j'}Pr_p[(i_j,i_{j+1})\in
F\wedge (i_{j'},i_{j'+1})\in F] \cr
 &\leq \epsilon\lceil\delta\sqrt{t}\rceil+2\Sigma_{j<j'}Pr_{(i_j,i_{j'})\in G^{j'-j}}[i_j\in S\wedge i_{j'}\in S]\cr
&\leq  \epsilon\lceil\delta\sqrt{t}\rceil+2\Sigma_{j<j'}\epsilon
d(\epsilon d+2\cdot 0.9^{j'-j}) \ \ \ \ \ \ \ \ \ \ \ \ \ \ \ \ \ \ \
\ \ \ \ \  \mbox{by\ Proposition\ \ref{ecc}}\cr
&\leq \epsilon\lceil\delta\sqrt{t}\rceil+2\epsilon^2d^2\lceil\delta
\sqrt{t}\rceil^2+40\epsilon d\lceil\delta\sqrt{t}\rceil\leq 50\epsilon
d\lceil\delta\sqrt{t}\rceil \ \ \ \ \ \ \mbox{using \ }\epsilon<
1/(d\sqrt{t})\rlap{\hbox to 24 pt{\hfill\qEd}}
}\]\smallskip

This concludes our formalization of the PCP theorem in the theory $PV_1$. It can be briefly summarized as follows. In Theorem \ref{pcp} we formulated the PCP theorem as a $\forall\Sigma^b_1$-formula. Thus, by $\forall\Sigma^b_1$-conservativity of $S^1_2$ over $PV_1$ we could afford to work instead in the theory $S^1_2$. Specifically, we used $\Pi^b_1$-LLIND induction available in $S^1_2$ to show that the PCP theorem is a consequence of a statement about CSP instances, Proposition \ref{cc}. Then we observed that the CSP formulation of the PCP theorem is a collorary of two propositions, Gap amplification \ref{cc1} and Alphabet reduction \ref{cc2}. The latter one was an application of the exponential PCP theorem in a scaled-down setting where we needed to count only sets of constant size, hence it was provable already in $PV_1$. The gap amplification was a consequence of a CL-reduction into nice CSP instances and Powering proposition \ref{pprop}. The reduction to nice instances used the $(n,d,\lambda)$-graphs which we constructed in Section \ref{rest}.  Section \ref{rest} contained the most challenging part where we needed to employ certain approximating tools to reason about algebraic definitions of pseudorandom constructions in $PV_1$. In the remaining part of the proof of the PCP theorem, including the powering proposition, we were mainly verifying step by step that the reasoning used in the standard proof does not exceed the possibilities of the theory $PV_1$.

\section{Acknowledgement}

I would like to thank Jan Kraj\'{i}\v{c}ek for many constructive discussions during the development of the paper and Sam Buss for detailed comments and suggestions which improved the quality of the manuscript. I would also like to thank Neil Thapen, Pavel Pudl\'{a}k and Emil Je\v{r}\'{a}bek for comments and suggestions during its seminar presentation. This research was supported by grants GA UK 5732/2014 and SVV-2014-260107.

\def\conver{
\medskip

We will also need a converse of Proposition \ref{expans}.

\begin{prop}(in $PV_1$)\label{expans} Let $G$ be a $d$-regular graph with $n$ vertices $V$ and edges $E$ such that for every $S\subseteq V, |S|\leq n/2$, $$|E(S,V-S)|\geq d\rho |S|$$ where $E(S,T)$ denotes the set of edges $(i,j)\in E$ with $i\in S, j\in T$. Then $\lambda(G)\leq 1-\rho^2/2$
\end{prop}

\proof Let $A$ be the random-walk matrix of $G$. (If $\lambda(G)=\lambda$, there is) Let $x\in Q^n$ be such that$\left<x,\mathbf{1}\right>=0$ (such that $||Ax||=\lambda ||x||$). 

Let $x=y+z$ where $y$ is equal to $x$ on the coordinates on which $x$ is positive and equal to 0 otherwise. Assume $y$ is nonzero on at most $n/2$ of its coordinates. Othersise take $-x$ instead of $x$. Let $Z=\Sigma_{i,j}A_{i,j}|y_i^2-y_j^2|$. It suffices to show the following two claims: 
\medskip

Claim 1: $Z\geq 2\rho ||y||^2$ 

Claim 2: $Z\leq (8(1-||Ax||))^{1/2}||y||^2$
\medskip

The proof of Claim 1: $x_i=|V-S|$ if $i\in S$ and $x_i=-|S|$ if $i\in V-S$. Put $Z:=\Sigma_{i,j} A_{i,j}(x_i-x_j)^2$ for the random-walk matrix $A$ of $G$. Then $Z=\frac{2}{d}|E(S,V-S)|(|S|+|V-S|)^2$. As $A$'s rows and columns sum up to one, we have also $$Z=\Sigma_{i,j} A_{i,j}x^2_i-2\Sigma_{i,j}A_{i,j}x_ix_j+\Sigma_{i,j} A_{i,j}x^2_j=2||x||^2_2-2\left<x,Ax\right>$$ Note that $\left<y,z\right>=0$ and so $\left<Ay,y\right>+\left<Az,y\right>+\left<Ay,z\right>=\left<Ax,x\right>-\left<z,z\right>=
\left<Ax,Ax\right>-||z||^2=\lambda^2||x||^2-||z||^2=\lambda^2||y||^2-(1-\lambda^2)||z||^2$

so $||Ax||_2= ||A\frac{x}{||x||_2}||_2||x||_2\leq \lambda ||x||_2$. By Cauchy-Schwarz inequality, $\left<x,Ax\right>\leq ||x||_2||Ax||_2$. Therefore, $$\frac{1}{d}|E(S,V-S)|(|S|+|V-S|)^2\geq (1-\lambda)||x||_2$$ It remains to observe that $||x||_2=|S||V-S|(|S|+|V-S|)$  \qed
}

\end{document}